\documentstyle[10pt]{article}

\title{Classification of Direct Limits of
Even Cuntz-Circle Algebras}
\author{ Huaxin Lin\\
 and \\
N. Christopher Phillips}
\date{}
\begin{document}
\maketitle

\newcommand{\beq}{\begin{equation}}
\newcommand{\eeq}{\end{equation}}
\newcommand{\bd}{\begin{displaymath}}
\newcommand{\ed}{\end{displaymath}}
\newcommand{\ben}{\begin{enumerate}}
\newcommand{\een}{\end{enumerate}}
\newcommand{\bde}{\begin{description}} 
\newcommand{\ede}{\end{description}}
\newcommand{\beqr}{\begin{eqnarray*}}
\newcommand{\eeqr}{\end{eqnarray*}}
\newcommand{\hs}[1]{\hspace{#1}}
\newcommand{\vs}[1]{\vspace{#1}}
\newcommand{\bc}{\begin{center}}
\newcommand{\ec}{\end{center}}
\newcommand{\bv}{\begin{verbatim}}
\newcommand{\ev}{\end{verbatim}}

\newcommand{\af}{\alpha}

\newcommand{\half}{\frac{1}{2}}
\newcommand{\p}{\partial}
\newcommand{\bt}{\beta}
\newcommand{\gm}{\gamma}
\newcommand{\dt}{\delta}
\newcommand{\ep}{\epsilon}
\newcommand{\zt}{\zeta}
\newcommand{\et}{\eta}
\newcommand{\tta}{\theta}
\newcommand{\ld}{\lambda}
\newcommand{\sm}{\sigma}
\newcommand{\ssm}{\sigma^2}
\newcommand{\vph}{\varphi}
\renewcommand{\i}{\subset}

\newcommand{\s}[1]{\noindent {\bf EXERCISE~#1}{}}

\setcounter{section}{-1}
\newtheorem{thm}{Theorem}[section]
\newtheorem{lem}{Lemma}[section]
\newtheorem{prop}{Proposition}[section]
\newtheorem{dfn}{Definition}[section]
\newtheorem{conj}{Conjecture}

\newcommand{\qed}{\hspace*{\fill}Q.E.D.}  
\newcommand{\cd}{\cdots}
\newcommand{\lr}{\longrightarrow}

\pagenumbering{arabic}
\parindent=0pt

\def\S{\hspace{6.4mm}}


\newcommand{\OA}[1]{{\cal O}_{#1}}
\newcommand{\QA}[1]{{\cal Q}_{#1}}
\newcommand{\SO}[1]{C(S^1) \otimes \OA{#1}}
\newcommand{\sO}[1]{S \otimes \OA{#1}}
\newcommand{\So}[1]{C_0 (S^1 \setminus \{1\}) \otimes \OA{#1}}
\newcommand{\KSO}[2]{K_{#1}(C(S^1) \otimes \OA{#2})}
\newcommand{\KsO}[2]{K_{#1}(S \otimes \OA{#2})}
\newcommand{\KO}[2]{K_{#1} (\OA{#2})}
\newcommand{\Ext}{{\rm Ext}}
\newcommand{\Hom}{{\rm Hom}}
\newcommand{\Tor}{{\rm Tor}}
\newcommand{\ext}[2]{\Ext^1_{{\bf Z}}({#1},{#2})}
\newcommand{\Class}{\mbox{\boldmath ${\cal C}$}}


\newcommand{\CA}{$C^*$-algebra}
\newcommand{\CSalg}{Cuntz-circle algebra}
\newcommand{\DECSalg}{direct limit of even Cuntz-circle algebras}
\newcommand{\aab}{approximately absorbing}
\newcommand{\aue}{approximate unitary equivalence}
\newcommand{\ayue}{approximately unitarily equivalent}
\newcommand{\mops}{mutually orthogonal projections}
\newcommand{\hm}{homomorphism}
\newcommand{\pisca}{purely infinite simple \CA}
\newcommand{\andeqn}{\,\,\,\,\,\, {\rm and} \,\,\,\,\,\,}
\newcommand{\QED}{\rule{1.5mm}{3mm}}


\newcommand{\arrow}{\rightarrow}
\newcommand{\tdsum}{\widetilde{\oplus}}
\newcommand{\pa}{\|}  
\renewcommand{\ep}{\varepsilon}
\newcommand{\id}{{\rm id}}
\newcommand{\aueeps}[1]{\stackrel{#1}{\sim}}
\newcommand{\aeps}[1]{\stackrel{#1}{\approx}}
\newcommand{\dirlim}{\displaystyle \lim_{\longrightarrow}}

\pagebreak


\pagebreak

\vspace{0.5 in}

\begin{center} {\LARGE Contents}  \end{center}
\vspace{1 in}

\S {\bf Abstract} \hfill 4

\S {\bf Introduction}                           \hfill     5

{\bf 1. Approximately absorbing homomorphisms}    \hfill     8

{\bf 2. Homotopies of asymptotic morphisms} \hfill 14

{\bf 3. Approximate unitary equivalence of homomorphisms} \hfill 37

{\bf 4. The existence theorem } \hfill 48

{\bf 5. The main results} \hfill 56

\S {\bf References }                              \hfill     69

\pagebreak

\vspace{0.1 in}

\begin{center}{\LARGE Abstract} \end{center}

\vspace{0.5 in}

	We prove a classification theorem for purely infinite simple
\CA s that is strong enough to show that the tensor products of two
different irrational rotation algebras with the same even Cuntz algebra
are isomorphic. In more detail, let $\Class$ be
the class of simple \CA s $A$ which are direct limits
$A\cong \dirlim A_k,$ in which each $A_k$ is a finite direct sum
of algebras of the form
$C(X)\otimes M_n\otimes \OA{m},$ where $m$ is even,
$\OA{m}$ is
the Cuntz algebra,
and $X$ is either a point, a compact interval, or the
circle $S^1,$ and each map $A_k\to A$ is approximately
absorbing. (``Approximately absorbing'' is defined in Section 1.)
We show that two unital \CA s $A$ and $B$ in $\Class$
are isomorphic if and only if
$$
(K_0(A),[1_A], K_1(A))\cong (K_0(B),[1_B],K_1(B)).
$$
This class is large enough to exhaust all possible $K$-groups: if
$G_0$ and $G_1$ are countable odd torsion (abelian) groups and
$g\in G_0,$ then there is a \CA\, $A$ in $\Class$ with
$(K_0(A),[1_A],K_1(A))\cong (G_0,g,G_1).$ The class $\Class$ contains
the tensor products of irrational rotation algebras with even Cuntz
algebras. It is also closed under the formation of
hereditary subalgebras, countable direct limits (provided that
the direct limit is simple), and tensor products with simple
AF algebras.

\vspace{0.5in}

\pagebreak

\pagestyle{myheadings}
\markboth{\hfill {\rm Huaxin Lin and N. Christopher Phillips}\hfill}{\hfill
{\rm Direct limits of Cuntz-circle algebras} \hfill}

\section{Introduction}

\vspace{\baselineskip}

We prove a classification theorem for simple direct limits of what
we call even \CSalg s: finite direct sums of algebras of the form
$C(X) \otimes M_n \otimes \OA{m}$, where $m$ is even, $\OA{m}$ is
the Cuntz algebra (first introduced in \cite{Cu2}),
and $X$ is either a point, a compact interval, or the
circle $S^1$. The unital version of our main theorem is:

\vspace{0.6\baselineskip}

{\bf Theorem A.} (Theorem 5.4)
 Let $A = \dirlim A_k$ and $B = \dirlim B_k$ be
simple separable unital \CA s, which are direct limits of even
\CSalg s. Assume that the \hm s $A_k \to A$ and $B_k \to B$ are
``approximately absorbing'' (defined in Section 1). Then
$A \cong B$ if and only if
 $( K_0 (A), [1_A], K_1 (A)) \cong ( K_0 (B), [1_B], K_1 (B))$.
In particular, if there are isomorphisms
$\alpha_0 : K_0 (A) \to K_0 (B)$ and $\alpha_1 : K_1 (A) \to K_1 (B)$
such that $\alpha_0 ([1_A]) = [1_B]$, then $A$ is isomorphic to $B$.

\vspace{0.6\baselineskip}

As a corollary, we obtain:

\vspace{0.6\baselineskip}

{\bf Theorem B.} (Corollary 5.12) Let $\theta_1$ and $\theta_2$ be
irrational numbers, and let $A_{\theta_1}$ and $A_{\theta_2}$ be
the corresponding irrational rotation algebras. Then for any even
$m$, we have
$A_{\theta_1} \otimes \OA{m} \cong A_{\theta_2} \otimes \OA{m}$.

\vspace{0.6\baselineskip}

This contrasts with the fact, due to Rieffel \cite{Rf} and
Pimsner and Voi\-cu\-les\-cu \cite{PV}
that
$A_{\theta_1} \cong A_{\theta_2}$ only when
  $\theta_1 = \pm \theta_2 \pmod{\bf Z}.$ (Theorem B was already
known for $m = 2$ \cite{Ln5}, and remains unknown for odd $m$.)
generally, tensor products of simple direct limits of circle
or interval algebras (whether of real rank 0 or 1) with even Cuntz
algebras are classified up to
isomorphism by their $K$-theory.
Further results along these lines are given in Section 5.

We do not know if maps from even \CSalg s to simple direct limits of
even \CSalg s are necessarily \aab . Therefore, we do not classify
arbitrary simple direct limits of even \CSalg s. However, there are
more general situations than Theorem B in which the \aab\  condition
is automatic. For example, in Section 5 we will define ``co-Cuntz
algebras'' $\QA{m}$, which are unital \pisca s satisfying
$K_0 (\QA{m}) = 0$ and $K_1 (\QA{m}) \cong {\bf Z}/(m - 1) {\bf Z}$.
We prove, without any hypotheses on the maps of the systems:

\vspace{0.6\baselineskip}

{\bf Theorem C.} (Theorem 5.24 (1))
Let $A = \dirlim A_k$ and $B = \dirlim B_k$ be simple unital direct
limits, in which the $A_k$ and $B_k$ are each finite direct sums
of matrix algebras over even Cuntz algebras and even co-Cuntz algebras.
Then $A \cong B$ if and only if
 $( K_0 (A), [1_A], K_1 (A)) \cong ( K_0 (B), [1_B], K_1 (B))$.

\vspace{0.6\baselineskip}

We also prove that our class is closed under the formation of
direct limits.

All feasible values of the invariant are actually realized:

\vspace{0.6\baselineskip}

{\bf Theorem D.} (Theorem 5.26 (1))  Let $G_0$ and $G_1$ be countable
odd torsion groups, and let $g_0 \in G_0$. Then there exists a unital
\CA , in the classes covered by Theorems A and C, such that
 $( K_0 (A), [1_A], K_1 (A)) \cong ( G_0 , g_0, G_1 )$.

\vspace{0.6\baselineskip}

Theorem A probably remains true if the condition on $X$ is relaxed
slightly, to allow arbitrary compact subsets of $S^1$ in the
definition of a \CSalg . Proving this would make an already long paper
even longer, and we would get no new values of the invariant.
So we don't do it.

Our results are part of the general classification program for simple
separable nuclear \CA s. This program was initiated by George Elliott
\cite{Ell1} many years ago (1976) with the classification of
AF algebras. Starting much more recently (about 1990), it has been
extended to the class of \CA s of real rank zero which are direct limits
of circle algebras by Elliott \cite{Ell2}, and to much larger classes
of (simple) stably finite direct limit algebras by Elliott
\cite{Ell3}, \cite{Ell4}, Su \cite{Su1}, \cite{Su2}, \cite{Su3},
Elliott and Gong \cite{EG1}, \cite{EG2},
  Elliott, Gong, Lin and Pasnicu \cite{EGLP}, etc.,
and to classes of \pisca s by Bratteli, Kishimoto, R\o rdam and
St\o rmer \cite{BKRS} and by R\o rdam  \cite{Rr1}, \cite{Rr2},
\cite{Rr3}. In addition to the classification of direct limits, one
of the striking successes of this program is the Elliott-Evans
 realization of the irrational rotation algebras as direct limits
of circle algebras \cite{EE}; we use this result to derive Theorem
B from Theorem A.

The classification program is the analog for
\CA s of the classification of hyperfinite factors
in the theory of von Neumann algebras. It has only recently become
apparent that some reasonable class of simple \CA s is in fact
classifiable. Classification of separable commutative \CA s
is equivalent to classification of second countable locally compact
Hausdorff spaces, a problem long considered out of reach, while
classification of commutative weak* separable von Neumann algebras
is fairly easy. Similarly, even ignoring the problem of classifying
spaces, type I \CA s are very much harder to classify than type I
von Neumann algebras. (The extension problem remains a major obstacle,
even in the presence of $KK$-theory.) These differences between the
\CA\  and von Neumann algebra situations have led to
the assumption that simple \CA s should be much
harder to classify than factors. (For example, in \cite{Ph2}, the
second author asked for a proof that
$A_{\theta_1} \otimes \OA{m} \not\cong A_{\theta_2} \otimes \OA{m}$
 for most values of $\theta_1$ and $\theta_2$.) The results of the
classification program described above, including those of this
paper, suggest that simple \CA s are in fact classifiable, at least
in the ``amenable'' case (or perhaps a large subset of it).

Independently of this work, R\o rdam has in Section 5 of \cite{Rr3}
introduced a ``classifiable class''
of separable unital \pisca s.
(The main object of R\o rdam's paper is the classification of certain
\pisca s $A$ for which $ K_1 (A)$ is torsion-free; the $K_1$-groups
of the algebras in this paper are all odd torsion groups.)
Algebras $A$ in R\o rdam's class are determined up to isomorphism by the
invariant
 $( K_0 (A), [1_A], K_1 (A))$. We certainly believe that our algebras
are in fact in his class.
Unfortunately, neither R\o rdam's methods nor ours seem to prove this.
In particular, his results do not seem to
help with the proof of our Theorem B, the
isomorphism of tensor products of different irrational rotation
algebras with the same even Cuntz algebra, or other similar results
in our Section 5.

Bratteli, Elliott, Evans, and Kishimoto are presently working on
a classification theorem presumably covering
a larger class of simple \CA s.

Our proof follows the standard outline, first introduced by Elliott
in \cite{Ell2}. (We use $KK$-classes in place of \hm s on $K$-theory,
as in \cite{Rr2}.) Thus, we construct approximate intertwinings
of two given direct systems, and for this we need an existence
theorem and a uniqueness theorem. The first three sections of this
paper contain the uniqueness theorem, the fourth contains the existence
theorem, and in the last section we put everything together.
Our uniqueness theorem would obviously have been impossible without
R\o rdam's results on even Cuntz algebras \cite{Rr1}. In fact, we
need his results not only in Section 3, but also for the preliminary
work in Sections 1 and 2, and even for the existence theorem in
Section 4. However, no proofs in this paper closely resemble the proof
of the main technical result in \cite{Rr1}.

In more detail, the outline of this paper is as follows. In Section
1 we define \aab\ \hm s, and give some conditions under which
\hm s from $\SO{m}$ to a \pisca\  are automatically \aab . It turns
out that we need the main technical theorem of \cite{Rr1} to
prove that anything at all is \aab . In Section 2 we prove a weak
version, sufficient for our purposes, of the following statement:
If $m$ is even, then two \hm s from $\SO{m}$ to a \pisca\  $A$,
with the same class in $KK^0 (\SO{m}, A)$, can be connected by a
``discrete homotopy'' of asymptotic morphisms. This section is the
longest in the paper; it also requires \cite{Rr1}. Section 3
contains the proof of the uniqueness theorem. It is based on the
absorption argument first introduced by the second author in
\cite{Ph1}, and requires results of the first author \cite{Ln5}
on approximately commuting unitaries to cope with the fact
that Section 2 yields only asymptotic morphisms.
It is this absorption argument that requires that \hm s be assumed
\aab . In Section 4 we prove the existence theorem: every
class in $KK^0 (\SO{m}, \SO{n})$, with $m$, $n$ even, is represented
by a \hm , unital if the $K$-theory data allows it. We use
the Universal Coefficient Theorem to combine various \hm s already
constructed in the literature. Here, too, we need to restrict to
even Cuntz algebras, since the construction, by Loring \cite{Lr3},
of one of the needed \hm s requires R\o rdam's classification
theorem in \cite{Rr1}. Finally, in Section 5 we put the existence and
uniqueness theorems together to prove our main theorem and derive
various corollaries. We also prove several auxiliary results, such
as Theorem C, and construct the co-Cuntz algebras $\QA{m}$.

We state here some terminology and notation that will be used
throughout this paper.

\vspace{0.6\baselineskip}

{\bf 0.1 Definition} A {\em \CSalg} is a \CA\  $A$ which has the form
\[ A \cong
 \bigoplus_{i = 1}^k M_{r(i)} \otimes C(X_i) \otimes \OA{m(i)} , \]
with each $X_i$ being a connected compact subset of
the circle $S^1$.
In this expression, $k$ and the  $m(i)$ are finite, and $m(i) \geq 2$.
An {\em even} \CSalg\  is one for which all the $m(i)$ above are even.

\vspace{0.6\baselineskip}

Cuntz-circle algebras are thus finite direct sums of matrix algebras
over $C(X) \otimes \OA{m}$, where $X$ is allowed to be homeomorphic to
a circle, a point, or a compact interval.

\vspace{0.6\baselineskip}

{\bf 0.2 Conventions}
(1) Throughout this paper, $u$ is the canonical
generating unitary in $C(S^1)$, or in $C(X)$ for any subset
$X \subset S^1$. Also, $s_1, \dots, s_m$ are the canonical generating
isometries of the Cuntz algebra $\OA{m}$; they thus satisfy
\[
s_j^* s_j = 1 \andeqn \sum_{j = 1}^m s_j s_j^* = 1.
\]

(2) Let $B$ be a \CA. We write
$B^+$ for the unitization of $B$ (the unit being added regardless of
whether or not $B$ already has one). We write $\tilde B$ for the
algebra which is equal to $B$ if $B$ is unital, and equal to
$B^+$ if $B$ is not unital.

(3) If $B$ is a unital \CA, then $U(B)$ is the unitary group of $B$ and
$U_0(B)$ is the connected component of $U(B)$ containing the identity
of $B.$

\vspace{0.6\baselineskip}

This work was done when the first author was in SUNY at
Buffalo. He would like to thank L. Coburn, J. Kraus, T. Natsume,
K. Olsen and J. Xia for their hospitality during his stay.
The first author would also like to thank George Elliott and Mikael
R\o rdam for useful conversations, and
 second author would like to thank Marius D\v{a}d\v{a}rlat and
Terry Loring. The second author would also like to thank SUNY
Buffalo for its hospitality during a visit in March 1993.

\vfill
\newpage

\section{ Approximately absorbing homomorphisms}

\vspace{\baselineskip}

In this section, we introduce \aab\ \hm s and give some of their
elementary properties. We have not proved that \hm s from
$\SO{m}$ to a \pisca\  $B$ are automatically \aab, so this condition
appears as a hypothesis in our most general theorems. However, we will
see in this section that
the \hm s in the direct systems corresponding to the most interesting
cases (tensor products of even Cuntz algebras with irrational rotation
algebras, Bunce-Deddens algebras, etc.) are automatically \aab .

R\o rdam's work (\cite{Rr1} and \cite{Rr2}) is already needed to
prove that
a \hm\  from $\OA{m}$ to $B$ is \aab. We will therefore need to assume
throughout this section that our Cuntz algebras are even.

We begin by establishing terminology and notation for \aue.

\vspace{0.6\baselineskip}

{\bf 1.1 Definition}. Let $A$ and $B$ be \CA s, let $G$ be a set of
generators of $A$, and let $\varphi$ and $\psi$ be two \hm s from $A$ to
$B$. We say that $\varphi$ and $\psi$ are {\em \ayue\  to within} $\ep$,
with respect to $G$, if there is a unitary $v \in \tilde{B}$ such that
\[
\| \varphi(g) - v \psi(g) v^* \| < \ep
\]
for all $g \in G$. We abbreviate this as
\[
\varphi \aueeps{\ep} \psi.
\]
(Note that we have suppressed $G$ in the notation.)
We say that $\varphi$ and $\psi$ are {\em \ayue} if
 $\varphi \aueeps{\ep} \psi$ for all $\ep > 0$.
(Of course, this notion does not depend on the choice of $G$.)

\vspace{0.6\baselineskip}

{\bf 1.2 Convention}
In the previous definition, if $A = \SO{m}$ (or $C(X) \otimes \OA{m}$
with $X \subset S^1$), then we will take
the generating set $G$ to be
\[
\{u \otimes 1\} \cup  \{ 1 \otimes s_j : j = 1, \dots, m\},
\]
unless otherwise specified.

\vspace{0.6\baselineskip}

{\bf 1.3 Definition} Let $A$ be any unital \CA ,
 and let $B$ be a \pisca .
 Let $\varphi, \psi: A \to B$ be two homomorphisms, and assume that
$\varphi(1) \neq 0$ and
$[\psi(1)] = 0$ in $K_0 (B)$. We define a homomorphism
$\varphi \tdsum \psi : A \to B$, well defined up to unitary equivalence,
by the following construction. Choose a projection $q \in B$ such that
$0 < q < \varphi(1)$ and $[q] = 0$. Since $B$ is purely infinite and
simple, there are partial isometries $v$ and $w$ such that
$  vv^* = \varphi(1) -q$, $v^*v = \varphi(1)$, $ww^* = q$,
 and $w^*w = \psi(1)$.
Now define $(\varphi \tdsum \psi)(a) = v\varphi(a)v^* + w\psi(a)w^*$
for $ a \in A$.

\vspace{0.6\baselineskip}

{\bf 1.4  Definition} Let $X$ be a compact space, let $B$ be a \pisca ,
 and let
 $\varphi: C(X)\otimes M_k(\OA{m})\to B$ be a homomorphism.
Let $X_0 \subset X$ be the closed set such that
$\ker(\varphi) =  C_0(X\setminus X_0)\otimes M_k(\OA{m})$.
Then $\varphi$ is {\em \aab} if $\varphi$ is \ayue\ to
$\varphi \tdsum \psi$
 for any \hm\ $\psi : C(X) \otimes M_k (\OA{m}) \to B$
of the form $\psi (f \otimes a) = \sum_{i=1}^{l} f(x_i) \psi_i (a)$,
where:

(1) $x_1, \dots, x_l \in X_0.$

(2) $\psi_1, \dots, \psi_l : M_k (\OA{m}) \to B$ are \hm s.

(3) $\psi_1 (1), \dots, \psi_l (1)$ are \mops\ in $B$.

(4) $[\psi_1] =  \cdots = [\psi_l] = 0 $ in $KK^0 (\OA{m}, B)$.

Note that, if we regard $C(X)\otimes M_k(\OA{m})$ as the
\CA\  of $M_k(\OA{m})$-valued continuous functions on
$X,$ then, for $f \in C(X) \otimes M_k (\OA{m}),$ we may also write
$\psi(f) = \sum_{i=1}^l \psi_i (f(x_i))$.

If $A$ is a \CSalg , then we say that $\varphi : A \to B$ is \aab\  if
the restriction of $\varphi$ to each summand of $A$ is \aab .

\vspace{0.6\baselineskip}

The next lemmas will be about \aab\  \hm s from
$C(X) \otimes M_k (\OA{m})$, but they have trivial generalizations to
\hm s from \CSalg s.

\vspace{0.6\baselineskip}

{\bf 1.5 Remark} It is clear from the definition that
$\varphi: C(X) \otimes M_k (\OA{m}) \to B$ is \aab\  if and only if
$\varphi$ is \aab\  when regarded as a \hm\  from
$C(X) \otimes M_k (\OA{m})$ to $\varphi(1) B \varphi(1)$.

\vspace{0.6\baselineskip}

The following lemma asserts that homomorphisms are \aab\ when $X$ is a
point and $m$ is even.

\vspace{0.6\baselineskip}

{\bf 1.6 Lemma} Let $m$ be even, let $B$ be a \pisca , and
 let $\varphi,\, \psi: M_k(\OA{m})\to B$ be two homomorphisms.
If $\varphi \neq 0$ and
$[\psi]=0$ in $KK(M_k(\OA{m}),B)$, then $\varphi$ is approximately
unitarily equivalent to $\varphi \tdsum \psi$.

\vspace{0.6\baselineskip}

{\it Proof:} Since $B$ is purely infinite, $\psi(1)$ is equivalent
to a subprojection of $\varphi(1)$. We may therefore
assume $\psi(1) \leq \varphi(1)$.
 Replacing $B$ by $\varphi(1)B\varphi(1)$,
we may assume that $B$ and $\varphi$ are unital.
Since $B$ is purely infinite and $M_k(\OA{m} )$ is a unital
direct limit of even Cuntz algebras,
the result now follows from Theorem 5.3 of \cite{Rr2}.  \QED

\vspace{0.6\baselineskip}

The next lemma provides a more convenient way to show that a \hm\  is
\aab .

\vspace{0.6\baselineskip}

{\bf 1.7 Lemma} Let the notation be  as in Definition 1.4,
with $X \i S^1$ and $m$ even.
Then  $\varphi$ is    \aab\ if and only if for every $\ep > 0$
and $\lambda_1, \dots, \lambda_l \in X_0$, there exist a unital
\hm\  $\sigma: M_k (\OA{m}) \to \varphi(1) B \varphi(1)$,
a unitary $v \in \varphi(1) B \varphi(1)$,
and nonzero mutually orthogonal projections
$q_1, \dots, q_l \in \varphi(1) B \varphi(1)$ such that:

(1)  For $j = 1, \dots, m$, and for
 each standard matrix unit $e \in M_k$,
\[ \pa \sigma(e \otimes s_j) - \varphi(1 \otimes e \otimes s_j) \pa
   < \ep. \]

(2) $\| v - \varphi(u \otimes 1) \| < \ep$.

(3) $v$ commutes with the range of $\sigma$ and with the $q_i$, and
each $q_i$ commutes with the range of $\sigma$.

(4) $q_i v q_i = \lambda_i q_i$ for $i = 1, \dots, l$.

\vspace{0.6\baselineskip}

{\em Proof:} Note that the conditions on $\sigma$ and $v$ say that
the map $f \otimes a \mapsto f(v)\sigma(a)$ defines a homomorphism
from $C(X) \otimes M_k(\OA{m})$ to $B$ which agrees with $\varphi$
to within $\ep$ on a particular set $G$ of generators.

Replacing $B$ by $\varphi(1) B \varphi(1)$, we may assume that $B$
and $\varphi$ are unital.

Suppose $\varphi$ is \aab .
  Choose nonzero mutually orthogonal projections
$p_1, \dots, p_l \in B$ whose classes in $K_0(B)$ are all zero. With the
help of $km$ \mops\ summing to $p_i$ whose $K_0$-classes are zero,
it is easy to construct a homomorphism $\psi_i : M_k(\OA{m}) \to B$
such that $\psi_i (1) = p_i$.
Define $\psi : C(X) \otimes M_k(\OA{m}) \to B$
by $\psi(f \otimes a) = \sum_i f(\lambda_i)\psi_i(a)$.  Since $\varphi$
is \aab , $\varphi$ is \ayue\ to
 $\varphi \tdsum \psi$. Let $w$ implement
this \aue\ to within $\ep$ on the generators listed in the statement
of the lemma, and set $v = w(\varphi \tdsum \psi)(u)w^*$
 and $q_i = w p_i w^*$.

Conversely, let the conditions of the lemma hold for $\varphi$, and let
$\psi$ be as in the definition of \aab . We follow the notation of that
definition, except that we call the points $\lambda_i$ instead of $x_i$.
Let $\ep > 0$. We want to show that $\varphi \aueeps{\ep}
 \varphi \tdsum \psi$ (with respect to the set  $G$
 of generators above).
Choose a homomorphism $\sigma$, a unitary $v$,
and projections $q_i$, as in the hypotheses, except using $\ep/3$
in place of $\ep$.
Define
$\varphi^{\prime}: C(X) \otimes M_k (\OA{m}) \to B$
by $\varphi^{\prime} (f \otimes a) = f(v)\sigma(a)$. Then
$\varphi^{\prime} \aueeps{\ep/3} \varphi,$  whence also
$\varphi^{\prime} \tdsum \psi \aueeps{\ep/3}
 \varphi \tdsum \psi.$
It therefore suffices to show that
$\varphi^{\prime} \tdsum \psi \aueeps{\ep/3} \varphi^{\prime}.$

To see this,
let $q_0 = \varphi^{\prime}(1) - q_1 - \cdots -q_l$,
define $\varphi_i^{\prime}(f \otimes a)= f(\lambda_i)q_i \sigma(a)$,
 and observe that
 $\varphi^{\prime}(b) = \sum_{i= 0}^l \varphi_i^{\prime}(b) $ .
Further define
 $\overline{\psi}_i (f \otimes a) = f(\lambda_i)\psi_i(a)$,
so that $\psi (b) = \sum_{i=1}^l \overline{\psi}_i(b)$.
It follows from the previous lemma that $ \varphi^{\prime}_i$
is \ayue\ to $\varphi^{\prime}_i \tdsum \overline{\psi}_i$.
Forming the (orthogonal) sum over $i = 1, \dots, l$, and adding
$\varphi_0^{\prime}$, gives the required \aue . \QED

\vspace{0.6\baselineskip}

{\bf 1.8 Corollary} Let $m$ be even.
Let $B$ and $C$ be purely infinite simple \CA s,
and let $\eta : B \to C$ be a nonzero homomorphism. If $X \subset S^1$
and
$\varphi: C(X) \otimes M_k(\OA{m}) \to B$ is \aab , then
$\eta \circ \varphi$ is also \aab .

\vspace{0.6\baselineskip}

{\em Proof:} Since $B$ is simple, $\eta$ is injective, and the condition
of the lemma is preserved under application of $\eta$. \QED

\vspace{0.6\baselineskip}

{\bf 1.9 Corollary}
Let $A$ be a simple separable \CA ,  obtained as a simple direct limit
 with no dimension growth in the sense of \cite{Ph4}.
 Let $m$ be even,
let $B$ be  a purely infinite simple \CA , and let $X \subset S^1$
 be compact.
Let $\iota: A \otimes \OA{m} \to B$ be an injective \hm ,
and let
$\mu: C(X) \otimes M_k \to A$ be a nonzero homomorphism. Then
$\iota \circ (\mu \otimes {\rm id}_{\OA{m}})$ is \aab.

\vspace{0.6\baselineskip}

{\it Proof:}  By the lemma, it suffices to show that for $\ep > 0$ and
$\lambda_1, \dots, \lambda_l \in X$, there exists unitary $v$
and nonzero \mops\ $q_1, \dots, q_l$ in $\mu(1)A\mu(1)$,
all commuting with
$\mu(1 \otimes M_k)$, such that
$\pa v - \mu(u \otimes 1) \pa < \ep$
 and $q_i v q_i = \lambda_i q_i$ for each $i$.
(The required objects in $B$ are then gotten by tensoring with the
identity in $\OA{m}$ and applying $\iota$.)
It actually suffices to choose a rank one projection $e \in M_k$, and
approximate $\mu(u \otimes e)$ by a unitary with this property in
$\mu(e) A \mu(e)$. We can then ignore the commutant condition.

If $X = S^1$ and ${\rm sp}(\mu (u \otimes e)) = S^1$ (evaluated in
$\mu(e) A \mu(e)$), and if $\mu(e) A \mu(e)$ is again a simple
direct limit with no dimension growth, then
the desired approximation is now Lemma 5.2 of \cite{Ph4}.
The proof
of the general case uses essentially the same reasoning as the
proof of that lemma. \QED

\vspace{0.6\baselineskip}

{\bf 1.10 Corollary} Let $m$ be even and let $B$ be a \pisca .
Let $X \subset S^1$ be compact, let
 $\varphi :  C(X) \otimes \OA{m} \to B$
be a homomorphism, and let $k \geq 1$. Then $\varphi$ is \aab\ if
 and only if
${\rm id}_{M_k} \otimes \varphi :
                   C(X) \otimes M_k (\OA{m}) \to M_k (B)$
is \aab .

\vspace{0.6\baselineskip}

{\em Proof:} Let $\varphi$ be \aab , and let $X_0 \subset X$
 be as before.
Let $\ep > 0$ and
$\lambda_1, \dots, \lambda_l \in X_0$. Choose a homomorphism
$\sigma_0 : \OA{m} \to B$, a unitary $v_0$, and projections $q_i^{(0)}$
for $\varphi$ as in the condition of the lemma. Then set
$\sigma = {\rm id}_{M_k} \otimes \sigma_0$, $q_i = 1 \otimes q_i^{(0)}$,
and $v = 1 \otimes v_0$. This shows that
 ${\rm id}_{M_k} \otimes \varphi$
is \aab .

Conversely,  let
${\rm id}_{M_k} \otimes \varphi$ be \aab , and let $\ep > 0$ and
$\lambda_1, \dots, \lambda_l \in X_0$. Without loss of generality
we may assume the $\lambda_i$ are distinct.
Let $\{e_{\mu \nu}\}$ be a system of matrix units in $M_k$. Choose
$\delta > 0$ such that if $\{f_{\mu}\}$ are projections in $M_k (B)$
such that $\|e_{\mu\mu} \otimes 1 - f_{\mu}\| < \delta$, then there is a
unitary $z \in M_k (B)$ such that $z f_{\mu} z^* = e_{\mu\mu} \otimes 1$
and $\| z - 1 \| < \ep/3$. Now choose
$\sigma : M_k (\OA{m}) \to M_k (B)$ and $v, q_1, \dots, q_l \in M_k(B)$
for ${\rm id}_{M_k} \otimes \varphi$
as in Lemma 1.7, using
 $\min(\delta, \ep/3)$ for $\ep$, and taking the standard matrix
units to be $\{e_{\mu\nu}\}$.
Replace $\sigma$, $v$, and $q_1, \dots, q_l$ by their conjugates by
$z$. This gives $\sigma$, $v$, and $q_1, \dots, q_l$ as in the lemma,
with norm estimate $\ep$, and with
$\sigma (e_{\mu\mu} \otimes 1) = e_{\mu\mu} \otimes 1 $.

Now set $\sigma_0 (a) = \sigma (e_{11} \otimes a)$
for $a \in \OA{m}$, and set
 $v_0 = e_{11} v e_{11}$ . For fixed $i$, choose $\mu$ such that
$( e_{\mu\mu} \otimes 1) q_i (e_{\mu\mu} \otimes 1) \neq 0 $,
and set
\[ q_i^{(0)} = ({\rm id}_{M_k} \otimes \varphi)(e_{1\mu} \otimes 1) q_i
({\rm id}_{M_k} \otimes \varphi)(e_{\mu 1} \otimes 1) . \]
It is easy to check
the commutation relations and norm estimates required to satisfy the
conditions of Lemma 1.7. (That the $q_i^{(0)}$ are mutually orthogonal
follows from the relations
$q_i^{(0)}v_0 = v_0q_i^{(0)} = \lambda_i q_i^{(0)}$ and the fact
that the $\lambda_i$ are distinct.)  This shows that $\varphi$
is \aab . \QED

\vspace{0.6\baselineskip}

{\bf 1.11 Lemma} Let $X \subset S^1$ be compact and connected, and let
$p \in A = C(X) \otimes M_n (\OA{m})$ be a nonzero projection.
Then there are $l \leq k$ and an isomorphism
$C(X) \otimes M_k (\OA{m}) \cong C(X) \otimes M_n (\OA{m})$ which sends
$C(X) \otimes M_l (\OA{m})$, regarded as a corner in
$C(X) \otimes M_k (\OA{m})$, onto $pAp$.

\vspace{0.6\baselineskip}

{\em Proof:} We identify $A$ with $C(X, M_n (\OA{m}))$. Theorem B
of \cite{Zh3} implies that $p$ is homotopic to, and hence unitarily
equivalent to, a constant projection $q$. (If $X = S^1$, this uses the
fact that $K_1 (\OA{m}) = 0$.)
Conjugating by this unitary, we may assume that $p$ itself is
a constant projection. This reduces us to consideration of the case
$X$ is a point. Thus, we assume $A = M_n (\OA{m})$. Without loss of
generality, we also assume $p \neq 1$.

Since $K_0 (\OA{m})$ is finite cyclic and generated by $[1]$, we can
find $0 < l< k$ such that $l [1] = [p]$ and $k [1] = n [1]$ in
$K_0 (\OA{m})$. Since $A = M_n (\OA{m})$ is purely infinite and simple,
there are $l$ orthogonal projections in $A$, each Murray-von Neumann
equivalent to $1_{\OA{m}}$, which sum to $p$, and $k - l$ such
projections which sum to $1 - p$. Since all $k$ projections are
Murray-von Neumann equivalent to  $1_{\OA{m}}$, they induce
an isomorphism $M_k (\OA{m}) \cong A$ which sends
$M_l (\OA{m})$ onto $pAp$. \QED

\vspace{0.6\baselineskip}

{\bf 1.12 Lemma} Let $B$ be a \pisca ,
let $X \subset S^1$ be compact and
have finitely many connected components, and let $m$ be even. If
$\varphi: C(X) \otimes M_k (\OA{m}) \to B$ is \aab , and
$p \in C(X) \otimes M_k (\OA{m})$ is a nonzero projection,
then $\varphi_0 = \varphi |_{p[C(X) \otimes M_k (\OA{m})]p}$
is also \aab .

\vspace{0.6\baselineskip}

{\em Proof:} We first prove this in the special case in which
$p$ is the characteristic function
$\chi_{Y}$ of some closed and open subset $Y \subset X$.

Let $F$ be a generating set for $C(Y) \otimes M_k (\OA{m})$, which we
regard as a subalgebra of $C(X) \otimes M_k (\OA{m})$, and let $G$ be
a generating set for $C(X) \otimes M_k (\OA{m})$ which contains $F$ and
$u \otimes 1$. Let $\ep > 0$, and let
$\psi : C(Y) \otimes M_k (\OA{m}) \to B$ be as in Definition 1.4.
When necessary, regard $\psi$ as defined on all of $A$ by setting it
equal to zero on $(1-p) A (1-p)$.
Choose $\dt > 0$ such that if $w_1$ and $w_2$ are unitaries in a
\CA\  with spectrum contained in $X$, and $\| w_1 - w_2 \| < \dt$,
then the projections
$\chi_{Y} (w_1)$ and $\chi_{Y} (w_2)$ are so close that there is a
unitary $z$ which conjugates one to the other and satisfies
$\| z - 1 \| < \ep / 4$.
Let $\rho = \min (\dt, \ep/2)$. Since $\varphi$ is \aab , we have
$\varphi \tdsum \psi \aueeps{\rho} \varphi$ with respect to $G$.
Conjugating the first of these
by a suitable unitary, we may assume they actually agree to within
$\rho$ on $G$. Now,
with $w_1 = (\varphi \tdsum \psi) (u \otimes 1)$  and
$w_2 = \varphi (u \otimes 1)$, conjugate $\varphi \tdsum \psi$
by $z$ as above.
The resulting \hm s agree exactly on
$p = \chi_Y (u \otimes 1)$ and to within $2 \ep /4 + \rho \leq \ep$
on $F$. Therefore the cutdowns by $\varphi (p)$ agree to within
$\ep$ on $F$. These cutdowns are $\varphi_0 \tdsum \psi$ and $\varphi_0$
respectively. So the special case is proved.

Combining the special case just proved with Remark 1.5, we see that
it now suffices to prove the lemma when $X$ is connected.  Lemma 1.11
reduces this proof to two applications of Corollary 1.10. \QED

\vspace{0.4\baselineskip}

To deal with maps from \CSalg s to nonsimple \CSalg s, we introduce
the following definition. We include the injectivity condition
because it is needed in Section 5.

\vspace{0.6\baselineskip}

{\bf 1.13   Definition} Let $D$ be any \CA , let $X$ be compact, and let
$\varphi : C(X) \otimes \OA{m} \to D$ be a \hm . Then $\varphi$ is
{\em permanently \aab} if whenever $B$ is a \pisca\  and
$\mu : D \to B$ is any nonzero \hm , then $\mu \circ \varphi$
is injective
and \aab .

\vspace{0.6\baselineskip}

To see that the definition is not vacuous, note
that any injective \aab\  \hm\  to a \pisca\  satisfies
the conditions.

The next two lemmas provide all the
permanently \aab\  \hm s we will need in this paper.

\vspace{0.6\baselineskip}

{\bf 1.14 Lemma} Let $A$ be a \CSalg , and let $\varphi: A \to D$
be a permanently \aab\  \hm . Let $E$ be a nonzero \CA , and let
$\psi: D \to E$ be any \hm\  such that $\psi(D)$ is not contained
in any proper ideal of $E$. Then $\psi \circ \varphi$ is
again permanently \aab .

\vspace{0.6\baselineskip}

{\em Proof:} Let $B$ be a \pisca , and let $\eta: E \to B$ be a
nonzero \hm . Then $\eta \circ \psi \neq 0$, so we apply the
definition to $\varphi$. \QED

\vspace{0.6\baselineskip}

{\bf 1.15 Lemma} Let $m$ be even, let $B$ be a purely infinite simple
\CA , and let $Y$ be a compact Hausdorff space.
Let $X$ be either $S^1$, a closed arc in $S^1$, or a point in $S^1$.
 Then
there exists a permanently \aab\  homomorphism
 $\varphi : C(X) \otimes M_n(\OA{m})  \to C(Y) \otimes B$ such
 that $[\varphi] = 0$ in
$KK^0 (C(X) \otimes M_n(\OA{m}), C(Y) \otimes B)$, and
such that for every $\ep > 0$ and any finite subset
$G \subset C(X)\otimes M_n(\OA{m}),$
there is a \hm\  $\psi : C(X) \otimes M_n(\OA{m})  \to C(Y) \otimes B$
satisfying:

(1) $\| \psi (g) - \varphi (g) \| < \ep$
for all $g\in G.$

(2) There are $x_1, \dots, x_l \in X$ and \hm s $\psi_1, \dots, \psi_l$
from $M_n (\OA{m})$ to $C(Y) \otimes B$,
whose classes in $KK^0 (M_n(\OA{m}), C(Y) \otimes B)$
are all zero, such that
$\psi_1(1),...,\psi_l(1)$ are \mops\  and such that
$\psi(f \otimes a) = \sum_{i=1}^{l} f(x_i) \psi_i (a)$ for $f \in C(X)$
and $a \in M_n (\OA{m})$.

(3) $\psi(1) = \varphi(1)$.

\noindent
If $[1] = 0$ in $K_0 (B)$, then $\varphi$ may be chosen unital.

\vspace{0.6\baselineskip}

{\em Proof:} By Corollary 1.10, we need only consider the case $n=1.$
(To see that this applies to the last statement, assume $[1] = 0$ in
$K_0 (B)$. Then $1_B$ is Murray-von Neumann equivalent to $1_{M_n (B)}$
in $M_n (B)$. Therefore $B \cong M_n (B)$, so that a unital \hm\  to
$C(Y) \otimes M_n (B)$ is the same as one to $C(Y) \otimes B$.)

Next, we observe that it suffices to prove the lemma with $B = \OA{2}$,
with $Y$ a one point space, and merely requiring injectivity in place
of the permanently \aab\  condition.
To see this, suppose that $\varphi_0 : C(X) \otimes \OA{m} \to \OA{2}$
has the required properties. It is easy to find a nonzero (hence
injective) \hm\  $\sigma_0 : \OA{2} \to B$, unital if
[1] = 0 in $K_0 (B)$. Define
$\sigma : \OA{2} \to C(Y) \otimes B$ by
$\sigma (a) = 1 \otimes \sigma_0 (a)$, and set
$\varphi = \sigma \circ \varphi_0$. Then
$[\varphi] = [\sigma] \times [\varphi_0] = 0$ in
$KK^0 (C(X) \otimes \OA{m}, C(Y) \otimes B)$. It is furthermore
clear that for every $\ep > 0$ and finite
$G \subset C(X) \otimes \OA{m}$
there is a \hm\  $\psi$ as required in the lemma. Finally, if
$C$ is some other \pisca , and
$\lambda: C(X) \otimes \OA{m} \to C(Y) \otimes B$ is any nonzero \hm ,
then $\lambda \circ \sigma$ is nonzero. Since $\OA{2}$ is simple,
$\lambda \circ \sigma$ is injective. The condition on the existence
of the map $\psi$ for $\varphi_0$, combined with an easy argument
using Lemma 1.7, now shows that $\lambda \circ \varphi$ is \aab .
Thus $\varphi$ is permanently \aab . This completes the reduction to
the case $B = \OA{2}$ etc.

 Since $K_* (\OA{2}) = 0$, the  Universal Coefficient Theorem \cite{RS}
implies that $KK^0 (\OA{2}, \OA{2}) = 0$. Therefore
 $KK^0 (A, \OA{2}) = 0$ for any separable nuclear \CA\  $A$.
Thus, in the special case
$B = \OA{2}$,  we can ignore
the requirement that \hm s be zero in $KK$-theory.

Let $D$ be the $2^{\infty}$ UHF algebra. We show that
if $X \subset S^1$ is as in the lemma, then $D$ contains
a unitary $v$ whose spectrum is $X$. If $X$ is a point, this is trivial.
Otherwise, note that $D$ is a non-elementary simple $C^*$-algebra.
By page 61 of \cite{AS}, there exists a selfadjoint $h\in D$ such that
${\rm sp}(h)=[0,1].$
Exponentiating, we obtain a unitary whose
spectrum is $S^1$ or any given closed arc in $S^1$.

Define $\mu: C(X) \to D$ by $\mu(u) = v$.
Further choose a unital \hm\ $\tau : \OA{m} \to \OA{2}$. (It is easy
 to find
such a thing.) Note that $D \otimes \OA{2} \cong \OA{2}$ by \cite{Rr1}.
 Define $\varphi$ to be the composite
$$
 C(X) \otimes \OA{m} \stackrel{\mu \otimes \tau}{\longrightarrow}
  D \otimes \OA{2} \stackrel{\cong}{\longrightarrow} \OA{2} .
$$
The approximating maps $\psi$ are obtained by replacing
$\mu (u)$ by approximating unitaries with finite spectrum, which
exist since $D$ is a UHF algebra.  \QED

\newpage

\section{ Homotopies of asymptotic morphisms}

\vspace{\baselineskip}

Results of D\v{a}d\v{a}rlat and Loring \cite{DL}
imply that if $\varphi_0$, $\varphi_1 : \So{m} \to B$
are two homomorphisms with the same class in
$KK^0 (\So{m} , B)$, and if $B$ is stable, then $\varphi_0$
and $\varphi_1$ are homotopic via asymptotic morphisms.  (Asymptotic
morphisms were defined in \cite{CH}.)
We are going to need a similar statement for \hm s from the unital
\CA\  $\SO{m}$. Of course, some conditions will have to be imposed,
some obvious (both \hm s nonzero), some not so obvious. Furthermore,
we don't actually need a homotopy of asymptotic morphisms.
What we need, and what we construct, is weaker in two ways. First,
for each value $\alpha$ of the homotopy parameter, instead of an
asymptotic morphism we merely produce a single linear map which
is multiplicative to within a prespecified
 $\ep > 0$ on our standard set of generators.
Second, the homotopy we construct does not agree exactly with the
original maps at the endpoints, but only to within $\ep$ on the
generators.  This section,
the longest one in the paper, is devoted to producing such a thing
under suitable conditions. This is done in Lemma 2.10.

We will need to know that two nonunital nonzero
 \hm s from $\OA{m}$ to $B$,
having the same class in $KK$-theory, are homotopic. This requires
R\o rdam's results \cite{Rr1}, so we will have to take $m$ even.
We will also have to assume $B$ is a \pisca\  (although not
necessarily stable). And we will need some other technical
conditions. Lemma 2.11, the last lemma in this section, shows how
to force some of them to hold.

The first part of this section is devoted to showing that
certain hereditary subalgebras of $C([0,1], B)$ have increasing
approximate identities consisting of projections. This fact is
an essential technical step in the proof of Lemma 2.10.
In fact, the hereditary subalgebras we consider can be shown, with
only a few more paragraphs, to be isomorphic to
$C([0,1], B) \otimes {\cal K}$.

It is probably true that, under suitable conditions, two
\hm s from $\SO{m}$ to $B$ with the same class in $KK^0 (\SO{m}, B)$
are homotopic via asymptotic morphisms. The conditions should be:
$B$ purely infinite simple, $m$ even (at least for now), and both
\hm s nonunital and ``absorbing up to homotopy'' in a suitable
sense. We don't need such a strong result, and so we don't try to
prove it, but this should be contrasted with the paper of
D\v{a}d\v{a}rlat and Loring \cite{DL}, in which the domain algebra is
never unital. Note that R\o rdam's work easily implies such a result for
$\OA{m}$ in place of $\SO{m}$, even using \hm s instead of
asymptotic morphisms. (See Lemma 2.9.)

We begin by summarizing in a convenient form some standard
approximation results in \CA s.

\vspace{0.6\baselineskip}

{\bf 2.1 Lemma} There exist functions $\beta_1$, $\beta_2$,
 $\beta_3$, $\beta_4$,
and $\beta_5^{(m)}$ from $[0, \infty)$ to $[0, \infty]$ which are
 nondecreasing
and satisfy $\lim_{t \to 0} \beta_i (t) = 0$, and which provide the
following estimates for approximation problems in a general \CA\ $A$.
In (5), we write simply $\beta_5$  when $m$ is understood. In all parts,
it is to be understood that when $\beta_i (\eta) = \infty$, the elements
claimed may in fact not exist.

(1) If $p_0 \in A$ is selfadjoint and satisfies
 $\| p_0^2 - p_0 \| < \eta$, then there exists a projection $p \in A$
such that $\| p - p_0 \| < \beta_1 (\eta)$.

(2) If $p$, $q \in A$ are projections such that $\| pq - q \| < \eta$,
then there exists a projection $p' \in A$ such that $p' \geq q$ and
$\| p' - p \| < \beta_2 (\eta)$.

(3) If $p$, $q \in A$ are projections such that $\| pq - q \| < \eta$
(as in (2)),
then there exists a projection $q' \in A$ and a unitary
path $t \mapsto v_t \in \widetilde{A}$ such that $p \geq q'$,
$\| q' - q \| < \beta_3 (\eta)$, $v_0 = 1$, $v_1 q' v_1^* = q$, and
$\| v_t - 1 \| < \beta_3 (\eta)$ for all $t$.

(4) If $p$, $q \in A$ are projections and $s_0 \in A$ satisfies
$\| p - s_0^* s_0 \| < \eta$ and $\| q - s_0 s_0^* \| < \eta$,
then there is a partial isometry $s \in A$, given by
$s = (p s_0 q)[(p s_0 q)^* (p s_0 q)]^{-1/2}$ (functional calculus
evaluated in $qAq$), which satisfies $s^*s = p$,
$ss^* = q$, and $\| s - s_0 \| < \beta_4 (\eta)$.

(5) If $A$ is unital, and $s_j$, $w_j \in A$ (for $j = 1, \dots, m$)
 are partial isometries
such that $s_j^* s_j = 1$,
 $\sum_i s_i s_i^* = 1$,
 the projections $q_j = w_j w_j^*$
and $q = \sum_i q_i$  satisfy $w_j^* w_j  = q$
for each $j$, and
\[
\| q_j s_j q - w_j \| < \eta,
\]
then there exist partial isometries $y_j$ such that
$y_j^* y_j = \sum_i y_i y_i^* = 1-q$ and
\[
\| s_j - (w_j + y_j) \| < \beta_5^{(m)} (\eta).
\]

\vspace{0.6\baselineskip}

{\it Proof:} Parts (1)--(4) are well known and have standard proofs,
which we omit. We sketch the proof of (5).

We start by observing that
\[
\|qs^*_j (1-q) s_jq\| \le \|qs^*_j (1-q_j) s_jq\|
  = \|q - qs^*_jq_js_jq\|<2\eta.
\]
Therefore
\[
\|(1-q) s_j q\| < (2 \eta)^{1/2}.
\]
We also estimate $\|q s_j (1-q)\|.$ We have
\[
\| s_k q - q_k s_k q\| =
    \|(s_k q - q_k s_k q  )^* (s_k q - q_k s_k q  )\|^{1/2} =
   \|q - q s^*_k q_k s_k q\|^{1/2} < (2 \eta)^{1/2} .
\]
So
\[
\| s_k q - w_k \| < \eta + (2 \eta)^{1/2} .
\]
Therefore
\[
\|q s_j (1-q)\| = \left\| \sum_{k=1}^m w_k w_k^* s_j (1-q)\right\|
  \le 2 m (\eta + (2\eta)^{1/2}) +
       \left\| \sum_{k=1}^m s_k q_k s_k^* s_j (1-q)\right\|.
\]
The last term is equal to $\|s_j q_j s_j^* s_j (1-q)\| = 0.$

The estimates on $\|(1-q)s_jq\|$ and $\|qs_j(1-q)\|$ show that $q$
approximately commutes with each $s_j.$

The desired result now follows from the fact that the defining
relations of $\OA{m}$ are exactly stable in the sense of Loring
\cite{Lr4}, \cite{Lr5}. This follows from Theorem 2.6 of \cite{Lr5}
and Corollary 2.24 of \cite{Bl0}. (Note that ``exactly semiprojective''
in \cite{Lr5} is the same as ``semiprojective'' in \cite{Bl0}.)
\QED

\vspace{0.6\baselineskip}

The following notation will be used throughout this section.

\vspace{0.6\baselineskip}

{\bf 2.2 Notation} If $B$ is a  \CA , $X$ is a compact Hausdorff space,
 and $D$
is a hereditary subalgebra of $C(X, B)$, we let
${\rm ev}_x : C(X, B) \to B$ be the evaluation map at $x \in X$, and
define
\[
D_x = {\rm ev}_x (D) = \{b(x): b \in D \}.
\]

\vspace{0.6\baselineskip}

{\bf 2.3 Lemma} Let the notation be as in 2.2. Then $D_x$ is a
hereditary $C^*$-subalgebra of $B$. (In particular, it is closed.)
Moreover, if $b \in C(X,B)$ satisfies $b(x) \in D_x$ for all
$x \in X$, then $b \in D$.

\vspace{0.6\baselineskip}

{\em Proof:} $D_x$ is the image of the \CA\ $D$ under the
\hm\  ${\rm ev}_x: C(X,B) \to B$, and is therefore a $C^*$-subalgebra of
$B$.  It remains to prove that $D_x$ is hereditary. So let $b_0 \in D_x$
and $a_0 \in B$ satisfy $0 \leq a_0 \leq b_0$. Choose $b \in D$ such
that $b(x) = b_0$; using standard manipulations we may assume that
$b \geq 0$. Now $\{b_0^{1/n}\}$ is an approximate identity for
the hereditary subalgebra of $B$ generated by $b_0$; in fact, for
any $c$ in this subalgebra, we actually have
$b_0^{1/n}cb_0^{1/n} \to c$. Now let $a \in C(X, B)$ be
the constant function with value $a_0$. Then
$b^{1/n}ab^{1/n} \in D$, and
$(b^{1/n}ab^{1/n})(x) \to a_0$. Therefore $a_0 \in D_x$.

It remains to prove the last statement. Let $\ep > 0$. For each
 $x \in X$ choose, using an approximate identity for $D$ and the first
part of the lemma, an element $d_x \in D$ such that
$\| (d_x b)(x) - b(x)) \| < \ep/2$. Choose an open set $U_x \subset X$
such that $\| (d_x b)(y) - b(y)) \| < \ep$ for $y \in U_x$. Cover
$X$ by finitely many of these sets, and let $\{ f_i: 1 \leq i \leq n\}$
be a partition of unity subordinate to this open cover. Let
${\rm supp}(f_i) \subset U_{x_i}$. Define
$d(x) = \sum_i d_{x_i}(x) (f_i b)(x)$. Then $d \in DB$ and
$\| d - b \| < \ep$. Therefore
 $b \in \overline{DB}$. Similarly
 $b \in \overline{BD}$.
Therefore  $b \in \overline{DB} \cap \overline{BD} = D$.
\QED

\vspace{0.6\baselineskip}

{\bf 2.4 Lemma} Let the notation be as in 2.2.
Let $Z \subset X$ be closed, and let $x_0 \not\in Z$.
Let $p \in D$ be a projection, and let $e \in D_{x_0}$
be a projection homotopic to $p(x_0)$ in $D_{x_0}$. Then there exists a
projection $q \in D$ such that $q(x) = p(x)$
for all $x \in Z$, $q(x_0) = e$, and $q$ is homotopic to $p$.

\vspace{0.6\baselineskip}

{\em Proof:} Without loss of generality, assume $B$ is unital.
Let $t \mapsto v(t)$ be a unitary path in
 $\widetilde{D_{x_0}} \subset B$
such that $v(0) = 1$, $v(1)p(x_0)v(1)^* = e$, and
$v(t) - 1 \in D_{x_0}$ for $t \in [0,1]$. Then $v$ can be
regarded as an element of $U_0(C([0,1], \widetilde{D_{x_0}}))$.
Since $C([0,1], \widetilde{D}) \to C([0,1], \widetilde{D_{x_0}})$
is surjective, there exists
 $w \in U_0 (C([0,1], \widetilde{D}))$ whose image
in $C([0,1], \widetilde{D_{x_0}})$ is $v$. An easy adjustment allows us
to assume that $w(0) = 1$ and $w(t) - 1 \in D$ for all $t$.
Now choose a continuous function
$f : X \to [0,1]$ such that $f|_Z = 0$ and $f(x_0) = 1$. Regarding
$w$ as a function from $X \times [0,1]$ to $B$, define
$c \in U(C(X,B))$ by $c(x) = w(x, f(x))$, and define
$q(x) = c(x)p(x)c(x)^*$.
The required homotopy from $p$ to $q$ will be given by
$t \mapsto w(x, tf(x))p(x)w(x, tf(x))^*$.

It remains only to show that these projections are actually in
$D$. This follows immediately from the fact
that the unitaries used to define them differ by 1 from
elements of $D$, which is a consequence of the last part of the
previous lemma.
\QED

\vspace{0.6\baselineskip}

{\bf 2.5 Lemma} Let $B$ and $D$ be as in Notation 2.2, and assume
 in addition that $X = [\alpha, \beta]$ is a closed interval,
$B$ is separable, purely infinite, and simple, and $D_t$ is
nonzero and nonunital for all $t$. Let $p \in D$,
 $e_\alpha \in D_{\alpha}$,
 and
$e_{\beta} \in D_{\beta}$ be projections such that,
 for $i = \alpha, \beta$
 we have
\[
e_i > p(i) \,\,\,\,\,\, {\rm and} \,\,\,\,\,\,
    [e_i - p(i)] = 0 \,\,\, {\rm in} \,\,\,  K_0 (B).
\]
Then there exists a projection $q \in D$ such that $q \geq p$,
 $q(\alpha) = e_\alpha$, and $q({\beta}) = e_{\beta}$.

\vspace{0.6\baselineskip}

{\em Proof:} Without loss of generality, assume $\alpha = 0$ and
$\beta = 1$. Further, replace $D$ by $(1-p)D(1-p)$. This allows us
to assume that $p = 0$. The new $D_t$ is still nonunital
and nonzero. (Note that $p(t)$ can't be an identity for $D_t$, since
$D_t$ doesn't have an identity.)

Each $D_t$ contains a positive element with norm greater than 1.
Therefore, for each $t$ the subalgebra
 $D$ contains a positive element $a$ such that
$\| a(t) \| > 1$. Using a partition of unity argument and Lemma 2.3,
we can produce $a \in D$ such that $a \geq 0$ and
 $\| a(t) \| > 1$ for all $t$.
Cutting down using functional calculus, we can assume $\| a(t) \| = 1$
for all $t$. Choose a partition
\[
0 = t_0 < t_1 < \cdots < t_n = 1
\]
such that $\| a(t) - a(t_i) \| < 1/16$ for $t \in [t_{i-1}, t_i]$.
Now $D_{t_i}$, being purely infinite and simple, has real rank $0$.
Therefore $a(t_i)$ can be approximated by selfadjoint elements with
finite spectrum. Since also $1 \in {\rm sp}(a(t_i))$, there is a
nonzero projection $f_i \in D_{t_i}$ such that
$\| f_i a(t_i) - f_i \| < 1/16$. Since  $f_i D_{t_i} f_i$
is purely infinite and $K_0 (f_i D_{t_i} f_i) = K_0 (D_{t_i}),$
there is a nonzero projection $f_i' \in f_i D_{t_i} f_i$
such that $[f_i']=0$. (See page 188 of \cite{Cu1}.)
Replacing $f_i$ by $f_i'$, we may assume
furthermore that $[f_i] = 0$ in $K_0 (D_{t_i})$.
We next observe that $x(t) = a(t)f_i a(t)$ satisfies
$x(t) \in D_t$ and
$\| x(t) - f_i \| < 1/4$ for $t \in [t_{i-1}, t_i]$.

We can therefore apply functional calculus to produce a continuous
function $q_i$ from $[t_{i-1}, t_i]$ to the projections in $B$
such that $q_i (t) \in D_t$ and $\| q_i (t) - f_i \| < 1/2$.
It follows that the classes of $q_i(t)$ and $f_i$ in $K_0 (B)$
are equal, and that $q_i(t) \neq 0$.
  Since the inclusion of $D_t$ in $B$ induces an isomorphism
on $K$-theory, it follows that $[q_i(t)] = 0$ in $K_0 (D_t)$.
In particular, $[q_i(t_i)] = [q_{i+1}(t_i)]$. Note that both
$q_i(t_i)$ and $q_{i+1}(t_i)$ are nontrivial. Therefore Theorem
1.1 of \cite{Zh4} implies that these two projections are homotopic.
Using the previous lemma, we now construct $q$ on $[0, t_1]$ so that
$q(0) = e_0 $ and $ q(t_1) = q_1 (t_1)$, on $[t_1, t_2]$
so that $q(t_1) = q_1 (t_1)$ and $q(t_2) = q_2 (t_2)$, etc.,
finishing by constructing $q$ on $[t_{n-1}, \frac{1}{2}(t_{n-1} + 1)]$
so that $q(t_{n-1}) = q_{n-1}(t_{n-1})$ and
$q(\frac{1}{2}(t_{n-1} + 1)) = q_n ( \frac{1}{2}(t_{n-1} + 1))$ and
on $[ \frac{1}{2}(t_{n-1} + 1), 1]$ so that
$q(\frac{1}{2}(t_{n-1} + 1)) = q_n ( \frac{1}{2}(t_{n-1} + 1))$ and
$q(1) = e_1$. Since the definitions agree on the points where the
intervals overlap, $q$ is in fact continuous. This is the desired
projection.
\QED

\vspace{0.6\baselineskip}

The hypotheses of the following proposition (as well as those of the
previous lemma) imply that $D$ is full in $C([0,1], B)$.
Therefore $D$ is stably isomorphic to  $C([0,1], B)$. In fact,
it turns out that
$D$ is isomorphic to $C([0,1], B) \otimes K$. (We won't prove
this, but it is an easy step from the conclusion of the next
proposition.) If we already knew this,  the
conclusions would be obvious. But we don't know how to prove such an
 isomorphism
except by using these results.

\vspace{0.6\baselineskip}

{\bf 2.6 Proposition}  Under the hypotheses of the previous lemma,
$D$ has an increasing approximate identity consisting of projections.

\vspace{0.6\baselineskip}

{\em Proof:} We will prove the following claim:
 For every positive element
$a \in D$, every projection $p \in D$, and every $\ep > 0$, there
exists a projection $q \in D$ such that
\[
\| qa-a \| < \ep \,\,\,\,\,\, {\rm and} \,\,\,\,\,\, \| qp - p \| < \ep.
\]
Given this, we first observe that the same result holds with the
stronger conclusion $q \geq p$ in place of $ \| qp - p \| < \ep$.
Indeed, with  $\beta_2 $ as in Lemma 2.1, we choose $\delta > 0$
such that $ \|a\| ( \delta + \beta_2 (\delta)) < \ep$, construct $q_0$
as above using $\delta$ in place of $\ep$, and use Lemma 2.1 (2)
to replace $q_0$ by $q$ such that $q \geq p$ and $\| q - q_0 \| <
\beta_2 (\delta)$. Using the stronger conclusion, we construct our
approximate identity by induction, letting $a$ and $\ep$ run
 independently through a countable dense subset of the positive part of
$D$ and the set $\{\frac{1}{n}: n \in {\bf N}\}$.

We now prove the claim. Without loss of generality, assume
 $\| a \| \leq 1$.

Since $D$ is separable, it has a strictly positive element $b$.
We first prove the following subclaim:
there exists a strictly decreasing sequence $\{\alpha_n\}$
of positive real numbers with $\alpha_n \to 0$ such that
\[
(\alpha_{n+1}, \alpha_{n}) \cap {\rm sp}(b(t)) \neq \emptyset
\]
for all $n$ and $t$. We construct $\alpha_n$ inductively, starting with
$\alpha_0 = 1$. Given $\alpha_n > 0$, we first observe that there is
$\beta_t \in  {\rm sp}(b(t))$
such that $\alpha_n > \beta_t > 0$. (Otherwise, since $b(t)$ is
strictly positive in $D_t$, this algebra would be unital.)
Continuity of the spectrum on selfadjoint elements ensures that
there is an open set $U_t \subset [0,1]$ such that for $s \in U_t$
we have
\[
\left(\frac{\beta_t}{2} ,  \frac{ \beta_t + \alpha_{n} }{2}\right)
        \cap {\rm sp}(b(s)) \neq \emptyset .
\]
Cover $[0,1]$ with finitely many of these sets, say $U_{t_1}, \dots,
U_{t_l}$, and choose $\alpha_{n+1}$
such that
\[
 0 < \alpha_{n+1} < \min(\alpha_n/2, \beta_{t_1}, \dots, \beta_{t_l}).
\]
This proves the subclaim.

Let $g_k : [0, \infty) \to [0,1]$ be the continuous function
which is $0$ on $[0, \alpha_{k+1}]$, $1$ on $[ \alpha_{k}, \infty)$,
and linear on $[\alpha_{k+1}, \alpha_k]$. Note that $g_{k+1} g_k = g_k$,
and that
the elements $g_k (b)$ form an approximate identity for
$D$.

Choose $\rho > 0$ such that
\[
6 \beta_3 (14 \rho) + 16 \rho < \ep \andeqn
3 \beta_3 (14 \rho) + 4 \rho < 1/2.
\]
(Here the function $\beta_3$ is from Lemma 2.1 (3).)
 Choose $k$
so large that
\[
\| g_l (b) a - a\| < \rho  \,\,\,\,\,\, {\rm and} \,\,\,\,\,\,
\| g_l (b) p - p\| < \rho
\]
for all $l \geq k$.
Next, choose a partition
\[
0 = t_0 < t_1 < \cdots < t_n = 1
\]
of $[0,1]$ such that, for all $t \in [t_{i-1}, t_i]$, we have
\[
\| g_{l} (b) (t) - g_{l} (b)(t_{i-1}) \| < \rho \andeqn
\| g_{l} (b) (t) - g_{l} (b)(t_i) \| < \rho
\]
for $l = k$ and $l = k + 3$, and
\[
\| a(t) - a(t_i) \| < \rho \andeqn
\| p(t) - p(t_i) \| < \rho.
\]

The relations
$g_{k+2}(b(t_i))g_{k+1}(b(t_i)) = g_{k+1}(b(t_i))$ and
$g_{k+1}(b(t_i))g_{k}(b(t_i)) = g_{k}(b(t_i))$
 imply,
using \cite{Bn2} and the fact that $D_{t_i}$
has real rank zero,
 the existence of a projection
$e_i^{(0)} \in D_{t_i}$ such that
$g_{k+2}(b(t_i)) \geq e_i^{(0)} \geq g_{k}(b(t_i))$.
The choice of the
sequence $\{\alpha_n\}$ implies that
${\rm sp}(g_{k+2}(b(t_i))) \cap (0,1) \neq \emptyset $.
Since $g_{k+2}(b(t_i)) e_i^{(0)} =e_i^{(0)}$, it follows that
$e_i^{(0)}$ is not an identity for the hereditary subalgebra $E_i$
 of $B$
generated by $g_{k+2}(b(t_i))$. Since this hereditary subalgebra
is purely infinite simple, all classes in its $K_0$-group are
represented by projections dominating $e_i^{(0)}$.
Therefore we can replace $e_i^{(0)}$ by a possibly larger projection
in $E_i$ whose class in $K_0 (E_i)$ is zero.
The  $K_0$-class of $e_i^{(0)}$ in $B$ is then zero as well.
Note that the equation
$g_{k+3}(b(t_i))g_{k+2}(b(t_i)) = g_{k+2}(b(t_i))$
implies that $g_{k+3}(b(t_i))x = g_{k+2}(b(t_i))$
for all $x \in E_i$. In particular, with our new choice of $e_i^{(0)}$,
we have this equation for $x = e_i^{(0)}$,
whence
\[
g_{k+3}(b(t_i)) \geq e_i^{(0)} \geq g_{k}(b(t_i)).
\]
Similarly, there exists a projection $e_i^{(1)} \in D_{t_i}$ such that
\[
g_{k+6}(b(t_i)) \geq e_i^{(1)} \geq g_{k+3}(b(t_i))
\]
and $[e_i^{(1)}] = 0$ in $K_0 (B)$. In particular,
 $ e_i^{(1)} \geq e_i^{(0)} $.

Temporarily fix $i$, and work over the interval $[t_{i-1}, t_{i+1}]$.
Observe that the element
$x(t) =  g_{k+3}(b(t))e_i^{(0)}g_{k+3}(b(t))$
is in the hereditary subalgebra $F_t$ generated by $g_{k+3}(b(t))$
and satisfies
$\| x(t) -  e_i^{(0)}\| < 2 \rho$.
Since $\rho < 1/4$, functional calculus yields a projection
$f_i^{(0)} (t)$,
depending continuously on $t$, such that
\[
\| f_i^{(0)} (t) -  e_i^{(0)}\| < 4 \rho.
\]
Similarly, we obtain a projection $f_i^{(1)} (t)$
 in the hereditary subalgebra generated by $g_{k+6}(b(t))$,
depending continuously on $t \in [t_{i-1}, t_{i+1}]$, such that
\[
\| f_i^{(1)} (t) -  e_i^{(1)}\| < 4 \rho.
\]

We now observe that
\[
\| f_i^{(1)} (t) f_i^{(0)} (t) -  f_i^{(0)} (t) \| < 3(4 \rho)
       < 14 \rho,
\]
one term $4 \rho$ coming from each replacement of an
$f_i^{(\nu)} (t)$ by an $e_i^{(\nu)}$.
We also observe that
\begin{eqnarray*}
\| e_i^{(1)} e_{i+1}^{(0)} - e_{i+1}^{(0)}\|
 & = & \| e_i^{(1)} g_{k+3}(b)(t_{i+1}) e_{i+1}^{(0)} -
                   g_{k+3}(b)(t_{i+1}) e_{i+1}^{(0)}\| \\
 & \leq &  \| e_i^{(1)} g_{k+3}(b)(t_{i+1}) -  g_{k+3}(b)(t_{i+1}) \| \\
 & \leq &  2 \| g_{k+3}(b)(t_{i+1}) -  g_{k+3}(b)(t_{i}) \|
                 < 2 \rho,
\end{eqnarray*}
since $e_i^{(1)} g_{k+3}(b)(t_{i}) =  g_{k+3}(b)(t_{i})$.
An estimate similar to the one at the beginning of this paragraph
therefore shows that
\[
\| f_i^{(1)} (t) f_{i+1}^{(0)} (t) -  f_{i+1}^{(0)} (t) \|
     < 14 \rho
\]
for  $t \in [t_{i}, t_{i+1}]$.

Using Lemma 2.1 (3), we obtain a projection $r_i (t) \in D_t$
 and a unitary path $s \mapsto v_s(t) \in \widetilde{D_t}$
for $t \in [t_{i-1}, t_{i}]$ and $s \in [0,1]$,
both varying continuously with $t$,
 such that
$r_i (t) \leq  f_{i}^{(1)} (t) $,
$v_0 (t) = 1$, $v_1(t) r_i (t) v_1(t)^* = f_{i}^{(0)} (t) $, and
\[
\| r_i (t) -   f_{i}^{(0)} (t) \| < \beta_3 (14 \rho)
 \andeqn
\| v_s (t) - 1 \| < \beta_3 (14 \rho).
\]
We similarly obtain $r_i'(t)$ and $v_s'(t)$ satisfying all the same
conditions, except with $f_i^{(\nu)} (t)$
replaced by $f_{i-1}^{(\nu)} (t)$.
We now define a continuous projection $q_i$ on $[t_{i-1}, t_{i}]$
by letting $t_i' = \frac{1}{2}(t_{i-1} + t_{i})$  and setting
\[
q_i ( t) =
      v_{1-\alpha}'(t)^* f_{i}^{(0)} (t) v_{1-\alpha}'(t)
  \,\,\,\,\,\, {\rm for} \,\,\, t = (1-\alpha ) t_{i-1} + \alpha t_i'
   \,\,\, {\rm and}
   \,\,\, \alpha \in [0,1]
\]
and
\[
q_i ( t) =
      v_{\alpha}(t)^* f_{i}^{(0)} (t) v_{\alpha}(t)
  \,\,\,\,\,\, {\rm for} \,\,\, t = (1-\alpha ) t_{i}' + \alpha t_i
   \,\,\, {\rm and}
   \,\,\, \alpha \in [0,1].
\]
This gives $q_i (t_i') =  f_{i}^{(0)} (t_i') $ (with either definition),
and
\[
q_i (t_{i-1}) = r_i'(t_{i-1}) \leq f_{i-1}^{(1)} (t_{i-1})
   \andeqn
q_i (t_{i}) = r_i(t_{i}) \leq f_{i}^{(1)} (t_i).
\]

Let $t \in [t_i', t_i]$. Then for suitable $\alpha \in [0,1]$, we
have
\begin{eqnarray*}
\| q_i (t) - e_i^{(0)} \| & \leq &
     \| q_i (t) - r_i (t) \| +  \| r_i (t) - f_i^{(0)} (t) \| +
             \| f_i^{(0)} (t)   - e_i^{(0)} \| \\
 & < & 2 \| v_{\alpha}(t) - 1 \| + \beta_3 (14 \rho) + 4 \rho
         < 3 \beta_3 (14 \rho) + 4 \rho.
\end{eqnarray*}
A similar estimate holds for $t \in [t_{i-1}, t_i']$. Since
$3 \beta_3 (14 \rho) + 4 \rho < \frac{1}{2}$, it follows that
$[q_i (t)] = 0$ in $K_0 (B)$, and hence also in $K_0 (D_t)$, for
all  $t \in [t_{i-1}, t_i]$.

Since $D_{t_i}$ is not unital, there is a projection $r_i \in D_{t_i}$
such that $r_i > f_i^{(1)} (t_i)$ and $[r_i] = 0$ in $K_0 (D_{t_i})$.
Lemma 2.5 provides a continuous projection
 $q_i' :  [t_{i-1}, t_i] \to B$ such that $q_i' (t) \in D_{t}$,
$q_i' (t_{i-1}) = r_{i-1} - q_{i}(t_{i-1})$, and
$q_i' (t_{i}) = r_{i} - q_{i}(t_{i})$. Now define
 $q(t) = q_i (t) + q_i'(t)$ for  $t \in [t_{i-1}, t_i]$.
Then $q$ is well defined and continuous, since at the overlap points
$t_i$ both definitions yield $r_i$. Furthermore, $q \in D$ by Lemma 2.3,
and $q(t) \geq q_i (t)$ whenever $t \in [t_{i-1}, t_i]$.

We now estimate, for  $t \in [t_{i-1}, t_i]$:
\begin{eqnarray*}
\lefteqn{\| q_i(t) a(t) - a(t) \|} \\
 & \leq & \| q_i(t)\| \|a(t) - a(t_i)\| +
    \| q_i(t)- e_{i}^{(0)}\| \|a(t_i)\|    \\
 & &  \mbox{} + \|e_{i}^{(0)}\| \|a(t_i) - g_k(b)(t_i)a(t_i) \|
            + \|  e_{i}^{(0)}g_k(b)(t_i) - g_k(b)(t_i) \| \|a(t_i) \| \\
 & & \mbox{} + \| g_k(b)(t_i) a(t_i) - a(t_i) \| + \|a(t_i) - a(t) \| \\
 & < & \rho + (3 \beta_3 (14 \rho) + 4 \rho) + \rho + 0 + \rho + \rho \\
 & = & 3 \beta_3 (14 \rho) + 8 \rho.
\end{eqnarray*}
Therefore
\begin{eqnarray*}
\lefteqn{\| q(t) a(t) - a(t) \|} \\
 & & \leq \| q(t)\| \|a(t) - q_i(t)a(t)\| +
       \|q(t)q_i(t) - q_i(t)\| \| a(t)\|+ \|q_i(t)a(t) - a(t)\|  \\
 & & < 3 (\beta_3 (14 \rho)+8 \rho) + 0 + (3 \beta_3 (14 \rho)+8 \rho)
          =  6 \beta_3 (14 \rho) + 16 \rho \leq \ep.
\end{eqnarray*}
A similar argument, with $p(t)$ in place of $a(t)$, shows that
\[
\| q(t) p(t) - p(t) \| < \ep.
\]
Thus, $q$ is the desired projection.
\QED

\vspace{0.6\baselineskip}

{\bf 2.7 Lemma} Let $B$ be a purely infinite simple \CA\  (not
necessarily unital). Let $p,q \in C([0,1], B)$ be nonzero projections
 whose $K$-theory classes are equal. Then $p$ is
Murray-von Neumann equivalent to $q$.

\vspace{0.6\baselineskip}

{\it Proof:} Since $B$ is purely infinite and simple, and since $p(0)$
and $q(0)$ are nonzero projections with the same class in $K$-theory,
there is $s(0) \in B$ such that
\[
s(0) s(0)^* = p(0) \andeqn s(0)^* s(0) = q(0).
\]
Standard methods give unitary paths $t \to u(t), v(t)$ in $\tilde{B}$
such that $u(0) = v(0) = 1$ and
\[
u(t) p(0) u(t)^* = p(t) \andeqn v(t) q(0) v(t)^* = q(t).
\]
The required partial isometry is given by $s(t) = u(t) s(0) v(t)^*$.
\QED

\vspace{0.6\baselineskip}

{\bf 2.8 Corollary} The approximate identity in Proposition
2.6 can be chosen so that the $K_0$-classes of all the projections in it
are trivial.

\vspace{0.6\baselineskip}

{\em Proof:} Let $\{p_k \}$ be an increasing approximate identity
consisting of projections. Without loss of generality, we may assume
it is strictly increasing. It suffices to find projections
$q_k \in D$ with $p_k \leq q_k \leq p_{k + 1}$ and $[p_k] = 0$ in
$K_0 (D)$. Lemma 2.7 implies that both $p_k$ and $p_{k + 1}$ are
Murray-von Neumann equivalent in C([0,1], B) to the constant projections
with values $p_k (0)$ and $p_{k + 1} (0)$ respectively. These
equivalences give an isomorphism
$p_{k + 1} D p_{k + 1} \to C([0,1],  p_{k + 1} (0) D_0 p_{k + 1} (0))$
(with $D_0$ as in Notation 2.2). Note that $p_k$ goes to the constant
projection with value $p_k (0)$. Since $B$ is purely infinite and
simple, we can choose a projection $e \in B$ such that
$p_k (0) \leq e \leq p_{k + 1} (0)$ and $[e] = 0$ in $K_0 (B)$.
Now take $q_k$ to be the inverse image in $D$ of $1 \otimes e$.
Note that $[1 \otimes e] = 0$ in
$K_0 (C([0,1],  p_{k + 1} (0) D_0 p_{k + 1} (0)))$,
so $[q_k] = 0$ in $K_0 (D)$.  \QED

\vspace{0.6\baselineskip}

{\bf 2.9 Lemma} Let $\psi_0, \psi_1 : \OA{m} \to B$ be two \hm s
from an even Cuntz algebra to a purely infinite simple \CA .
If $[\psi_0] = [ \psi_1 ]$ in $KK^0 (\OA{m}, B)$, and if
$\psi_0$ and $\psi_1$ are both nonunital and
nonzero, then $\psi_0$ is homotopic to $\psi_1$.

\vspace{0.6\baselineskip}

{\em Proof:} The condition $[\psi_0]=[\psi_1]$ in
$KK^0(\OA{m}, B)$ implies that $[\psi_0(1)]=[\psi_1(1)]$
in $K_0(B).$ Therefore $\psi_0 (1)$ is homotopic to $\psi_1 (1)$.
It follows that there is a path $t \mapsto W_t$ of unitaries in
 $\tilde{B}$
such that $W_0=1$ and $W_1^*\psi_1(1)W_1=\psi_0(1).$
So we may assume that $\psi_0(1)=\psi_1(1).$  We can now reduce to
the case that $B$ is unital: if not, replace $B$ by $pBp$ for some
projection $p \in B$ with $p > \psi_0 (1)$.

Theorem 3.6 of \cite{Rr1} provides
a unitary $V\in \psi_0(1)B\psi_0(1)$ such that
\[
\|V^*\psi_1(s_j)V-\psi_0(s_j)\|<1/(2m)
\]
for $j=1,2, \dots ,m.$
Since $B$ is purely infinite simple, there is a unitary
$V_0\in (1-\psi_0(1))B(1-\psi_0(1))$ such that
$[V+V_0]=0$ in $K_1(B).$ Replacing $V$ by $V + V_0$,
we may assume that $V$ is a unitary
in $B$ and is in the identity component of $U(B)$.
Therefore $\psi_1$  is homotopic to $V^*\psi_1V.$
 Thus, without loss of generality,
we may assume that
\[
\|\psi_1(s_j)-\psi_0(s_j)\|<1/(2m)
\]
for $j=1,2, \dots ,m.$
Then (compare \cite{Rr1}, 3.3) the unitary
$U\in \psi_0(1)B\psi_0(1)$ given by
\[
U = \sum_{j} \psi_1 (s_j) \psi_0 (s_j)^*
\]
satisfies
\[
\psi_1(s_j)=U\psi_0(s_j)
\]
for all $j$. The formula for $U$ implies that $\|U - 1 \| < 1/2$, so
there is a unitary path $t \mapsto U_t$ with $U_0 = 1$ and $U_1 = U$.
Define the required homotopy by $\psi_t (s_j) = U_t \psi_0 (s_j)$.
\QED

\vspace{0.6\baselineskip}

{\bf 2.10 Lemma} Let $B$ be purely infinite, simple, and separable. Let
$ t \mapsto \varphi_t$ be an asymptotic morphism from $\So{m}$ to
$C([0,1], B)$, with $m$ even.
Let $\psi_0, \psi_1 : \SO{m} \to B$ be homomorphisms.
Let $M_m \subset \OA{m}$ be the $m \times m$ matrix subalgebra generated
by the elements $s_i s_j^*$.  Assume:

(H1) Each $\varphi_t$ is linear and *-preserving, and
 $\sup_t \| \varphi_t \| < \infty$. (We demand neither contractivity nor
 positivity.)

(H2) $\varphi_t|_{C_0 (S^1 \setminus \{1\}) \otimes M_m}$ is a
 \hm\ for each $t$.

(H3) Whenever $a \in \So{m}$ is actually in
 $C_0 (S^1 \setminus \{1\}, p_i \OA{m} p_j)$,
then $\varphi_t (a^*) \varphi_t (a)$ is in the hereditary subalgebra of
$C([0,1], B)$
generated by $\varphi_t (C_0 (S^1 \setminus \{1\}, {\bf C}p_j))$.

(H4) For $i = 0,1$, the \hm s $\psi_i$ satisfy
${\rm ev}_i \circ \varphi_t = \psi_i |_{ \So{m}}$
for all $t$.

(H5) The \hm s $\psi_0$ and $\psi_1$ are both nonunital, and satisfy
$[\psi_0 |_{ {\bf C} \otimes \OA{m}}] =
      [\psi_1 |_{ {\bf C} \otimes \OA{m}}]$
in $KK^0 (\OA{m}, B)$.

(H6)
For all $\alpha \in [0,1]$ and $t \geq 0$, the spectrum
${\rm sp}({\rm ev}_{\alpha} (\varphi_t((u-1) \otimes s_1 s_1^*)))$ is
equal to the entire circle with radius $1$ and center $-1$.

Then for all $\ep > 0$ and all $T$, there is $t \geq T$, a unital
\hm\  $\rho : \OA{m} \to C([0,1], B)$, and a  unitary
$v \in \rho(1)C([0,1], B)\rho(1)$ such that:

(C1) $\| v \rho(s_j) - \rho(s_j) v \| < \ep$.

(C2) $\| (v - \rho(1)) - \varphi_t((u-1) \otimes 1) \|
           < \ep $ and
   $\| (v - \rho(1)) \rho(s_j) - \varphi_t((u-1) \otimes s_j) \|
           < \ep $.

(C3) $\| {\rm ev}_i (v) - \psi_i (u \otimes 1) \| < \ep$  and
 $\| {\rm ev}_i (\rho(s_j)) - \psi_i (1 \otimes s_j) \| < \ep$
 for $i = 0,1$.

\vspace{0.6\baselineskip}

Lemma 2.11 below will show that an arbitrary homotopy can be modified
in such a way that the hypotheses (H1)--(H3) of this lemma are
satisfied. In applications, hypothesis (H6) will be achieved by
forming the direct sum with a suitable homomorphism.

This lemma gives conditions under which two homomorphisms from
$\SO{m}$ to $B$, with the same class in $KK$-theory, can be
(almost) connected
by a homotopy of ``discrete asymptotic morphisms'' (indexed by
a discrete set rather than $[0, \infty)$). Presumably
one can actually get a homotopy of proper asymptotic morphisms.
Doing this would, however, make an already messy proof even worse,
and we do not need the stronger result.

The proof of this lemma is very long and technical. The basic idea,
however, is quite simple, and we therefore describe it before we begin.

For simplicity, assume in this sketch that $\varphi$ is actually a
\hm\  from $C_0 (S^1 \setminus \{1\}) \otimes \OA{m}$
to $C([0,1]) \otimes B$.
We really want a \hm\  from $\SO{m}$ to  $C([0,1]) \otimes B$, but
we will have to settle for an approximate \hm . The missing part is
a \hm\  from $\OA{m}$ to  $C([0,1]) \otimes B$.
Let $T = \{ \zeta - 1 : \zeta \in S^1 \setminus \{ 1 \} \}$, so that
$T \cup \{ 0 \} = {\rm sp}(u-1)$. To construct our \hm , let
$g, h \in C_c (T)$ satisfy $h(\zeta) \in T$ and
$h(\zeta) \approx \zeta$ for $\zeta \in T$,
$0 \leq g \leq 1$, and $gh = h$. Choose a ``large'' projection
$e_1$ in the hereditary $C^*$-subalgebra
$D$ generated by $\varphi ( h(u-1) \otimes p_1)$. (This subalgebra has
an increasing approximate identity of projections by Lemma 2.5 and
the hypothesis (H6).) Set $p_j = s_j s_j^* \in \OA{m}$. Note that
the element  $\varphi ( g(u-1) \otimes p_j)$ acts as an identity
for $D$. Thus, we can define projections
\[
e_j = \varphi ( g(u-1) \otimes s_j s_1^* ) e_1
                       \varphi ( g(u-1) \otimes s_1 s_j^*)
\andeqn   e = \sum_{j = 1}^m e_j,
\]
and partial isometries
\[
t_j = e_j \varphi ( g(u-1) \otimes s_j).
\]
Now $h(u-1)$ is close to $u-1$, and $e_1$ is fairly far out in an
approximate identity for an algebra which contains
$\varphi ( h(u-1) \otimes p_1)$. Therefore the elements we have
constructed approximately commute with
$\varphi ( (u-1) \otimes 1) + 1$. Furthermore, it is quite easy to
extend the $t_j$ to have larger initial and final projections
and still approximately commute with
$\varphi ( (u-1) \otimes 1) + 1$.

Unfortunately, there is a problem. We certainly have
$t_j t_j^* = e_j$, but there is no reason to have
$t_j^* t_j = e$. The best we can say is that $t_j^* t_j$ and $e$ are
close, not in norm but only in the strict topology. (Both are
``large'' in a certain hereditary subalgebra.) In order to get
objects satisfying the relations for $\OA{m}$, we need to
correct for this error. To have room for the correction, we
must repeat the argument of the previous paragraph using a
second pair of functions $g', h' \in C_c (T)$. These functions
are to satisfy the same relations as $g$ and $h$, and they are
supposed to be ``larger''. The precise condition is
$h' (\zeta) = \zeta$ on the support of $g$.
Let $e_j'$ and $e' = \sum_j e_j'$ be the new, larger, projections.
It is crucial for the proof that $e'$ strictly dominates
both $e$ and the common initial projection $t_j^* t_j$ of the
original set of partial isometries.

The correction procedure makes the proof much more complicated.
(Among other problems, we can't exactly get the domination referred
to at the end of the previous paragraph. We must settle for
an $\ep$ approximation.)
Even worse, we can't work with a \hm\ $\varphi$, but only with
an approximate \hm\  $\varphi_t$ taken from the given asymptotic
morphism. As a result, errors accumulate with every step in
the construction. The proof ends with a device to get things to
match up at the endpoints of the homotopy, which we have not
discussed in this sketch.

\vspace{0.6\baselineskip}

{\it Proof of Lemma 2.10:} We divide the proof into nine steps, of
which Steps 1 and 6 are further subdivided.

{\em Step 1:} Preliminaries. The purpose of this step is to
set things up for the rest of the proof.

{\em Step 1.1:} We begin by making
several reductions and definitions.

We usually write $\alpha$ for an element of $[0, 1]$, and use
function notation for the dependence of an element of $C([0, 1], B)$
on $\alpha$.
We use subscripts for the parameter $t$ of the asymptotic morphism.

 For simplicity, we assume that $B$ is unital.

 Since the conclusion of the lemma only involves approximations
to within $\ep,$ we can make a preliminary small perturbation of
$t\to \varphi_t.$ By making such a perturbation,
 we may assume that
${\rm ev}_{\alpha}\circ\varphi_t$ is a constant function of
 $\alpha$ for $\alpha$ in some neighborhood of
0 and also for $\alpha$ in some neighborhood of 1. We now
reparametrize the interval $[0,1]$ so as to be able to
assume that these neighborhoods are $[0,1/3]$ and $[2/3,1],$
respectively.

Set $p_j=s_js_j^* \in \OA{m}$.

{\em Step 1.2:} We now choose some small numbers $\eta$, $\dt$,
 $\dt'$, and $\dt''$, some useful functions defined on a circle
in ${\bf C}$, and a suitable large value $t$ of the parameter in our
asymptotic morphism.

Using the stability of the defining relations for $\OA{m}$
(see the proof of Lemma 2.1 (5)), choose $\eta > 0$ such that
whenever $A$ is a \CA\  and $\sigma_0, \sigma_1 : \OA{m} \to A$
are \hm s such that $\| \sigma_0 (s_j)- \sigma_1 (s_j) \| < \eta$
for $j = 1, \dots, m$, then $\sigma_0$ is homotopic to $\sigma_1$.

Let
 $M_0 = \sup\{\| \varphi_t \|: t \in [0,1] \}.$
Note that $M_0 \geq 1$.
For $\dt > 0$, define, using the functions
$\bt_1, \bt_2, \dots$ from Lemma 2.1:
\begin{eqnarray}
M &  = &  M_0 + \dt,     \nonumber \\
\dt' & =  & \left[ M^2 \bt_3 (\bt_2 (\dt) + \dt) +
       (M^2 + 1) \bt_1 (M^2 \dt) + 2 (M^2 + 1) \dt \right]^{1/2},
          \label{defn2}   \\
{\rm and}  & &   \nonumber  \\
\dt'' & = & \max \left( 2 \dt, \,\,
       \bt_2 \left( \bt_2 (\dt') + (M^2 + 2) \bt_1 (M^2 \dt) +
                                       2 M^2 \dt \right)
                  + \bt_1 (M^2 \dt) \right).   \nonumber  \\
 & & \label{defn3}
\end{eqnarray}
Now choose $\dt>0$ so that:
\begin{eqnarray}
 & & \bt_2 (2 \bt_2 (\dt)) + \bt_2 (\dt) < 1/2,   \label{init1} \\
 & & \bt_4 \left(\dt'' + \bt_2 (\dt') + \bt_1 (M^2 \dt) \right)
            < \infty,   \label{init2} \\
 & &  4 M^3 (m + 1) \bt_4 (\dt'') + 2 \bt_4 (2 M^3 \dt)
     + 10 M^3 (m + 1) \dt < \ep,   \label{init3} \\
{\rm and}  & &   \nonumber  \\
 & &  \bt_5 (\dt + \bt_4 (\dt'')) < \min (\eta, \ep).   \label{init4}
\end{eqnarray}

Define $T = \{\zt-1: \zt \in S^1 - \{1\} \}$. (Thus, if $v$ is unitary,
then ${\rm sp}(v - 1) \subset T \cup \{0\}.$)
Choose functions
 $g, g', h, h' \in C_c (T)$
satisfying the following properties:
\begin{eqnarray*}
 & & 0\le g, g' \le 1,\,\,\,\,\,\, g' g = g,  \\
 & & h(T), h' (T) \subset T \cup \{0\} \andeqn
        |h  (\zt)-\zt|, |h' (\zt)-\zt| <\delta,  \\
 & & g h = h, \,\,\,\,\,\, g' h' = h',  \andeqn
     h' (\zt)=\zt \,\,\,{\rm on}\,\,\, {\rm supp}(g).
\end{eqnarray*}
Note that $g'$ and $h'$ are not the derivatives of $g$ and $h$.

Choose, and fix, $t$ so large that the estimates
\beq
\left\| \varphi_t \left(\prod_{i = 1}^l x_i \right) -
      \prod_{i = 1}^l  \varphi_t ( x_i) \right\| < \dt   \label{product}
\eeq
hold for $1 \leq l \leq 7$ and $x_1, \dots, x_l$ of the form
$a \otimes b$, with $a$ or $a^*$ in
\[
\{ g(u - 1), \,  g(u - 1)^{1/3}, \, g'(u - 1), \,
   g'(u - 1)^{1/3}, \, h(u - 1), \, h'(u - 1) \}
\]
and with $b$ or $b^*$ in $\{ 1, s_j, p_j \}$.
(To get such an estimate for products of three or more elements, we use
$\| \varphi_t \| \leq M_0$ for all $t$. Note that our hypotheses imply
$\varphi_t (x) \varphi_t (y) = \varphi_t (xy)$ for $x,\, y$ both
of the form $a \otimes b,$ with $a$ as above and $b \in \{1, p_j\}$.)

{\em Step 1.3:} We now choose elements of $C([0,1],B)$ to be used to
build the partial isometries and matrix units that we will need.

Define the following elements of $C([0,1],B):$
\begin{eqnarray*}
a_j & = & \varphi_t(h  (u-1)\otimes p_j),   \\
w_j & = &  \varphi_t (g   (u-1)^{1/3} \otimes p_j)
        \varphi_t(g  (u-1)^{1/3} \otimes s_j)
        \varphi_t (g   (u-1)^{1/3} \otimes 1) ,  \\
\tilde{w}_{ij} & = & \varphi_t(g (u-1)\otimes s_i s_j^*).
\end{eqnarray*}
Further define $a_j'$, $w_j'$, and $\tilde{w}_{ij}'$ in the same way,
but using $g'(u-1)$ and $h'(u-1)$ in place of $g(u-1)$ and $h(u-1)$.
The first and last factors in the definition of $w_j$ and $w_j'$ ensure
that, for example, $w_j w_j^*$ is in the hereditary subalgebra generated
by $\varphi_t (g(u - 1)) \otimes p_j).$ The powers of $g(u - 1)$ and
 $g'(u - 1)$ are chosen to as to have $\tilde{w}_{ij}$
 close to $w_i w_j^*$
and $\tilde{w}_{ij}'$ close to $w_i' (w_j')^*$.

The estimate (\ref{product}) implies relations of the form
\begin{eqnarray}
  & &  \|a_j w_j - w_j \left(\sum_{i=1}^m a_i \right)\| <2\dt,
      \nonumber  \\
  & & \|\tilde{w}_{ij}' w_j - w_i\| < 2\dt ,\,\,\,\,\,\,
      \|w_i w_j^* - \tilde{w}_{ik} \tilde{w}_{kj}^*\| < 2\dt,
       \,\,\,\,\,\, {\rm etc.,}    \label{stdest}
\end{eqnarray}
since both terms in each difference differ by less than $\dt$ from
$\varphi(x)$ for some $x.$ Among
$\tilde{w}_{ij},\tilde{w}_{ij}'$, $a_j$, and $a_j',$
 one actually gets
equalities at these places, for example:
\beq
a_i\tilde{w}_{ij}=\tilde{w}_{ij}a_j, \,\,\,\,\,\,
\tilde{w}_{ij}\tilde{w}_{kl}=\dt_{jk}\tilde{w}_{ii}\tilde{w}_{il}
       \,\,\,\,\,\, {\rm etc.}    \label{stdeq}
\eeq
Also note that
\[
\tilde{w}_{ji}=\tilde{w}_{ij}^* \andeqn
            \tilde{w}_{ji}'=(\tilde{w}_{ij}')^*.
\]
Since $\| \varphi_t \| \leq M_0,$ we have $\|a_j\|\leq 2M_0 \leq 2M$ for
each $j,$ and all the other elements listed above have norm at most
either  $M_0 \leq M$ or $M_0 + \dt = M$.

{\em Step 1.4:} We define two important hereditary subalgebras
$D$ and $D'$, which will play a crucial role in the proof.

Let $D$ be the hereditary $C^*$-subalgebra of $C([0,1],B)$ generated
by $a_1.$ Then for each $\alpha\in [0,1],$ the subalgebra
$D_{\alpha}\i B$ (defined as in Notation 2.2) is the hereditary
$C^*$-subalgebra generated by
\[
{\rm ev}_{\alpha}(a_1)={\rm ev}_{\alpha}(\varphi_t (h(u-1)\otimes p_1))=
h({\rm ev}_{\alpha}\circ\varphi_t((u-1)\otimes p_1)).
\]
(This uses the fact that $\varphi_t|_{C_0(S^1 - \{1\})\otimes M_m}$
is a homomorphism.)
Since ${\rm ev}_{\alpha}\circ\varphi_t((u-1)\otimes p_1)$ has full
spectrum by assumption (full here means equal to
 $T \cup \{0\} = \{\zt-1:\zt\in S^1\},$
the largest it can be), and since the range of $h$ is
$T \cup \{0\},$ the element
 ${\rm ev}_{\alpha}(a_1)$ also has full spectrum.
It follows that each $D_{\alpha}$ is nonzero and nonunital.
Therefore Proposition 2.6 and Corollary 2.8 imply that $D$
has an increasing approximate identity consisting of projections
whose classes in $K_0 (D)$ are trivial.

The same reasoning applies to the hereditary subalgebra $D'$ generated
by $a_1'.$ Note that $D\i D'.$

Since $a_1$ and $a_1'$ are normal, we have
$$
D=\overline{ a_1C([0,1],B)a_1^*}=\overline{ a_1^*C([0,1],B)a_1}
$$
and similarly for $a_1'$ and $D'.$

We observe the following important property
 of the hereditary subalgebras
$D$ and $D':$

(M) If $x\in D$, then $\tilde{w}_{11}x=x\tilde{w}_{11}=x.$

(M$'$) If $x\in D',$ then  $\tilde{w}_{11}'x=x\tilde{w}_{11}'=x.$

Property (M) follows by a standard argument from the relations
$\tilde{w}_{11} a_1 = \tilde{w}_{11} a_1 = a_1$
and $\tilde{w}_{11} a_1^* = a_1^* \tilde{w}_{11} = a_1^*.$
Property (M$'$) is similar. Similar arguments also give the following
properties:

(Z) If $x\in D$, then $\tilde{w}_{ij} x = x \tilde{w}_{ji} = x w_j = 0$
for arbitrary $i$ and for $j \neq 1$.

(Z$'$) If $x\in D'$, then
$\tilde{w}_{ij}' x = x \tilde{w}_{ji}' = x w_j' = 0$
for arbitrary $i$ and for $j \neq 1$.

{\em Step 2:} We now construct a first (lower) level of matrix units
$\{e_{ij}\}$ in $D$. These look like they should span the copy of
$M_m$ inside a homomorphic image of $\OA{m}$, but unfortunately things
are not that easy.

Choose a projection $e_1\in D$ (from the approximate
identity obtained above) such that $[e_1]=0$ in
$K_0(D)$ and
\beq
\|e_1 a_1 - a_1\|, \, \|a_1 e_1 - a_1\| < \dt.   \label{e1}
\eeq
Define
\[
e_{ij}=\tilde{w}_{i1}e_1\tilde{w}_{1j},\,\,\,\,\,\,\,
       e_j=e_{jj},\andeqn e=\sum_{j=1}^me_j.
\]
We claim that $(e_{ij})_{i,j=1}^m$ is a system of matrix units in
$C([0,1],B).$ To prove this,
 we first note that clearly $e_{ij}^*=e_{ji}.$
Furthermore, using the properties (M) and (Z) from Step 1.4,
\[
e_{ij} e_{kl} =
  \tilde{w}_{i1} e_1  \tilde{w}_{1j}  \tilde{w}_{k1} e_1  \tilde{w}_{1l}
    = \dt_{jk}  \tilde{w}_{i1} e_1  \tilde{w}_{11}^2 e_1  \tilde{w}_{1l}
     = \dt_{jk} e_{il}.
\]

{\em Step 3:} We now construct a second level of matrix units
$\{e_{ij}'\}$ in $D'$. They need to be enough bigger than the first
ones to allow room for a correction for the failure of
$w_j^* e_j w_j$ to be close to $e$. To ensure this,
we first construct projections $r$, $f$, and $e_1'$ in $D'$ such that
$\| f - w_1 e w_1^* \|$ is small, $e_1, f \leq r$, and
$r < e_1'$. We will need $r$ in the next step, but all the other
intermediate projections constructed in this step can be discarded
after it is done.

The definition of $w_1$ implies that $w_1 e w_1^* \in D'.$ Furthermore,
we have the following estimate, in which the second step uses
properties (M)
and (Z), and the last step is similar to (\ref{stdest}):
\begin{eqnarray*}
\lefteqn{
\| (w_1 e w_1^*)^2 - w_1 e w_1^* \| \leq M^2 \| e w_1^* w_1 e -e \|
}   \\
 & &
= M^2 \left\| e w_1^* w_1 e -
              e\left( \sum_j  \tilde{w}_{jj}^2 \right) e \right\|
   \leq
M^2 \left\| w_1^* w_1 - \left( \sum_j  \tilde{w}_{jj}^2 \right) \right\|
           < M^2 \dt.
\end{eqnarray*}
Therefore there exists a projection $f_0 \in D'$ such that
\[
\| f_0 - w_1 e w_1^* \| < \beta_1 (M^2 \dt).
\]

Since $D'$ has an increasing approximate identity
of projections with trivial $K_0$-classes, we can
choose a projection $r_0\in D'$ such that
$[r_0] = 0$ in $K_0 (D')$ and
\beq
\|r_0 e_1 - e_1\| < \dt \andeqn \|r_0 f_0 - f_0\| < \dt,
          \label{rf-f}
\eeq
and then choose a projection $r_1 \in D'$ such that $r_1 > r_0$
and $[r_1] = 0$ in $K_0 (D')$. It follows that there is
a projection $r\in D'$ such that $r\ge e_1$ and
\beq
\|r-r_0\|<\beta_2(\dt). \label{r-r0}
\eeq
We then have  $\| r f_0 - f_0 \| < \bt_2 (\dt) + \dt.$ Therefore there
is a projection $f \leq r$ such that
$\| f_0 - f \| < \bt_3 (\bt_2 (\dt) + \dt),$ and it follows that
\beq
\| f - w_1 e  w_1^* \| <
    \bt_1 (M^2 \dt) + \bt_3 (\bt_2 (\dt) + \dt).   \label{f}
\eeq
We also have $\| r_1 r - r \| \leq 2 \| r - r_0 \| < 2 \bt_2 (\dt).$
Using Lemma 2.1 (2) and combining the resulting estimates with ones we
already have (including (\ref{init1})), we obtain a projection
$e_1' \in D'$ such that $e_1' \geq r$ and
\[
\| (e_1' - r) - (r_1 - r_0) \| < \bt_2 (2 \bt_2 (\dt)) + \bt_2 (\dt)
          < 1/2.
\]
Since $r_1 > r_0$, it follows that $e_1' > r.$ Similar arguments show
that the classes of $f$, $r$, and $e_1'$ in $K_0 (D')$ are all zero.

We now define
\[
e_{ij}'= \tilde{w}_{i1}'e_1' \tilde{w}_{1j}', \,\,\,\,\,\,
          e_j'=e_{jj}',\andeqn e'=\sum_je_j'.
\]
Then $\{e_{ij}'\}$ is a system of matrix units by the same argument
as for $\{e_{ij}\}.$ (See Step 2.)
 We have $e_1'\ge e_1$ by construction, so
\[
e_j e_j' =
 \tilde{w}_{j1} e_1 \tilde{w}_{1j} \tilde{w}_{j1}' e_1' \tilde{w}_{1j}'
 = \tilde{w}_{j1} e_1 e_1' \tilde{w}_{1j}' = e_j,
\]
that is, $e_j \ge e_j'.$ It follows that $e' \geq e.$
Since we actually have $e_1' > e_1$, we in fact
get $e' > e.$ Furthermore, similar
arguments show
\[
e_ie_{ij}'=e_{ij}=e_{ij}'e_j.
\]

{\em Step 4:}
We now
construct a projection  $p > e$ such that
$[p] = 0$ in $K_0 (C([0,1], B))$ and
$\| p - (w_1')^* e_1' w_1' \|$ is small. In order
to ensure that $p > e$, we use an intermediate projection
$p_0 \geq e$ such that
$\| p_0 - (w_1')^* r w_1' \|$ is small.

We estimate, using reasoning similar to that for (\ref{stdest}):
\begin{eqnarray*}
\lefteqn{
\|((w_1')^* r w_1')^2 - (w_1')^* r w_1'\|
  \le  M^2 \|r w_1' (w_1')^* r - r\|}   \\
 & & = M^2 \|r w_1' (w_1')^* r - r ( \tilde{w}_{11}')^2 r\|
       \le M^2 \|w_1' (w_1')^* - ( \tilde{w}_{11}')^2\| < M^2 \dt.
\end{eqnarray*}
Therefore there is a projection $p_1$ such that
\[
\| p_1 - (w_1')^* r w_1' \| < \bt_1 (M^2 \dt).
\]
Similarly, there is a projection $p_2$ such that
\[
\| p_2 - (w_1')^* e_1' w_1' \| < \bt_1 (M^2 \dt).
\]
Using $f \leq r$ in the third step and (\ref{f}) in the fourth
step, we estimate:
\begin{eqnarray}
\lefteqn{
\| p_1 e - e \|^2   =   \| e (1 - p_1) e \|
     \leq  \|p_1  - (w_1')^* r w_1' \| + \|e - e (w_1')^* r w_1' e \|
             \nonumber   }   \\
  & < & (1 + M^2)  \bt_1 (M^2 \dt) + \| e (1 - (w_1')^* f w_1') e \|
                \nonumber  \\
 & < & (1 + M^2)  \bt_1 (M^2 \dt) +  M^2 \bt_3 (\bt_2 (\dt) + \dt) +
           \| e (1 - (w_1')^* w_1 e w_1^* w_1') e \|  \nonumber \\
 & &               \label{pe}
\end{eqnarray}
To estimate the last term, we first observe that
\[
e \sum_j \varphi_t (g(u-1)^{1/3} \otimes p_j) =
      \sum_j e_j (\tilde{w}_{jj})^{1/3} = e =
            e \sum_j \varphi_t (g'(u-1)^{1/3} \otimes p_j).
\]
Therefore $e w_1^* = e (w_1')^*.$ So we can replace $w_1$
by $w_1'$ in the last term of (\ref{pe}). Now
\begin{eqnarray*}
\|e (w_1')^* w_1' - e\| & = &
\|e \left(\sum_j \tilde{w}_{jj} \right) (w_1')^* w_1' -
                          e \left(\sum_j \tilde{w}_{jj} \right)\|
  \\
 & \le & \|e\| \| \left(\sum_j \tilde{w}_{jj} \right) (w_1')^* w_1'-
             \left(\sum_j \tilde{w}_{jj} \right)\|.
\end{eqnarray*}
We have $\sum_j \tilde{w}_{jj} = \varphi_t (g_1 (u-1) \otimes 1)$.
Furthermore,
$\left(\sum_j \tilde{w}_{jj} \right) (w_1')^* w_1'$
is the product of 7 factors
$\varphi_t (x_1) \cdots \varphi_t (x_7)$, with the $x_i$ in the list of
elements to which (\ref{product}) applies, and
$x_1 \cdots x_7 = g_1 (u - 1) \otimes 1$. Therefore
\[
\|e (w_1')^* w_1' - e \| \leq 2 \dt.
\]
Applying this inequality and its adjoint, we can replace the middle
factor in the last term of (\ref{pe}) by $1 - e$ at a cost of
 $2 (M^2 + 1) \dt.$
Therefore, using (\ref{defn2}), we get
\[
\| p_1 e - p_1 \| <
     [M^2 \bt_3 (\bt_2 (\dt) + \dt) + (M^2 + 1)  \bt_1 (M^2 \dt) +
                2 (M^2 + 1) \dt]^{1/2} = \dt'.
\]
It follows that there is a projection $p_0 \geq e$ such that
$\| p_0 - p_1 \| < \bt_2 (\dt'),$ whence
\beq
\| p_0 - (w_1')^* r w_1' \| <  \bt_2 (\dt') + \bt_1 (M^2 \dt).
            \label{p0}
\eeq
We now have
\beqr
\lefteqn{ \| p_2 p_0 - p_0\|  \leq  2 \|p_0 - (w_1')^* r w_1'\|   }  \\
 & & \mbox{} + \|(w_1')^* r w_1'\| \|p_2 - (w_1')^* e_1' w_1'\|
         + \|(w_1')^* e_1' w_1' (w_1')^* r w_1' - (w_1')^* r w_1'\|.
\eeqr
The reasoning of (\ref{stdest}) shows that
$\| w_1' (w_1')^* - (\tilde{w}_{11}')^2 \| < 2 \dt;$ furthermore,
$e_1' (\tilde{w}_{11}')^2 r = e_1' r = r$. So the last term above has
norm at most $2 M^2 \dt$, and we get
\[
\| p_2 p_0 - p_0\| <
   2 (\bt_2 (\dt') + \bt_1 (M^2 \dt)) + M^2 \bt_1 (M^2 \dt) + 2 M^2 \dt.
\]
Applying Lemma 2.1 (2) in the usual way, we find a projection
$p \geq p_0$ such that
\begin{eqnarray}
\lefteqn{\| p - (w_1')^* e_1' w_1' \| \leq
        \| p - p_2 \| + \| p_2 - (w_1')^* e_1' w_1' \|} \nonumber  \\
   &  & < \bt_2 \left( \bt_2 (\dt') + (2 + M^2) \bt_1 (M^2 \dt)
                           + 2 M^2 \dt\right) + \bt_1 (M^2 \dt)
          \leq \dt''.           \label{p}
\end{eqnarray}
(See (\ref{defn3}) for the definition of $\dt''$.) We note that
\[
\|(p - p_0') - (w_1')^* (e_1' - r) w_1'\| <
   \dt'' + \bt_2 (\dt') + \bt_1 (M^2 \dt)
\]
by (\ref{p0}) and (\ref{p}), and
$\| w_1' (w_1')^* - (\tilde{w}_{11}')^2 \| < 2 \dt$ implies
(using (\ref{defn3}))
\[
\|(e_1' - r) w_1' (w_1')^* (e_1' - r) - (e_1' - r)\| < 2 \dt
       \leq \dt'' + \bt_2 (\dt') + \bt_1 (M^2 \dt).
\]
Since $\bt_4 (\dt'' + \bt_2 (\dt') + \bt_1 (M^2 \dt)) < \infty$
by (\ref{init2}), Lemma 2.1 (4) implies that $p - p_0'$
is Murray-von Neumann equivalent to $e_1' - r,$ and therefore
nonzero.
Similar estimates show $p$ is Murray-von Neumann equivalent to $e_1'$
and $p_0'$ is Murray-von Neumann equivalent to $r.$
It follows that $[p] = [p_0'] = 0$ in $K_0 (C([0,1], B)).$

{\em Step 5:} We now construct a \hm\
$\rho_0 : \OA{m} \to C([0,1]) \otimes B.$ This \hm\  will be one
direct summand in the \hm\  we need for the proof of the lemma.

Since $p > p_0'$ and $p_0' \ge e,$
 we get $p > e;$ also recall from above that
$e' > e.$ Since all $K_0$-classes are trivial,
Lemma 2.7 implies
that $p - e$ and $e' - e$ are Murray-von Neumann equivalent. Thus, there
exists a partial isometry $d \in C([0,1]) \otimes B$ such that
\[
d^* d = e', \,\,\,\,\,\, d d^* = p,  \andeqn d e = e d = e.
\]
Furthermore, the estimates
\[
\|p - (w_1')^* e_1' w_1' \| < \dt''
\andeqn  \|e_1' w_1' (w_1')^* e_1' - e_1' \| < 2 \dt \leq \dt''
\]
imply the existence of a partial isometry $\overline{w}_{1}$ such that
\beq
\overline{w}_{1} \overline{w}_{1}^* = e_1', \,\,\,\,\,\,
   \overline{w}_{1}^* \overline{w}_{1} = p, \andeqn
\|e_1' w_1' - \overline{w}_{1}\| < \bt_4 (\dt'').
   \label{wbar}
\eeq
Note that
\[
(\overline{w}_{1} d)^* (\overline{w}_{1} d) = e'
       \andeqn (\overline{w}_{1} d )(\overline{w}_{1} d)^* = e_1'.
\]

Define
\[
z_j = e_{j1}' \overline{w}_{1} d.
\]
Then
\[
z_j^* z_j = d^* \overline{w}_{1}^* e_{1j}' e_{j1}' \overline{w}_{1} d
         = d^* \overline{w}_{1}^* e_1' \overline{w}_{1} d = e'
\]
and
\[
z_j z_j^* = e_{j1}' \overline{w}_{1} d d^* \overline{w}_{1}^* e_{1j}'
             = e_{j1}' e_1' e_{1j}' = e_j'.
\]
Thus, there exists a homomorphism $\rho_0: \OA{m}\to C([0,1],B)$
such that $\rho_0 (s_j) = z_j$ and $\rho_0 (1) = e'.$

{\em Step 6:}
Next, we construct a unitary $v_0 \in e'C([0,1],B)e'$
which approximately commutes with the range of $\rho_0.$
For later use, we actually do the following more general calculation.
Let $\tilde{d}$ be a partial isometry and let $\tilde{e}'$ be a
projection such that:
\beq
\tilde{d} \tilde{d}^* \geq p, \,\,\,\,\,\,
     \tilde{d} e = e \tilde{d} = e, \,\,\,\,\,\,
      \tilde{e}' \geq \tilde{d}^* \tilde{d},  \andeqn
      \tilde{e}' \geq  e'. \label{tilde}
\eeq
Set
\beq
\tilde{c} = \tilde{e}' + \sum_{j = 1}^m e_j a_j e_j.  \label{ctilde}
\eeq
Define $\tilde{z}_j = e_{j1}' \overline{w}_1 \tilde{d}$. (That is, use
$\tilde{d}$ in place of $d$ in the definition of $z_j$.)
Then $\tilde{c}$ is close to a unitary
$\tilde{v}_0 = \tilde{c} (\tilde{c}^* \tilde{c})^{-1/2}
        \in \tilde{e}' C([0,1],B) \tilde{e}'$
which approximately
commutes with $\tilde{z}_j$ for all $j$.

For the special case $\tilde{d} = d$ of immediate interest, we set
$\tilde{c} = c$ and $\tilde{v}_0 = v_0.$

A number of the intermediate estimates will be needed in Step 7.

{\em Step 6.1:}
We show that $\tilde{c}$ is approximately unitary.

First, observe that
\beqr
\|a_j e_j - a_j\| & = &
    \| \tilde{w}_{j1} a_1 \tilde{w}_{1j} e_j -
                                 \tilde{w}_{j1} a_1 \tilde{w}_{1j}\|
= \| \tilde{w}_{j1} a_1 e_1 \tilde{w}_{1j} -
                                 \tilde{w}_{j1} a_1 \tilde{w}_{1j}\|
                 \\
 & \le &  M^2 \|a_1 e_1 - a_1\| < M^2 \dt.
\eeqr
Since $\|a_j^*\| \le 2M$, we get
$\|a_j^* e_j a_j - a_j^* a_j\| < 2 M^3 \dt. $
Combining this with the fact that the $e_j$ are orthogonal projections
dominated  by $\tilde{e}'$  and the fact (similar to (\ref{stdeq})) that
$a_j^* a_j - a_j - a_j^* = 0$, we obtain:
\begin{eqnarray}
\|\tilde{c}^* \tilde{c} - \tilde{e}'\| & = &
             \| \sum_j e_j (a_j^* e_j a_j - a_j - a_j^*) e_j\|
               \nonumber    \\
    &  \le & \max_j \|a_j^* e_j a_j - a_j - a_j^*\| <
     2 M^3 \dt.    \label{c*c}
\end{eqnarray}
A similar estimate shows that
\beq
\|\tilde{c} \tilde{c}^* - \tilde{e}'\|<2 M^3 \dt.    \label{cc*}
\eeq

{\em Step 6.2:} We prove that $\tilde{c}$ commutes with $e_{j1}'$, the
first factor in the definition of $\tilde{z}_j.$

First, certainly $\tilde{e}' \geq e' = \sum_i e_{ii}'$ commutes with
$e_{j1}'$.

For the other part, we start with the equation
$\tilde{w}_{1j} a_j \tilde{w}_{j1} = \tilde{w}_{11} a_1 \tilde{w}_{11}$,
which follows from the same reasoning as (\ref{stdeq}). Therefore
\[
e_j a_j e_{j1} =
(\tilde{w}_{j1} e_1 \tilde{w}_{1j}) a_j
                            (\tilde{w}_{j1} e_1 \tilde{w}_{11})
 = (\tilde{w}_{j1} e_1 \tilde{w}_{11}) a_1
                            (\tilde{w}_{11} e_1 \tilde{w}_{11})
 = e_{j1} a_1 e_{1}.
\]
Now one checks that
\[
\left( \sum_{i = 1}^m e_i a_i e_i \right) e_{j1}' = e_j a_j e_{j1}
          \andeqn
e_{j1}' \left( \sum_{i = 1}^m e_i a_i e_i \right) = e_{j1} a_1 e_{1}.
\]
So $\tilde{c}$ commutes with $e_{j1}'$.

{\em Step 6.3:} We next show that $\tilde{c}$ approximately commutes
with $\overline{w}_1 \tilde{d}$, the other factor in
the definition of $\tilde{z}_j.$

Again, $\tilde{e}'$ is easy. We have $\tilde{e}' \geq e' \geq e_1'$, so
$\tilde{e}' \overline{w}_1 \tilde{d} = \overline{w}_1 \tilde{d}$,
and also
 $\tilde{e}' \geq \tilde{d}^* \tilde{d}$, so
$\overline{w}_1 \tilde{d} \tilde{e}' = \overline{w}_1 \tilde{d}$.
So we get exact commutativity here:
\beq
\tilde{e}' \overline{w}_1 \tilde{d} -
         \overline{w}_1 \tilde{d} \tilde{e}' = 0.    \label{comm1}
\eeq

The other part is longer. We begin by using the relations
$e_j \overline{w}_1 = e_j e_1' \overline{w}_1 =
      \dt_{1j} e_1 \overline{w}_1$
and $\tilde{d} e_j = e_j \tilde{d} = e_j$,
to get
\beqr
\tilde{c} \overline{w}_1 \tilde{d} - \overline{w}_1 \tilde{d} \tilde{c}
           & = &
  \left( \sum_{j = 1}^m e_j a_j e_j \right) \overline{w}_1 \tilde{d}
   - \overline{w}_1 \tilde{d} \left( \sum_{j = 1}^m e_j a_j e_j \right)
         \\
 & = &  e_1 a_1 e_1 \overline{w}_1 \tilde{d} -
     \overline{w}_1 \left( \sum_{j = 1}^m e_j a_j e_j \right) \tilde{d}.
\eeqr
Using $\| \overline{w}_1 - e_1' w_1' \| < \bt_4 (\dt'')$,
$e_1 e_1' = e_1$, and $\| a_j \| \leq 2 M$, we get
\begin{eqnarray}
\| \tilde{c} \overline{w}_1 \tilde{d} -
                            \overline{w}_1 \tilde{d} \tilde{c}  \|
 &  < &  2 M (m + 1) \bt_4 (\dt'')   \nonumber  \\
 & & \mbox{} +   \| e_1 a_1 e_1 w_1' \tilde{d} -
    e_1' w_1' \left( \sum_{j = 1}^m e_j a_j e_j \right) \tilde{d} \|.
          \label{mess}
\end{eqnarray}
Now (\ref{e1}) implies $\| e_1 a_1 e_1 - a_1 \| < 2 \dt,$ and
\[
e_j a_j e_{j} - a_j =
(\tilde{w}_{j1} e_1 \tilde{w}_{1j}) a_j
                            (\tilde{w}_{j1} e_1 \tilde{w}_{1j})
     - \tilde{w}_{j1} a_1 \tilde{w}_{1j}
= \tilde{w}_{j1} (e_1 a_1 e_1 - a_1) \tilde{w}_{1j},
\]
so
\beq
\| e_j a_j e_{j} - a_j \| < 2 M^2 \dt.  \label{ejajej}
\eeq
We apply this to each of the $m + 1$ $a$'s in (\ref{mess}), and use
$\| \tilde{d} \| \leq 1$ and $\| w_1' \| \leq M,$ to get:
\beqr
\| \tilde{c} \overline{w}_1 \tilde{d} -
                        \overline{w}_1 \tilde{d} \tilde{c}  \|
 &  < &  2 M (m + 1) \bt_4 (\dt'') + 2 M^3 (m + 1) \dt  \nonumber  \\
 &  & \mbox{} + \| a_1 - e_1' a_1 \| \|w_1'\|
    + \| e_1' \| \| a_1 w_1' - w_1' \sum_j a_j \|.
\eeqr
Now
\[
\| a_1 w_1' - w_1' \sum_j a_j \|  =
      \| a_1 w_1' - \sum_j w_1' a_j \| < (m + 1) \dt,
\]
by reasoning similar to that which gave (\ref{stdest}). (One needs to
use (\ref{product}) on products of four factors.) Furthermore,
\beq
\| a_1 - e_1' a_1 \| = \|a_1 (1 - e_1') a_1\|^{1/2}
     \leq \|a_1 (1 - e_1) a_1\|^{1/2} < \dt,  \label{e1'a1}
\eeq
by (\ref{e1}). So
\begin{eqnarray}
\| \tilde{c} \overline{w}_1 \tilde{d} -
                      \overline{w}_1 \tilde{d} \tilde{c}  \|
 &  < &  2 M (m + 1) \bt_4 (\dt'') + 2 M^3 (m + 1) \dt
                 + M \dt + (m + 1) \dt  \nonumber  \\
 & \leq & 2 M^3 (m + 1) (\bt_4 (\dt'') + 2 \dt).       \label{comm2}
\end{eqnarray}

{\em Step 6.4:} We now combine our estimates to get an estimate on
the commutator of $\tilde{z}_j$ with
the unitary $\tilde{c} (\tilde{c}^* \tilde{c})^{-1/2}$.

Lemma 2.1 (4) and the estimates (\ref{c*c}) and (\ref{cc*}) show that
the unitary
$\tilde{v}_0  = \tilde{c} (\tilde{c}^* \tilde{c})^{-1/2}
      \in \tilde{e}'C([0,1],B)\tilde{e}'$
satisfies
\beq
\| \tilde{v}_0 - \tilde{c} \| < \bt_4 (2 M^3 \dt).  \label{v0}
\eeq
It is also important to notice that
\beq
\tilde{v}_0 (\tilde{e}' - e) = (\tilde{e}' - e) \tilde{v}_0
      = \tilde{e}' - e   \label{v0e'}.
\eeq
This follows from the analogous fact for $\tilde{c}$.

{}From (\ref{comm1}), (\ref{comm2}), and the definition of
$\tilde{z}_j$, we get
\beq
\| \tilde{c} \tilde{z}_j - \tilde{z}_j \tilde{c} \| <
        2 M^3 (m + 1) (\bt_4 (\dt'') + 2 \dt).  \label{comm3}
\eeq
Therefore
\beq
\| \tilde{v}_0 \tilde{z}_j - \tilde{z}_j \tilde{v}_0 \| <
        2 M^3 (m + 1) (\bt_4 (\dt'') + 2 \dt) + 2 \bt_4 (2 M^3 \dt)
               < \ep.  \label{v0zj}
\eeq
(The last step is a weaker relation than (\ref{init3}).)

{\em Step 7:} Continuing with the hypotheses and notation of Step 6,
we verify the following estimates, which are analogous to those in
the conclusion (C2):
\[
\|(\tilde{v}_0-\tilde{e}')-\vph_t((u-1)\otimes 1)\|,\,\,
       \|(\tilde{v}_0-\tilde{e}')\tilde{z}_j-\vph_t((u-1)\otimes s_j)\|
< \ep.
\]

For the first, we observe
\beqr
\lefteqn{
\|\tilde{v}_0 - \tilde{e}' - \vph ((u-1) \otimes 1)\|
        }   \\
 & \le &  \|\tilde{v}_0 - \tilde{c}\| +
                  \|\tilde{c} - \tilde{e}' - \sum_j a_j\|
+ \|\vph_t (h(u-1) \otimes 1) - \vph_t ( (u - 1) \otimes 1)\|
\\
 & < & \bt_4 (2 M^3 \dt)  + \sum_j \|e_j a_j e_j - a_j\| +
                            \| \varphi_t \|  \|h(u - 1) - (u - 1)\|
\\
 & <  & \bt_4 (2 M^3 \dt) + 2 M^2 m \dt + M \dt < \ep.
\eeqr
(The last step is a weaker relation than (\ref{init3}).)

For the second estimate, we do the following approximations,
in which we give the justifications and errors afterwards:
\beqr
(\tilde{v}_0 - \tilde{e}') \tilde{z}_j &
   \approx & (\tilde{c} - \tilde{e}') \tilde{z}_j
   \approx  \tilde{z}_j (\tilde{c} - \tilde{e}')
   = e_{j1}' \overline{w}_1 \tilde{d} \sum_i e_i a_i e_i    \\
 & \approx &  e_{j1}' w_1' \sum_i e_i a_i e_i
   \approx  e_{j1}' w_1' \sum_i a_i         \\
 & = &  \tilde{w}_{j1}' e_{1}' \tilde{w}_{11}' w_1' \sum_i a_i
   \approx  \tilde{w}_{j1}' e_{1}' a_1 \tilde{w}_{11}' w_1'
   \approx \tilde{w}_{j1}'  a_1 \tilde{w}_{11}' w_1'  \\
 & \approx &  \varphi_t (h(u - 1) \otimes s_j)
   \approx \varphi_t ((u - 1) \otimes s_j).
\eeqr
The errors are, in order (including the equality signs),
$\bt_4 (2 M^3 \dt)$ (by (\ref{v0})),
$ 2 M^3 (m + 1) (\bt_4 (\dt'') + 2 \dt)$\,
    (by (\ref{comm3})),
0 \,\,(by \,the\, definitions\, of\, $\tilde{z}_j$\, and\, $\tilde{c}$),\\
$2 M m \bt_4 (\dt'')$ (by (\ref{wbar}) and because $\tilde{d} e_i = e_i$
  and $e_{j1}' e_1' = e_{j1}'$),
$M m (2 M^2 \dt)$ (by (\ref{ejajej})),
0,
$M (m + 1) \dt$ (by reasoning similar to (\ref{stdest}), using
    (\ref{product}) on products of 5 factors),
$M^3 \dt$ (by (\ref{e1'a1})),
$\dt$ (by (\ref{product}), applied to a product of 6 factors),
and $M \dt$ (since $\| \varphi_t \| \leq M$ and by the choice of $h$).
Adding up the errors, and replacing the result by a larger but
simpler expression, we find that
\begin{eqnarray}
\lefteqn{ \|(\tilde{v}_0 - \tilde{e}') \tilde{z}_j -
      \varphi_t ((u - 1) \otimes s_j) \| }   \nonumber   \\
 & & <
  4 M^3 (m + 1) \bt_4 (\dt'') + \bt_4 (2 M^3 \dt)
         + 10 M^3 (m + 1) \dt < \ep,
\end{eqnarray}
by (\ref{init3}).

{\em Step 8:} In this step, we fill out $\rho_0$ and $v_0$ to produce
the homomorphism and unitary we actually need. Because we have worked
in a ``large'' hereditary subalgebra, we can extend $v_0$ by taking it
to be 1 outside $e' C([0,1], B)e'$, and make fairly arbitrary
definitions of partial isometries, without destroying the estimates
we obtained in the previous step. We need this flexibility to match
things up properly at the endpoints of our homotopy.

Recall our initial assumption that
 ${\rm ev}_{\alpha}\circ\vph_t$ is a
constant function of $\alpha$ for $\alpha\in [0,1/3]$ and also for
$\alpha\in [2/3,1].$ All our estimates therefore still hold if we
redefine every element of $C([0,1],B)$ that appeared above to
take same value on $[0,1/3]$ that it already does at 1/3, and to
take the same value on $[2/3,1]$ that it already does at 2/3.
We thus assume that everything is constant on $[0, 1/3]$ and also
on $[2/3,1].$

The projections $e'(0)$ and $p(0)$ are both elements of the
hereditary subalgebra of $B$ generated
 by ${\rm ev}_0\circ\vph_t((u-1)\otimes 1),$
since ${\rm ev}_0\circ\vph_t$ is a homomorphism. Since
${\rm ev}_0\circ\vph_t=\psi_0|_{C_0(S^1 - \{1\})\otimes \OA{m}},$ it
follows that
\[
\psi_0(1)>e'(0) \andeqn \psi_0(1)>p(0).
\]
Similarly
\[
\psi_1(1)>e'(1) \andeqn  \psi_1(1)>p(1).
\]
Note that
\[
[\psi_0(1)]=[\psi_1(1)] \,\,\,{\rm in}\,\,\, K_0(B),
\]
by (H5), so
\[
[1-\psi_0(1)]=[1-\psi_1(1)]\andeqn [\psi_0(1)-e'(0)]=[\psi_1(1)-e'(1)].
\]
Furthermore,
$1-\psi_0(1),$ $1-\psi_1(1),$ $\psi_0(1)-e'(0)$
and $\psi_1(1)-e'(1)$ are all nonzero. Also $e'(\alpha)-e(\alpha)$ is
always nonzero. Standard methods thus yield a unitary path
$\alpha\to x (\alpha)$ such that
$x (\alpha)$ is constant on $[0,1/3]$ and on $[2/3,1],$
$x (\alpha)e'(\alpha)x (\alpha)^*$ and $x (\alpha)e(\alpha)x (\alpha)^*$
are constant, and
\[
x(0) \psi_0 (1) x(0)^* =x(1) \psi_1 (1) x(1)^*.
\]
Conjugating everything by this path (and then forgetting it), we may
assume, in addition to everything else, that $\psi_0(1)=\psi_1(1)$
and that $e'$ and $e$ are constant.
Let $q_0$ be the constant projection with value $\psi_0 (1)=\psi_1(1)$.

To construct $\rho,$ we begin by choosing a path
$\alpha \to y(\alpha)$ in $U_0 (q_0(0)Bq_0(0)),$
 defined for $\alpha\in [0,1/3],$
such that
\[
e(\alpha) y(\alpha) = y(\alpha) e(\alpha) = e(\alpha)
    \andeqn y(0) e'(0) = d(0).
\]
(This is possible because the
partial isometry $d(0) - e(0)$ from $e'(0)-e(0)$ to
$p(0)-e(0)$ can be extended to a unitary
 in $U_0 \left([q_0 (0)-e(0)]B[q_0 (0)-e(0)] \right).$)
Next, observe that, at $\alpha=0$ (which we suppress in the notation),
\beqr
\lefteqn{
\|e_j'\psi_0(1\otimes s_j)y e'-z_j\|=
    \| (e_j'  \tilde{w}_{j1}')
                \left( \tilde{w}_{1j}' \psi_0 (1\otimes s_j) d \right) -
     e_j' \overline{w}_1 d\|
                 }           \\
 & = & \| e_{j1}' \varphi_t ( g' (u - 1) \otimes s_1) d -
     e_j' \overline{w}_1 d\|
                            \\
 & \leq &  \| e_{j1}' \| \left[ \| \varphi_t ( g' (u - 1) \otimes s_1) -
                                            w_1' \| +
           \| e_1' w_1' - \overline{w}_1 \| \right] \| d \|
                            \\
 & < & \dt + \bt_4 (\dt'').
\eeqr
This uses, among other things,
$e_j' \tilde{w}_{j1}' = e_{j1}' e_1' = e_{j1}'$ and (\ref{wbar}).
Therefore Lemma 2.1 (5)
 yields a \hm\  $\tau_0: \OA{m} \to B$ such that
 $\tau_0(1)=q_0 (0)-e'(0)$ and
\beq
          \|\tau_0(s_j)+z_j (0) -\psi_0(1\otimes s_j)y(0)\|
<\bt_5(  \dt +  \bt_4 (\dt'')) < \min(\ep, \eta).
        \label{tau}
\eeq
(See (\ref{init4}).)
For $\alpha\in [0,1/3],$ we now define $\tau_{\alpha} = \tau_0$ and
\[
\rho(s_j)(\alpha)=(\tau_0(s_j)+z_j(0))y(\alpha)^* =
    (\tau_{\alpha} (s_j)+z_j({\alpha} ))y(\alpha)^*  .
\]

A similar construction at $\alpha=1$ yields $y(\alpha)$ for
$\alpha\in [2/3,1]$ and a homomorphism $\tau_1:\OA{m}\to B$ such that
$\tau_1(1)= q_0 (1) -e'(1)$ and
\beq
\|\tau_1(s_j)+z_j(1)-\psi_1(1\otimes s_j)y(1)\|<
    \min(\ep, \eta).     \label{tau2}
\eeq
For $\alpha\in [2/3,1]$ we define $\tau_{\alpha} = \tau_1$ and
\[
\rho(s_j)(\alpha)=(\tau_{\alpha} (s_j)+z_j({\alpha} ))y(\alpha)^*.
\]

In $KK^0(\OA{m},B),$ we now have
\[
[{\rm ev}_0\circ\rho_0]=[{\rm ev}_1\circ \rho_0].
\]
Furthermore, the estimates (\ref{tau}) and (\ref{tau2}), and
the choice of $\eta$ at the beginning of the proof, imply that
for $i = 0, 1$ we have
$[{\rm ev}_i \circ\rho_i ]=[ \psi_i  |_{{\bf C} \otimes \OA{m}}]$
in $KK^0(\OA{m},B).$ Therefore
\beqr
[\tau_{1/3}]+[{\rm ev}_{1/3}\circ\rho_0] & = &
        [{\rm ev}_{1/3}\circ\rho]=
        [{\rm ev}_0 \circ\rho]=
        [\psi_{0}|_{{\bf C}\otimes \OA{m}}]
              \\
 & = & [\psi_{1}|_{{\bf C}\otimes \OA{m}}]
   =[{\rm ev}_{1}\circ\rho]=[\tau_{2/3}]+[{\rm ev}_{2/3}\circ\rho_0].
\eeqr
Thus $[\tau_{1/3}]=[\tau_{2/3}].$ Since $m$ is even,
$(1 -e'(0))B( 1 - e'(0))$ is purely infinite, and
$\tau_{1/3} (1), \tau_{2/3} (1) < 1 - e'(0),$
Lemma 2.9 implies that $\tau_{1/3}$ is
homotopic to $\tau_{2/3}.$ We let
 $ \alpha \to \tau_{\alpha}$ be a homotopy, defined for
$\alpha\in [1/3,3/2],$ with $\tau_{1/3}$ and $\tau_{2/3}$ as
already defined.  We define
\beq
\rho(s_j)(\alpha)=\tau_{\alpha}(s_j) +z_j(\alpha)  \label{rho2}
\eeq
for $\alpha\in [1/3,2/3].$

It is now easy to define the unitary $v$.
Let $q (\alpha) = \rho(1) (\alpha)$, which is equal to
$q_0 (\alpha)$  for $\alpha \not\in [1/3, 2/3]$, and in any case
is equal to $\tau_{\alpha} (1) + e_0' (\alpha)$.  Set
\beq
v=v_0+ q -e'.       \label{v}
\eeq

{\em Step 9:}
We now verify the conclusion of the lemma for these choices of $v$ and
 $\rho.$

We start with (C1).  On $[1/3,2/3],$ we  know that $\tau (s_j)$
commutes with $q - e'$, so (\ref{v0zj}), (\ref{rho2}), and (\ref{v})
imply
\[
\| v \rho (s_j) - \rho (s_j) v \| = \|v_0 z_j - z_j v_0\| < \ep.
\]
On $[0,1/3],$ we first observe that the relations
\[
ye=ey=e\andeqn v(q-e)=(q -e)v=q -e
\]
(see (\ref{v0e'})) imply that $y^*$ commutes with $v.$ Therefore
\beqr
v \rho (s_j) - \rho (s_j) v &
     = & v (\tau(s_j) + z_j) y^* - (\tau(s_j) + z_j) y^* v   \\
 & = &  [v (\tau(s_j) + z_j) - (\tau(s_j) + z_j) v] y^*.
\eeqr
The term in brackets has norm at most $\ep$ for the same reason as
 above,
and $\|y^*\|=1,$ so the estimate is verified here too. The same argument
shows that it holds on $[2/3,1]$ as well.

Next we do (C2). We have $v-\rho(1)=v_0-e',$ so
\[
\|(v-\rho(1) ) -\vph_t((u-1)\otimes 1)\|=
        \|v_0-e'-\vph_t((u-1)\otimes 1)\|<\ep,
\]
by Step 7. On the interval $[1/3, 2/3]$, we have
$(v - \rho(1)) \rho(s_j) = (v_0 - e') \rho_0 (s_j)$, so the inequality
$\|(v - \rho(1)) \rho(s_j) - \vph_t ( (u - 1) \otimes s_j)\| < \ep$
follows immediately from Step 7. On $[0, 1/3]$, we have
\beqr
\lefteqn{
(v - \rho(1)) \rho(s_j) - \vph_t ( (u - 1) \otimes s_j)
}  \\
 & & = (v_0-e')(\tau(s_j)+z_j)y^*-\vph_t((u-1)\otimes s_j)
            \\
 & & = (v_0-e')z_jy^*-\vph_t((u-1)\otimes s_j).
\eeqr
Now $z_jy^*$ differs from $z_j$ only that,
 in the definition, $d$ has been
replaced by $ \tilde{d} = dy^*.$
But $dy^*$ satisfies the properties required for $\tilde{d}$
in Steps 6 and 7, so Step 7 again gives
\[
\|(v-\rho(1))\rho(s_j)-\vph_t((u-1)\otimes s_j)\|<\ep.
\]
The same argument applies to $[2/3, 1]$.

Finally, we do (C3). Using (\ref{tau}), we have
\[
\| \rho_i (s_j) - \psi_i ( 1 \otimes s_j) \|
  = \| (\tau_i (s_j) + z_j (i)) y(i)^* - \psi_i (1 \otimes s_j)\| <\ep.
\]
Also,
\[
\|\psi_i(u\otimes 1)-v(i)\|=\|\vph_t((u-1)\otimes 1)(i)-(v(i)-q(i))\|
\]
\[
\le \|\vph_t((u-1)\otimes 1)-(v-\rho(1))\|<\ep,
\]
as has already been shown in the proof of (C2).  \QED

\vspace{0.6\baselineskip}

The purpose of the following lemma is to show that we can always assume
that hypotheses (H2) and (H3) of Lemma 2.10 hold.

\vspace{0.6\baselineskip}

{\bf 2.11 Lemma} Let $ t \mapsto \varphi_t$ be an asymptotic morphism,
assumed linear, contractive, and *-preserving, from
$\So{m}$ to a \CA\ $A$.
Let $M_m \subset \OA{m}$ be the subalgebra generated by the elements
$s_i s_j^*$.
Then there exists a constant $M$ and an
 asymptotic morphism  $ t \mapsto \psi_t$ from
$\So{m}$ to  $A$ such that:

(1) Each $\psi_t$ is linear, *-preserving, and satisfies
$\|\psi_t \| \leq M$. (We demand neither contractivity nor positivity.)

(2) $\psi_t|_{C_0 (S^1 \setminus \{1\}) \otimes M_m}$
is a \hm\ for each $t$.

(3) Whenever $a \in \So{m}$ is actually in
 $C_0 (S^1 - \{1\}, p_i \OA{m} p_j)$,
then $\psi_t (a^*) \psi_t (a)$ is in the hereditary subalgebra of
$A$
generated by $\psi_t (C_0 (S^1 \setminus \{1\}, {\bf C}p_j))$.

(4)  $\lim_{t \to \infty} (\psi_t(a) - \varphi_t(a)) = 0$ for all
$ a \in \So{m}$.

\noindent
Moreover, if $A$ has the form $C([0,1], B)$ for some \CA\ $B$,
and if ${\rm ev}_i \circ \varphi_t$ is a \hm\ for $i=0,1$ and all $t$,
then  $ t \mapsto \psi_t$ can be chosen to satisfy
 ${\rm ev}_i \circ \psi_t = {\rm ev}_i \circ \varphi_t$
 for $i=0,1$ and all $t$.

\vspace{0.6\baselineskip}

{\em Proof:} We will only do the case involving homotopies with
\hm s at the ends of the homotopy. For simplicity, we will
write
$\varphi_{\alpha, t}$ for
${\rm ev}_{\alpha} \circ \varphi_t$, and similarly for $\psi$ (as we
construct it).
 By a reparametrization of the homotopy, depending on $t$,
we can assume without changing the limiting behavior at $\infty$ that
\[
\varphi_{\alpha, t} = \varphi_{0, t} \,\,\,\,\,\, {\rm for} \,\,\,
        0 \leq \alpha \leq \min(1/t, 1/3)
\]
and
\[
\varphi_{\alpha, t} = \varphi_{1, t} \,\,\,\,\,\, {\rm for} \,\,\,
        \max(1-1/t, 2/3) \leq \alpha \leq 1.
\]
{}From here, we will only construct $\psi_{\alpha, t}$ for $t$ greater
than some $T$, omitting the argument needed to extend the construction
back to smaller values of $t$ while still preserving the condition
$\psi_{i, t} = \varphi_{i, t}$ for $i = 0,1$.

Recall (Definition 3.4 of \cite{Lr4})
that a system $(G,R)$ of \CA\  generators and relations is
{\em exactly stable} if for each $\ep > 0$ there is $\delta > 0$ such
that there is a
 \hm\ $\sigma_{\delta} : C^* (G,R) \to C^* (G, R_{\delta})$
(notation explained below)
with
\[
\| \sigma_{\delta} (g) - g^{(\delta)} \| < \ep \andeqn
    \pi_{\delta} ( \sigma_{\delta} (g) ) = g
\]
for each $g \in G$. Here $C^*(G,R)$ is the (not necessarily unital)
universal \CA\  on the generators $G$ satisfying the relations $R$.
The relations $R_{\delta}$ are the same as the relations $R$, except
``softened'' by $\delta$. If, following Loring's conventions,
$G = \{g_1, \dots, g_{\nu}\}$ and $R$ consists of the relations
$\| g_i \| \leq 1$ for $i = 1, \dots, {\nu}$ and relations of the form
$p_j (g_1, \dots, g_{\nu}) = 0$ for polynomials $p_j$
($j = 1, \dots, \mu$) in ${\nu}$ noncommuting
variables and their noncommuting adjoints (making $2{\nu}$ noncommuting
variables in all), then $R_{\delta}$ consists of the relations
$\| g_i \| \leq 1 + \delta$  and
$\|p_j (g_1, \dots, g_{\nu})\| \leq \delta$.
Unless otherwise specified, we will in fact follow this convention.
To avoid confusion, if $g \in G$ we denote the corresponding
generator of $C^* (G, R_{\delta})$ by $g^{(\delta)}$. We furthermore
let $\pi_{\delta} : C^* (G, R_{\delta}) \to C^*(G,R)$ be the
canonical map sending $g^{(\delta)}$ to $g$.

One easily checks that $C_0 (S^1 \setminus \{1\})$ is generated by
exactly stable relations. It therefore follows from
Theorem 5.7 of \cite{Lr4} that $C_0 (S^1 \setminus \{1\}) \otimes M_m$
is generated by exactly stable relations. Let $(G,R)$ be an exactly
stable system of generators and relations for this algebra, of the
form described above.
For any \CA\ $A$ and $a_1, \dots, a_{\nu} \in A$, we define
$\delta (a_1, \dots, a_{\nu} )$ to be the smallest number $\delta$
for which the $a_i$ satisfy the relations $R_{\delta}$. That is,
$\delta (a_1, \dots, a_{\nu} )$ is the smallest number $\delta$
such that
\[
\|a_i \| \leq 1+\delta \andeqn \|p_j (a_1, \dots, a_{\nu} ) \|
        \leq \delta
\]
for all $i$ and $j$. Now choose a strictly decreasing sequence
$\ep_0 > \ep_1 > \ep_2 > \cdots >0$ such that $\ep_n \to 0$ and:

(1) There is a
\hm\  $\sigma_{\ep_{n+1}} : C^* (G,R) \to C^* (G, R_{\ep_{n+1}})$
with
\[
\| \sigma_{\ep_{n+1}} (g) - g^{(\ep_{n+1})} \| < \ep_n \andeqn
    \pi_{\ep_{n+1}} ( \sigma_{\ep_{n+1}} (g) ) = g
\]
for $g \in G$.

(2) If $a_1, \dots, a_{\nu} \in A$ satisfy the relations $R$, and
$\| a_i - b_i\| < \ep_{n+1}$ for all $i$, then
$\delta (b_1, \dots, b_{\nu} ) \leq \ep_n$.

For simplicity we write $\sigma_n$ for $\sigma_{\ep_n}$, and we define
$\pi_n$, $g^{(n)}$, and $R_n$ following the same convention.
We further let $\pi_{m,n}:  C^* (G, R_{m}) \to C^*(G,R_n)$ be the
obvious projection map. We let
 $\delta_{\alpha,t} =
   \delta(\varphi_{\alpha,t}(g_1), \dots, \varphi_{\alpha,t}(g_{\nu}))$,
and, for  $\delta_{\alpha,t} \leq \ep_n$,
  we let $\kappa_n^{\alpha,t}  :C^* (G, R_{n}) \to B$
be the map sending  $g^{(n)}$ to  $\varphi_{\alpha,t}(g)$.
Note that  $\delta_{\alpha,t} $ is a continuous function of
$\alpha$ and $t$.

Since $t \mapsto \varphi_t$ is an asymptotic morphism, it follows that
$\delta_{\alpha,t} \to 0$ uniformly in $\alpha$ as $t \to \infty$.
Choose $T$ so large that $\delta_{\alpha,t} \leq \ep_4$
whenever $t \geq T$.

Let $\alpha$ and $t\geq T$ satisfy
$\ep_{2n} \geq  \delta_{\alpha,t} \geq \ep_{2n+2}$ for some $n$.
Necessarily $n \geq 2$. Let
$s \in [0,1]$ satisfy
 $\delta_{\alpha,t} = s \ep_{2n} + (1-s) \ep_{2n+2}$.
Observe that
\[
\|\sigma_{2n}(g_i) - (s g_i^{(2n)} + (1-s) \sigma_{2n}(g_i) ) \| \leq
           \ep_{2n-1}.
\]
 Thus the elements $s g_i^{(2n)} + (1-s) \sigma_{2n}(g_i)$ satisfy the
relations $R_{2n-2}$, and so there is a homomorphism
$\tau_{2n,s} : C^* (G, R_{2n-2}) \to C^*(G,R_{2n})$  sending the
generators $g_i^{(2n-2)}$ to these elements. Now define
$\psi_{\alpha,t}^{(0)} : C_0 (S^1 \setminus \{1\}) \otimes M_m \to B$
to be the composite
\[
C_0 (S^1 \setminus \{1\}) \otimes M_m \stackrel{\cong}{\longrightarrow}
   C^*(G,R) \stackrel{\sigma_{2n-2}}{\longrightarrow}
   C^*(G,R_{2n-2}) \stackrel{\tau_{2n, s}}{\longrightarrow}
   C^*(G,R_{2n}) \stackrel{\kappa_{2n}^{\alpha,t} }{\longrightarrow} B.
\]
If instead  $\delta_{\alpha,t} =0$, then define
$\psi_{\alpha,t}^{(0)} = \varphi_{\alpha,t} $.

We claim that $(\alpha,t) \mapsto \psi_{\alpha,t}^{(0)} $
is continuous (in the topology of pointwise convergence), and that
for $a \in C_0 (S^1 \setminus \{1\}) \otimes M_m$ we have
$\psi_{\alpha,t}^{(0)}(a) - \varphi_{\alpha,t}(a) \to 0$
uniformly in $\alpha$ as $t \to \infty$.
This claim will follow if we can show that the two possible
definitions of $\psi_{\alpha,t}^{(0)} $ agree when
$\delta_{\alpha,t} = \ep_{2n+2}$ for some $n$, and that
for $g \in G$ we have
$\| \psi_{\alpha,t}^{(0)}(g) - \varphi_{\alpha,t}(g) \|  < 2 \ep_{2n-3}$
when
$\ep_{2n} \geq  \delta_{\alpha,t} \geq \ep_{2n+2}$.
(The second estimate implies both continuity when
$\delta_{\alpha,t} = 0$ and the desired limiting behavior as
$t \to \infty$. Note that standard arguments show it suffices to
check estimates of this sort on a set of generators.)

We do the first part of this claim first. From the definition
based on the interval
$[\ep_{2n+2}, \ep_{2n}]$, we obtain
\[
\psi_{\alpha,t}^{(0)} = \kappa_{2n}^{\alpha,t}  \circ
    (\sigma_{2n} \circ \pi_{2n-2}) \circ \sigma_{2n-2} =
                  \kappa_{2n}^{\alpha,t}  \circ \sigma_{2n}.
\]
{}From the definition based on the interval
$[\ep_{2n+4}, \ep_{2n+2}]$, we obtain
\[
\psi_{\alpha,t}^{(0)} = \kappa_{2n+2}^{\alpha,t}  \circ
     \pi_{2n, 2n+2} \circ \sigma_{2n} =
                  \kappa_{2n}^{\alpha,t}  \circ \sigma_{2n}.
\]

We now do the second part of the claim.
Let $\ep_{2n} \geq  \delta_{\alpha,t} \geq \ep_{2n+2}$, and
let $s$ be as before.
For $i = 1 , \dots, \nu$ we have
\begin{eqnarray*}
\| \lefteqn{\psi_{\alpha,t}^{(0)}(g_i) - \varphi_{\alpha,t}(g_i) \| } \\
 & &  =
\| \kappa_{2n}^{\alpha,t} \circ \tau_{2n, s} \circ \sigma_{2n-2}(g_i) -
      \kappa_{2n}^{\alpha,t} (g_i^{(2n)}) \|  \\
 & & \leq  \| \tau_{2n, s} \circ \sigma_{2n-2}(g_i) -
                             \tau_{2n, s}(g_i^{(2n-2)}) \| +
           \|  \tau_{2n, s}(g_i^{(2n-2)}) - g_i^{(2n)} \|  \\
 & & \leq  \| \tau_{2n, s} \| \ep_{2n-3} + (1-s) \ep_{2n-1}
            < 2 \ep_{2n-3},
\end{eqnarray*}
as desired.

We now have to extend $\psi_{\alpha,t}^{(0)}$ to the whole of
$\So{m}$ in such a way that the conditions (1), (3), and (4) are
satisfied.
Let $e_{ij} = s_i s_j^*$; this defines a system of matrix units in
$M_m$. Let $\omega_0$ be any state on $e_{11} \OA{m} e_{11}$. Using the
isomorphism $M_m \otimes e_{11} \OA{m} e_{11} \cong \OA{m}$,
we obtain from $\omega_0$ a bounded linear *-preserving map
$\omega: \OA{m} \to M_m$ such that $\omega|_{M_m} = {\rm id}_{M_m}$.
Define
$\gamma_0, \gamma_{ij}: \So{m} \to \So{m}$
by
$\gamma_{0} (a)(\zeta) = \omega(a(\zeta))$ and
$\gamma_{ij} (a)(\zeta) =
   e_{ii}a(\zeta)e_{jj} - \omega(e_{ii}a(\zeta)e_{jj})$.
These maps are the projections for a Banach space internal direct sum
decomposition
\[
\So{m} = C_0 (S^1 \setminus \{1\}) \otimes
                          M_m \oplus \bigoplus_{i,j} Q_{ij},
\]
where $Q_{ij}$ is the range of $\gamma_{ij}$.

For $s \in [1/\pi, \infty)$ define
 $f_s \in  C_0 (S^1 \setminus \{1\}) $ by
\[
f_s (\exp(i \theta)) = \left\{ \begin{array}{ll}
          s\theta          & 0 \leq \theta \leq 1/s  \\
          1                & 1/s \leq \theta \leq 2 \pi - 1/s  \\
          s( 2 \pi-\theta) & 2 \pi - 1/s \leq \theta \leq 2 \pi.
\end{array} \right.
\]
Further let $f_{\infty}$ be the constant function $1$. Note that
$ s \mapsto f_s$, for $s \in [1/\pi, \infty)$, is a continuously indexed
approximate identity for $C_0 (S^1 \setminus \{1\}) $.

We will choose a suitable continuous function
$s : [0,1] \times [T, \infty] \to [1/\pi, \infty]$ such that
$s(\alpha, t) = \infty$ when $\alpha = 0$ or $1$ but not otherwise.
We then define $\psi$ as follows. For
$\frac{1}{2}\min(1/t, 1/3) \leq \alpha \leq \frac{1}{2}\max(1-1/t, 2/3)$
we set
\[
\psi_{\alpha, t} (a) = \left\{ \begin{array}{ll}
   \psi_{\alpha, t}^{(0)} (a)
                    &  a \in C_0 (S^1 \setminus \{1\}) \otimes M_m \\
   \psi_{\alpha, t}^{(0)} (f_{s(\alpha,t)} \otimes e_{ii})
            \varphi_{\alpha, t}(a)
            \psi_{\alpha, t}^{(0)} (f_{s(\alpha,t)} \otimes e_{jj})  &
                                     a \in Q_{ij}.
\end{array} \right.
\]
If instead $0 \leq \alpha \leq \frac{1}{2}\min(1/t, 1/3)$ or
$\frac{1}{2}\max(1-1/t, 2/3) \leq \alpha \leq 1$, then we define
\[
\psi_{\alpha, t} (a) = \left\{ \begin{array}{ll}
   \psi_{\alpha, t}^{(0)} (a)
                    &  a \in C_0 (S^1 \setminus \{1\}) \otimes M_m \\
   \varphi_{\alpha, t}((f_{s(\alpha,t)} \otimes e_{ii})a
                            (f_{s(\alpha,t)} \otimes e_{jj}))
                               &     a \in Q_{ij}.
\end{array} \right.
\]
The definitions agree on the overlaps because, if
$\alpha = \frac{1}{2}\min(1/t, 1/3)$ or
$\alpha = \frac{1}{2}\max(1-1/t, 2/3)$, then
$ \psi_{\alpha, t}^{(0)} =
 \varphi_{\alpha, t}|_{C_0 (S^1 \setminus \{1\}) \otimes M_m}$ and
$ \varphi_{\alpha, t}$ is a homomorphism. There  is also no problem with
the case $s(\alpha,t) = \infty$, because that only occurs when
$(\alpha,t)$ is in the interior of the set for which the second
definition applies.
With this choice, the conditions (1)--(3) of the conclusion are
 satisfied, regardless of the choice of $s$. (The constant
$M$ depends on the norms of the maps involved in the Banach space
direct sum definition used above.) It remains only to
choose $s$ properly so as to ensure that
conclusion (4) is satisfied, and so
that if $a$, $b \in \So{m}$ then
 $\| \psi_t(ab) - \psi_t(a) \psi_t(b) \| \to 0$ as $ t \to \infty$.

Before starting this, we first observe that
the convergence as $t \to \infty$,
\[
\| \psi_{\alpha, t}^{(0)} (a) - \varphi_{\alpha, t} (a) \| \to 0
      \andeqn
\| \varphi_{\alpha, t} (ab) -
     \varphi_{\alpha, t} (a)\varphi_{\alpha, t} (b) \| \to 0
\]
(for $a \in C_0 (S^1 \setminus \{1\}) \otimes M_m$ in the first
expression and $a$, $b \in \So{m}$ in the second),
is uniform in $\alpha \in [0,1]$
 and also in $a$ and $b$ as long as $a$ and $b$
are restricted to compact subsets of the appropriate domains.
This follows from pointwise convergence together with the fact that the
maps involved are linear and bounded uniformly in $t$ and $\alpha$.

Choose a countable subset $\{a_1, a_2, \dots \}$ of
$\So{m}$ whose linear span is dense in $\So{m}$ and such that
each $a_k$ is either in some $Q_{ij}$ or in
$C_0 (S^1 \setminus \{1\}) \otimes M_m$.
 Choose a strictly increasing sequence $s_1 < s_2 < \cdots$
of positive real numbers such that $s_n \to \infty$ and
\[
\| (f_s \otimes 1) a_k (f_s \otimes 1) - a_k \| < 1/n
\]
for $s \geq s_n$ and $k = 1, \dots, n$. Now choose
a strictly increasing sequence $t_1 < t_2 < \cdots$
of positive real numbers such that $t_n \to \infty$, and such that
the following estimates are satisfied.
For $k = 1, \dots, n$, if
$a_k \in C_0 (S^1 \setminus \{1\}) \otimes M_m$,
then we require
\[
\| \psi_{\alpha, t}^{(0)} (a_k) - \varphi_{\alpha, t} (a_k) \| \leq 1/n
\]
for $t \geq t_n$ and all $\alpha$,
while if $a_k \in Q_{ij}$, then we require
\[
\| \psi_{\alpha, t}^{(0)} (f_s \otimes e_{ii}) -
    \varphi_{\alpha, t} (f_s \otimes e_{ii}) \| \leq 1/n, \,\,\,\,\,\,
\| \psi_{\alpha, t}^{(0)} (f_s \otimes e_{jj}) -
    \varphi_{\alpha, t} (f_s \otimes e_{jj}) \| \leq 1/n,
\]
and
\[
\| \varphi_{\alpha, t} ((f_s \otimes e_{ii}) a_k (f_s \otimes e_{jj})) -
     \varphi_{\alpha, t} (f_s \otimes e_{ii})  \varphi_{\alpha, t} (a_k)
                     \varphi_{\alpha, t} (f_s \otimes e_{jj})\| \leq 1/n
\]
for $t \geq t_n$, $s \in [s_{n - 1}, s_{n + 2}]$, and all $\alpha$.
Now let $s$ be a continuous function satisfying:

(1) $s(\alpha, t) \geq s_n$ whenever $t \geq t_n$.

(2) $s(\alpha, t) = \infty$ for $\alpha = 0$ or $1$ and all $t$.

(3) $s(\alpha, t) \leq s_{n + 2}$ whenever $t \leq t_{n + 1}$ and
$\frac{1}{2}\min(1/t, 1/3) \leq \alpha \leq
 \frac{1}{2}\max(1-1/t, 2/3)$.

It is now immediate that if
$a_k \in C_0 (S^1 \setminus \{1\}) \otimes M_m$,
then
$\| \psi_{\alpha, t}^{(0)} (a_k) - \varphi_{\alpha, t} (a_k) \|  \to 0$
as $t \to \infty$, uniformly in $\alpha$. On the other hand, if
$a_k \in Q_{ij}$, then
$(f_s \otimes e_{ii}) a_k (f_s \otimes e_{jj}) =
                (f_s \otimes 1) a_k (f_s \otimes 1)$.
Using this and the fact that $\varphi_{\alpha, t}$ is contractive,
estimates show that for $n \geq k$ we have
\[
\| \psi_{\alpha, t}^{(0)} (a_k) - \varphi_{\alpha, t} (a_k) \|  \leq
     2(\|a_k\| + 1)/n
\]
when the first definition of $\psi_{\alpha, t}$ applies, and
\[
\| \psi_{\alpha, t}^{(0)} (a_k) - \varphi_{\alpha, t} (a_k) \|  \leq 1/n
\]
when the second definition applies.
Thus, in either case,\\
$\| \psi_{\alpha, t}^{(0)} (a_k) - \varphi_{\alpha, t} (a_k) \|  \to 0$
as $t \to \infty$, uniformly in $\alpha$.

It follows from linearity that conclusion (4) is satisfied for $a$ in
a dense subset of $\So{m}$. Since $\sup_t \|\varphi_t \|$ and
$\sup_t \|\psi_t \|$ are finite, a standard argument implies that it in
fact holds for all $a \in \So{m}$.

It remains only to show that
 $\| \psi_t(ab) - \psi_t(a) \psi_t(b) \| \to 0$ as $ t \to \infty$.
The following computation, together with the estimates
$\| \psi_t (a) \| \leq M \|a\|$ and
$\| \varphi_t (a) \| \leq  \|a\|$ for all $t$, shows that the required
condition follows from conclusion (4):
\begin{eqnarray*}
\lefteqn{\| \psi_t(ab) - \psi_t(a) \psi_t(b) \|} \\
   & & \leq \| \psi_t(ab) - \varphi_t(ab) \| +
      \| \varphi_t(ab) - \varphi_t(a) \varphi_t(b) \|   \\
   & & \mbox{} + \| \varphi_t(a) - \psi_t(a)\| \| \varphi_t(b) \| +
      \| \psi_t(a) \| \| \varphi_t(b) - \psi_t(b) \|.
\end{eqnarray*}
\QED

\newpage

\section{ Approximate unitary equivalence of homomorphisms}

\vspace{\baselineskip}

The purpose of this section is to show that if $X \subset S^1$ and
 $\varphi$ and $\psi$ are two injective unital
 \aab\ homomorphisms from $C(X) \otimes M_k (\OA{m})$ (with
$m$ even) to
a purely infinite simple \CA\ $B$, and if $\varphi$ and $\psi$
have
 the same class in $KK$-theory, then $\varphi$ is \ayue\ to $\psi$.
Our first lemma enables us to reduce to the case $k = 1$. This is
followed by a number of technical lemmas which combine to give the
proof in this case. We then define $\ep$-approximately injective
\hm s, and prove several results which enable us to say something
useful about \hm s which are not injective.

In this section,
we will use the definition of $KK$-theory in terms of asymptotic
morphisms \cite{CH}, rather than in terms of Kasparov bimodules
as in \cite{Ks}. Corollary 9 of \cite{CH} shows that both definitions
give the same group when both variables are separable and the first
variable is nuclear. Since $\OA{m}$ is separable and nuclear, these
conditions will hold in all cases of interest to us.

Throughout this section, all \CA s will be assumed separable. Usually,
 $B$ will be a purely infinite
simple \CA\  and $m$ will be a positive even integer.

\vspace{0.6\baselineskip}

{\bf 3.1 Lemma} Suppose that, for a fixed even $m$ and a fixed subset
$X \subset S^1$, the algebra $C = C(X) \otimes \OA{m}$ has the
property that, for any purely infinite simple \CA\ $B$,
 two injective unital \aab\  homomorphisms
$\varphi, \psi: C \to B$, satisfying
 $[\varphi] = [\psi]$ in $KK^0(C,B)$, are necessarily \ayue .
 Then any matrix algebra $M_k(C)$
has the same property.

\vspace{0.6\baselineskip}

{\em Proof:} Let $C$ be as in the statement, and let
$\varphi, \psi: M_k (C) \to B$ be two unital \aab\  homomorphisms
 satisfying $[\varphi] = [\psi]$ in $KK^0(C,B)$.
Let $\{e_{\mu\nu}\}$ be a system of matrix units in $M_k$.
Then in particular
 $[\varphi(e_{11} \otimes 1)] = [\psi(e_{11} \otimes 1)]$ in $K_0 (B)$.
Since $B$ is purely infinite,
 there exists a unitary $z_0 \in B$ such that
$z_0 \varphi(e_{11} \otimes 1) z_0^* =  \psi(e_{11} \otimes 1) $.
Define
\[ z = \sum_{\mu = 1}^k  \psi(e_{\mu 1} \otimes 1) z_0
                                   \varphi(e_{1 \mu} \otimes 1). \]
Then $z$ is a unitary satisfying
$z \varphi(e_{\mu \nu} \otimes 1) z^* =  \psi(e_{\mu \nu} \otimes 1) $.
Replacing $\varphi$ by $z\varphi(-)z^*$, we can assume that
$\varphi|_{M_k \otimes 1} = \psi|_{M_k \otimes 1} $.
Set
 $ B_0 = \varphi(e_{11} \otimes 1) B \varphi(e_{11} \otimes 1) $.
Identifying $B$ with $M_k (B_0)$
in the obvious way, we can now assume that
$\varphi = {\rm id}_{M_k} \otimes \varphi_0$ and
$\psi = {\rm id}_{M_k} \otimes \psi_0$ for suitable unital homomorphisms
$\varphi, \psi : C \to B_0$.

We now claim that
 $[\varphi_0] = [\psi_0]$ in $KK^0(C,B_0)$.
Let $\sigma \in KK^0 ({\bf C}, M_k)$ be the class of the map which
sends $1$ to a rank one projection. Then $\sigma$ has an inverse
$\tau \in KK^0 (M_k, {\bf C})$. It follows that
\[ [\varphi_0] = (\sigma \times {\rm id}_C) \times
 [ {\rm id}_{M_k} \otimes \varphi_0] \times (\tau \times {\rm id}_{B_0})
= (\sigma \times {\rm id}_C) \times [\varphi] \times
                                      (\tau \times {\rm id}_{B_0}) \]
as elements of $KK^0 (C, B_0)$. The same calculation applies to
$\psi_0$ and $\psi$. Therefore  $[\varphi_0] = [\psi_0]$.

Clearly $\varphi_0$ and $ \psi_0$ are both injective,
 and they
are \aab\ by Corollary 1.10. The hypothesis on $C$ implies that
$\varphi_0$ is \ayue\ to $\psi_0$. It follows immediately that
$\varphi$ is \ayue\ to $\psi$. \QED

\vspace{0.6\baselineskip}

{\bf 3.2 Lemma}
For any $\ep>0,$ there exists $\delta>0$ such that: If $B$ is a purely
infinite simple \CA , and if
$v \in B$ is unitary and $s_1, \dots, s_m\in B$
are isometries with
\[
 \sum_{j=1}^m s_j s_j^*=1,  \,\,\,\,\,\,
\|v s_j - s_j v\|<\delta,  \,\,\,\,\,\,{\rm and}\,\,\,\,\,\,
v\in U_0(B),
\]
then there are $\lambda_1, \dots,\lambda_l\in S^1$ and
mutually orthogonal projections $p_1, \dots,p_l$ in $B$ such that
\[
\|v-\sum_{i=1}^l\lambda_i p_i\|<\ep   \,\,\,\,\,\,
{\rm and}   \,\,\,\,\,\,
\left[\sum_{j=1}^m s_j p_i s_j^* \right]=[p_i]
 \,\,\, {\rm in}\,\,\, K_0(B).
\]

\vspace{0.6\baselineskip}

{\it Proof:}  Since $v\in U_0(B)$, we use \cite{Ph1} to choose
 $\lambda_1, \dots ,\lambda_l\in S^1$ and \mops\  $p_1, \dots, p_l$
 such that
\[ \|v-\sum_{i=1}^l\lambda_ip_i\|<\ep/4 .\]
What we need to do here is to choose the $p_i$ properly so as to satisfy
the other requirement.

We first consider the following situation: $v$ is actually in $qBq$ for
some projection $q$ satisfying $m[q] = [q]$ in $K_0 (B)$, and $v$ is
approximated as above in such a way that $p_i \in qBq$ and
\[
|\lambda_i-\lambda_{i+1}|<\ep/2 \,\,\, {\rm for} \,\,\, i=1, \dots ,l.
\]
However, no assumption is made on $\| v s_j - s_j v \|$.

Note that the existence of the isometries $s_j$ implies that
 $m[1] = [1]$. Furthermore, if  ${\rm sp}(v)=S^1$ and
the $\lambda_i$ are ordered cyclically, then the estimate above must
hold. Thus, this situation covers the case
$v \in B$ and ${\rm sp}(v)=S^1$.

Since $qBq$ is purely infinite and simple, there is a nonzero
projection $p_1'\le p_1$ such that $[p']=[q]$ in $K_0(qBq)$.
There is then a nonzero
projection $p_2''\le p_2$ such that $[(p-p_1')+p_2'']=[q]$; set
$p_2'=(p-p_1')+p_2''.$ Proceeding inductively, we obtain
\mops\  $p_i'$ for $i = 1, \dots, l-1$ such that $[p_i'] = [q]$,
 $p_i'\le p_{i-1}+p_i$  for $i = 1, \dots, l-1$,
and $\sum_{i=1}^{l-2}p_i\le \sum_{i=1}^{l-1}p_i'.$ (We take $p_0 = 0$.)
Finally, let $p_l'= q - \sum_{i=1}^{l-1}p_i'$.

We have $m[q] = [q]$ and
 $m[p_i']=[p_i']$ for $i = 1, \dots, l-1$. Therefore
\[
m[p_l'] = m([q]-\sum_{i=1}^{l-1}[p_i'])=[q]-\sum_{i=1}^{l-1}[p_i']
  =[p_l']
\]
too.  Thus
\[
\left[\sum_{j=1}^m s_kp_i's_k^*\right]=m[p_i']=[p_i']
\]
for all $i$.
Since $p_i'\le p_{i-1}+p_{i}$ and $ |\lambda_i-\lambda_{i+1}|<\ep/2 $,
we conclude that
\[
\|v-\sum_{i=1}^l\lambda_ip_i'\| \leq \|v-\sum_{i=1}^l\lambda_ip_i\|
  + \| \sum_{i =1}^l\lambda_ip_i - \sum_{i=1}^l\lambda_ip_i'\| < \ep.
\]
This completes the proof in this situation.

Now we prove the lemma. Let $X= {\rm sp}(v).$
 We may assume (using the remark above) that
$X\not=S^1.$ Write
$X=\coprod_{n=1}^N X_n ,$ the disjoint union of closed subsets $X_n$,
in such a way that:

(1) ${\rm dist}(X_k, X_n) \geq \ep/6$ for $k \neq n$, and

(2) For fixed  $n$, there are  $\lambda_1', \dots ,\lambda_l'\in X_n$
such that $|\lambda_i' -\lambda_{i+1}' |<\ep/6$  for $i=1, \dots, l$,
and such that every $\lambda \in X_n$ is within $\ep/6$ of some
$\lambda_i$.

Let $q_n$ be the spectral projection for $v$ corresponding to $X_n$.
Then $ v=\sum_{n=1}^N v_n $ with $v_n = q_n v q_n$ (so that
${\rm sp}(v_n)=X_n$ in $q_n B q_n$).
For $\delta$ small enough (depending only on $\ep$ and $m$),
if $\| v s_j - s_j v \| < \delta$, then
\[
\|q_n s_j - s_j q_n \|<1/m \,\,\, {\rm for} \,\,\, j = 1, \dots, m
     \,\,\, {\rm and}\,\,\, n = 1, \dots, N.
\]
It follows that
\[
\left\|q_n - \sum_{j=1}^m s_j q_n s_j^* \right\|<1.
\]
Therefore
\[
\left[\sum_{j=1}^m s_j q_n s_j^* \right]=[q_n],
\]
whence $m[q_n] = [q_n]$. We now approximate each $v_n$ to within $\ep/6$
 by a unitary $\sum_i \lambda_i p_i $ with finite spectrum
in $q_n B q_n$. We assume the numbers $\lambda_i$ are ordered
 cyclically. Since there can be no gaps in
${\rm sp}(v_n)$ of length greater than $\ep/6$, it follows that
$| \lambda_i - \lambda_{i + 1} | < \ep/2$. It now
 follows from the special
situation considered above that there are \mops\ $p_i' \in q_n B q_n$
such that
\[ \left\| v_n - \sum_i \lambda_i p_i' \right\| < \ep \,\,\,\,\,\,
   {\rm and} \,\,\,\,\,\,
   \left[ \sum_j s_j p_i' s_j^* \right] = [p_i']. \]
We get the conclusion of the lemma by summing over $n$. \QED

\vspace{0.6\baselineskip}

{\bf 3.3 Lemma} Let $m$ be even.
For any $\ep>0$ there is $\delta>0$  such that if
a unital homomorphism $\varphi$ from $\OA{m}$ to a purely infinite
simple \CA\ $B$ and a unitary
 $v \in U_0 (B)$ satisfy:
\[ \| v \varphi(s_j) - \varphi(s_j) v \| < \delta \,\,\,\,\,\,
     {\rm for} \,\,\, j = 1, \dots, m, \]
then there is a \hm\ $\psi : \SO{m} \to M_2 (B)$ such that:

(1) $\psi(u \otimes 1)$ has finite spectrum.

(2) $\| v \oplus v^* - \psi(u \otimes 1)\| < \ep$.

(3) $\| \varphi(s_j) \oplus \varphi(s_j) - \psi(1 \otimes s_j )\| < \ep$
for $j = 1, \dots, m$.

\vspace{0.6\baselineskip}

{\em Proof:} Let $\delta_1$, $\delta_2 > 0$. (We will choose $\delta_1$
and $\delta_2$ later.)

By the previous lemma, there are \mops\
$p_1, p_2,$
$\dots, p_l$ and a unitary
 $v_0 = \sum_{i = 1}^l \zeta_i p_i \in B$ with finite spectrum
such that $\| v_0 - v\| < \delta_1/(2m)$ and, for each $i$, the
projection
\[ q_i = \sum_{j = 1}^m \varphi(s_j) p_i \varphi(s_j)^* \]
satisfies $[q_i] = [p_i]$. Since $B$ is purely infinite, $p_i$ is
 unitarily equivalent to $q_i$. Furthermore,
\[ \sum_{i = 1}^l q_i = \sum_{i = 1}^l p_i =  1, \]
so in fact there is a unitary $U \in B$ such that $U^* q_i U = p_i$
for all $i$. Choose $z \in U(p_1 B p_1)$ such that $[z] = -[U]$
in $K_1 (p_1 B p_1) \cong K_1 (B)$. Replacing $U$ by
$U[z + (1-p_1)]$, we may assume $[U] = 0$ in $K_1 (B)$, and thus that
$U \in U_0 (B)$.

Set $t_j^{(i)}=U^*\varphi(s_j)p_i.$ Then we have
\[
(t_j^{(i)})^*t_j^{(i)}=p_i\,\,\,\,\,\,{\rm and}
     \,\,\,\,\,\,\sum_{j=1}^m t_j^{(i)}(t_j^{(i)})^*=p_i.
\]
For $i\neq i'$ and any $j, j'$, we further have
\[
t_j^{(i)}t_{j'}^{(i')}=t_{j'}^{(i')}t_j^{(i)}=0.
\]
Set $t_j = \sum_{i=1}^l t_j^{(i)}$. Then the $t_j$ are isometries
 which commute with $v_0$ and
which generate a $C^*$-subalgebra of $B$ isomorphic to $\OA{m}$.
Thus we can define $\psi_1 : \SO{m} \to B$
by $\psi_1 (u) = v_0$ and  $\psi_1 (s_j) = t_j$,
and also define  $\psi_2 : \SO{m} \to B$
by $\psi_2 (u) = v_0^*$ and  $\psi_2 (s_j) = t_j$.

Define $\lambda : B \to B$ by
  \[ \lambda(a)=\sum_{j=1}^m \varphi(s_j) a \varphi(s_j)^*. \]
If $\delta$ is chosen less than $\delta_1/(4m)$, then
$\|v_0\varphi(s_j)-\varphi(s_j)v_0\| < \delta_1/m$ for $ j=1,\dots, m$.
It follows that
$\|\lambda(v_0)-v_0\|<\delta_1$.
We compute
\[
U^*\lambda(v_0)U =
 U^* \left(\sum_{j=1}^m\varphi(s_j)v_0\varphi(s_j^*) \right)U
=U^* \left(\sum_{j=1}^m\zeta_j q_j \right)U=\sum_{j=1}^m\zeta_j p_j=v_0.
\]
So
\[
\|U^*v_0U-v_0\| = \|U^*v_0U-U^*\lambda(v_0)U\| < \delta_1. \]
Theorem 2.6 of \cite{Ln6} implies that, if $\delta_1$ is small enough,
then there are two commuting unitaries $U_0, V_0 \in B$ such that
\[\|U- U_0\|<\delta_2/2 \,\,\,\,\,\,{\rm and}
            \,\,\,\,\,\, \|V - v_0\|<\delta_2/2 . \]

Define $\tilde{U} = U_0 \oplus U_0, \tilde{V} = V_0 \oplus V_0^*
 \in M_2 (B)$.
Let $\mu : C(S^1 \times S^1) \to M_2(B)$ be the homomorphism sending the
two canonical generators of $C(S^1 \times S^1)$ to
  $\tilde{U}$ and $\tilde{V}$.
 Notice that ${\tilde U}$ commutes with the path
\[
\alpha \mapsto \begin{array}{c}\left(\begin{array}{cc}
V_0 & 0\\
0 & 1\end{array}\right) \left (\begin{array}{cc}
\cos(\alpha) & \sin(\alpha) \\
-\sin(\alpha) & \cos(\alpha) \end{array}\right)
\left(\begin{array}{cc}
V_0^* & 0\\
0 & 1 \end{array}\right)
\left (\begin{array}{cc}
\cos(\alpha) & -\sin(\alpha)\\
\sin(\alpha) & \cos(\alpha) \end{array}\right)
\end{array}
\]
from $1$ to $\tilde{V}$.
Therefore $\mu$ is homotopic to a \hm\ with a nontrivial kernel.
Letting $d: K_0 (C(S^1 \times S^1)) \to C(X, {\bf Z})$ be the
dimension map, it follows that $\mu_*|_{\ker(d)} = 0$. Since
both  $\tilde{U}$ and $\tilde{V}$ are in $U_0 (M_2 (B))$, it also
follows that $\mu_*$ is zero on $K_1$. Theorem 4.14 of \cite{Ln4}
now implies that $\tilde{U}$ and $\tilde{V}$  can be arbitrarily
closely approximated by commuting unitaries
 $\tilde{U'}$ and $\tilde{V'}$
with finite spectrum.  We require
\[ \|\tilde{U'}- \tilde{U} \|<\delta_2/2 \,\,\,\,\,\,{\rm and}
            \,\,\,\,\,\, \|\tilde{V'} - \tilde{V} \|<\delta_2/2 . \]

We now follow the proof of Lemma 3.1 of \cite{Ln5}. We take
the sets of isometries to be
$\{U^*\varphi(s_j) \oplus U^*\varphi(s_j)\}$ (for $\{s_j\}$) and
 $\{\varphi(s_j) \oplus \varphi(s_j)\}$ (for $\{t_j\}$),
the unitary $u$ there to be $U \oplus U$, and the unitary $w$ there
to be $v_0 \oplus v_0^*$.
 We have unitaries
 $\tilde{U'}$ and $\tilde{V'}$ with finite spectrum such that
\[ \| U \oplus U - \tilde{U'} \|<\delta_2 \,\,\,\,\,\,{\rm and}
       \,\,\,\,\,\, \|\ v_0 \oplus v_0^* - \tilde{V'} \|<\delta_2 . \]
Taking $\delta_2$ as the $\sigma$ at the beginning of the proof of
Lemma 3.1 of \cite{Ln5}, if $\delta_2$ is small enough, there exists
a
 unitary $W\in A$ such that
\[
\|W^*(\varphi(s_j)\oplus \varphi(s_j))W -
      U^*\varphi(s_j)\oplus U^*\varphi(s_j)\|<\ep/2
\]
for $j=1, \dots, m$ and
\[
\|W(v_0\oplus v_0^*)-(v_0\oplus v_0^*)W\|<\ep/2.
\]

Define $\psi(a) = W (\psi_1 (a) \oplus \psi_2 (a)) W^*$. We claim
this is the required homomorphism. Since $v_0$ has finite spectrum,
so does $\psi(u) = W (v_0 \oplus v_0^*) W^*$.
The estimates above show that
\[ \| W (v_0 \oplus v_0^*) W^* - v \oplus v^* \| < \ep/2
        + \delta_1/(2m) \]
(since $\| v - v_0 \| < \delta_1/(2m)$) and
\[ \| \psi(s_j) -  \varphi (s_j) \oplus \varphi (s_j) \| =
 \| W(U^*\varphi(s_j)\oplus U^*\varphi(s_j))W^* -
               \varphi (s_j) \oplus \varphi (s_j) \| < \ep/2.  \]
We may assume $\delta_1$ has been chosen to satisfy
$\delta_1/(2m) < \ep/2$, and so we are done.
\QED

\vspace{0.6\baselineskip}

The following lemma essentially says that asymptotic \hm s from
$\SO{m}$ can be required to be \hm s when restricted to
$C(S^1) \otimes 1$ and $1 \otimes \OA{m}$.

\vspace{0.6\baselineskip}

{\bf 3.4 Lemma}
Let $t \mapsto \varphi_t$ be a unital asymptotic morphism from $\SO{m}$
to a unital $C^*$-algebra $B$. (That is, we assume $\varphi_t (1) \to 1$
as $t \to \infty$.) Then there exists a continuous unitary path
$t \mapsto v_t$ in $B$ and a
continuous path $t \mapsto \psi_t$ of unital \hm s
from $\OA{m}$ to $B$ such that the expressions
\[ \| v \psi_t (s_j ) - \psi_t (s_j) v \| , \,\,\,\,\,\,
   \|  v - \varphi_t (u \otimes 1)  \|  \,\,\,\,\,\, {\rm and}
     \,\,\,\,\,\,   \|  \psi_t (s_j ) - \varphi_t (1 \otimes s_j)  \|
\]
all converge to $0$ as $t \to \infty$.

\vspace{0.6\baselineskip}

{\em Proof:}
 It suffices to show that $v_t$ and $\psi_t$ can be
constructed for $t \geq t_0$ for some $t_0$. (We then just use the
values at $t_0$ for $t < t_0$.)

Note that $\varphi_t (u \otimes 1) \varphi_t (u^* \otimes 1) - 1 \to 0$,
since $t \to \varphi_t$ is a unital asymptotic \hm . We choose $t_0$
so large that $\varphi_t (u \otimes 1)$ is invertible for $t \geq t_0$,
and set
\[ v_t = \varphi_t (u \otimes 1)
       [ \varphi_t (u \otimes 1)^* \varphi_t (u \otimes 1)]^{-1/2}. \]
Certainly $\|  v - \varphi_t (u \otimes 1)  \| \to 0$ as $t \to \infty$.

Similarly, if $t \geq t_0$ and $t_0$ is sufficiently large, we can use
the exact stability of the defining relations of $\OA{m}$, as in the
last paragraph of the proof of Lemma 2.1, to construct $\psi_t (s_j)$
from the elements $\varphi_t ( 1 \otimes s_j )$.
Exact stability implies in particular that
$\| \psi_t (s_j) - \varphi_t (1 \otimes s_j)\| \to 0$ as $t \to \infty$.
The relation $\| v \psi_t (s_j ) - \psi_t (s_j) v \| \to 0$ as
$t \to \infty$ now follows from the other two relations and the fact
that $t \mapsto \varphi_t$ is an asymptotic \hm . \QED

\vspace{0.6\baselineskip}

{\bf 3.5 Lemma} Let $m$ be even, and let $B$ be a purely infinite simple
\CA .
Let $\varphi_0$ and $\varphi_1$ be two injective
 \aab\  nonunital homomorphisms
 from $\SO{m}$ to
$B$ such that $[\varphi_0] = [\varphi_1]$ in $KK^0 (\SO{m}, B)$.
Then for any $\ep> 0$ and $\delta > 0$ there are $L$,
partial isometries $v_0, \dots, v_L \in B$, and \hm s
$\psi_0, \dots, \psi_L : \OA{m} \to B$ such that:

(1) $v_l v_l^* = v_l^* v_l = \psi_l (1)$ for all $l$.

(2) $\| v_l \psi_l (s_j) -  \psi_l (s_j) v_l \| < \delta$
 for all $l$ and $j$.

(3) $v_0 = \varphi_0 (u \otimes 1)$, $v_L = \varphi_1 (u \otimes 1)$,
$\psi_0 = \varphi_0 |_{1 \otimes \OA{m}}$, and
$\psi_L = \varphi_1 |_{1 \otimes \OA{m}}$.

(4) $\| v_l - v_{l-1} \| < \ep$ and
 $\| \psi_l (s_j) - \psi_{l-1} (s_j) \| < \ep$ for all $l$ and $j$.

\vspace{0.6\baselineskip}

Essentially, this lemma promises the existence of a discrete version
of a homotopy from $\varphi_0$ to $\varphi_1$ via
approximate homomorphisms from $\SO{m}$ to $B$ which are \hm s when
restricted to $C(S^1) \otimes 1$ and $1 \otimes \OA{m}$.

\vspace{0.6\baselineskip}

{\em Proof:}
Let us say that \hm s
 $\varphi_0$ and $\varphi_1$ are connected by an $(\ep, \delta)$-chain
 if there exist $L$, $v_0, \dots, v_L$, and
$\psi_0, \dots, \psi_L$, such that the conclusions (1)--(3) hold.
Note that this relation is an equivalence relation, and that homotopic
\hm s are necessarily related in this manner.

 Since $[\varphi_0(1)]=[\varphi_1(1)]$ and both are not
the identity of $B,$ there is a unitary path $t\to U_t\in B$
such that $U_0=1$ and $U_1^*\varphi_1(1)U_1=\varphi_0(1).$
We may replace $\varphi_1$ by the homotopic \hm\ $U_1^*\varphi_1(-)U_1$,
and thus assume that $\varphi_0(1)=\varphi_1(1).$
Let $q\in (1-\varphi_0(1))B(1-\varphi_0(1))$ be a nonzero projection
such that $q\not=1-\varphi_1(1)$. Let $h\in qBq$ be selfadjoint and
satisfy ${\rm sp}(h)=[0,1].$ (Since $B$ is a non-elementary simple
$C^*$-algebra, such an element exists by page 61 of \cite{AS}.)
Then the hereditary $C^*$-subalgebra $B_1$ generated by $h+\varphi_0(1)$
is nonunital but $\sigma$-unital. It follows from
\cite{Bn1} and Theorem 1.2 (i) of \cite{Zh1}
that $B_1\cong B\otimes {\cal K}.$

The maps $\varphi_0 |_{\So{m}}$ and $\varphi_1 |_{\So{m}}$,
 regarded as \hm s
from $\So{m}$ to $B_1$, have the same class in
$KK$-theory.
Since $B_1$ is stable,
it follows from Corollary 5.2 of \cite{DL} that there is
a homotopy $t \to \mu^{(0)}_t$ of asymptotic morphisms from $\So{m}$ to
$B \otimes K$ with $\mu^{(0)}_t = \varphi_t |_{\So{m}}$ for $t = 0,1$.
Applying the version of Lemma 2.11 for homotopies, we may assume that
$\mu^{(0)}_t$ in addition has the properties (1)--(3) of the conclusion
of that lemma. (Note that we do not get to assume that $\mu^{(0)}_t$
is contractive or positive.)
Take a nonzero projection
 $p \in  B_2 = (1-q - \varphi_1(1))B(1-q-\varphi_1(1))$
such that $p\not=1-q-\varphi_1(1).$
Now let
 $\nu: \SO{m} \to pB_2 p$
be the (nonunital) homomorphism constructed in Lemma 1.15. Set
 $\mu_t = \mu^{(0)}_t + \nu |_{\So{m}}$.
Then the asymptotic morphism $t \to \mu_t$ satisfies the hypotheses
of Lemma 2.10, with the $\psi_i$ there taken to be
$\varphi_i + \nu$. (Note that the full spectrum condition, hypothesis
(H6) of Lemma 2.10, is satisfied because $\nu$ is injective.)

By applying Lemma 2.10, we see that  $\varphi_0+\nu$ and $\varphi_1+\nu$
are connected by an $(\ep, \delta)$-chain.
Therefore the proof will be complete if we show that, for $i=0,1$, the
\hm s $\varphi_i+\nu$ and $\varphi_i$
are connected by an $(\ep, \delta)$-chain.

Since $\varphi_i$ is \aab, we have
$$
\varphi_i\aueeps{\ep} \varphi_i+\nu.
$$
So
 $\varphi_i$ and $W^*(\varphi_i+\nu)(-)W$, for some unitary
$W\in (p+q+\varphi_0(1))B(p+q+\varphi_0(1)),$
are connected by an $(\ep, \delta)$-chain (with $L = 1$).
Since $B$ is purely infinite, there is a unitary
$V\in (1-(p+q+\varphi_i(1))\tilde{B}(1-(p+q+\varphi_i(1))$ such that
$[V]=-[W]$ in $K_1(B).$ Set $U=W+V.$ Then $U\in U_0(B).$ Therefore
$W^*(\varphi_i+\nu)(-)W = U^*(\varphi_i+\nu)(-)U$
is actually homotopic  to $\varphi_i.$  \QED

\vspace{0.6\baselineskip}

{\bf 3.6 Lemma} Let $m$ be even, and let $B$ be a purely infinite simple
\CA . Let $X \subset S^1$ be a closed proper subset, and
let $\varphi_0$ and $\varphi_1$ be two injective
nonunital homomorphisms
 from $C(X) \otimes \OA{m}$ to
$B$ such that $[\varphi_0] = [\varphi_1]$ in
 $KK^0 (C(X) \otimes \OA{m}, B)$.
Then for any $\ep> 0,$  there are $L$,
partial isometries $v_0, \dots, v_L \in B$, and \hm s
$\psi_0, \dots, \psi_L : \OA{m} \to B$  such that:

(1) $v_l v_l^* =v_l^* v_l = \psi_l (1)$ for all $l$.

(2) $ v_l \psi_l (s_j)=  \psi_l (s_j) v_l $
 for all $l$ and $j$.

(3) $v_0 = \varphi_0 (u \otimes 1)$, $v_L = \varphi_1 (u \otimes 1)$,
$\psi_0 = \varphi_0 |_{1 \otimes \OA{m}}$, and
$\psi_L = \varphi_1 |_{1 \otimes \OA{m}}$.

(4) $\| v_l - v_{l-1} \| < \ep$ and
 $\| \psi_l (s_j) - \psi_{l-1} (s_j) \| < \ep$ for all $l$ and $j$.

(5) ${\rm sp}(v_l) \subset X$ for all $l$.

\vspace{0.6\baselineskip}

{\it Proof:}
Fix $\ep > 0$.

We first consider the case $X$ finite, say
 $X = \{\lambda_1, \dots, \lambda_N\}$.
In this case, we will actually construct a homotopy
 $t \mapsto \varphi_t$ from $\varphi_0$ to $\varphi_1$.
The required $v_l$ and $\psi_l$ will then be given by
\[
v_l = \varphi_{t_l} (u \otimes 1) \andeqn
    \psi_l(a) = \varphi_{t_l} (1 \otimes a)
\]
for a sufficiently fine partition
$0 = t_0 < t_1 < \cdots < t_L = 1$ of $[0,1]$.

Set $e_n = \chi_{\{\lambda_n\}}$. For $i=0,1$ and $n = 1, \dots, N$
define $\psi_i^{(n)}: \OA{m} \to B$ by
$\psi_i^{(n)} (a) = \varphi_i (e_n \otimes a)$.
Since $[\varphi_0]=[\varphi_1]$ in $KK^0 (C(X)\otimes \OA{m},B),$
we have
$[\varphi_0(e_n)]=[\varphi_1(e_n)].$ Therefore there
 is a path of unitaries $t \mapsto w_t$
in $B$ such that
\[
w_0^*\varphi_0(e_n)w_0=\varphi_0(e_n) \andeqn w_1^*\varphi_0(e_n)w_1
=\varphi_1(e_n).
\]
So, without loss of generality, we may further assume that
$\varphi_0(e_n)=\varphi_1(e_n).$
Now the condition that $[\varphi_0]=[\varphi_1]$
in $KK^0 (C(X)\otimes \OA{m} ,B)$
implies that
\[
[\psi_0^{(n)}]=[\psi_1^{(n)}]
          \,\,\,\,\,\,\, {\rm in} \,\,\, KK^0 (\OA{m}, B).
\]

Lemma 2.9 implies that $\psi_0^{(1)}$ is homotopic to $\psi_1^{(1)}$.
Call the homotopy $t \mapsto \psi_t^{(1)}$.
Choose a unitary path $t \mapsto v_t^{(1)}$ such that $v_0^{(1)} = 1$
 and
$(v_t^{(1)})^* \psi_0^{(1)}(1) v_t^{(1)} =\psi_t^{(1)}(1)$.
Then $v_1^{(1)}$ commutes with $\varphi_1 (e_1)$.
Now use Lemma 2.9 to produce a homotopy
 $t \mapsto  \psi_t^{(2) \prime}$
of \hm s from $\OA{m}$ to
$(1 - \psi_0^{(1)}(1)) B (1 - \psi_0^{(1)}(1))$, such that
$\psi_0^{(2) \prime} = \psi_0^{(2)}$ and
$\psi_1^{(2) \prime} = v_1^{(1)} \psi_1^{(2)}(-) (v_1^{(1)})^*$.
Define $\psi_t^{(2)} = (v_t^{(1)})^* \psi_t^{(2) \prime}v_t^{(1)}$.
 This yields
a homotopy from $\psi_0^{(2)}$ to $\psi_1^{(2)}$ whose range is
 orthogonal to that of $\psi_t^{(1)}$ at each $t$. Now let
$t \mapsto v_t^{(1)}$ be a unitary path which conjugates
$\psi_0^{(1)}(1) + \psi_0^{(2)}(1)$ to
 $\psi_t^{(1)}(1) + \psi_t^{(2)}(1)$, and construct a homotopy
 $ t \mapsto \psi_t^{(3)}$, etc. Then define the required homotopy from
$\varphi_0$ to $\varphi_1$ by
$\varphi_t (f \otimes a) = \sum_n f(\lambda_n) \psi_t^{(n)} (a)$.

Now consider the case in which $X$ contains no arc of length greater
than $\ep/2$. We can then write $X = \coprod_{n=1}^N X_n$, where each
$X_n$ is closed and contained in an arc of length at most $\ep/2$.
Let $e_n = \chi_{X_n}$, and let $u_n = e_n u e_n$. Then the
$e_n$ are \mops\ which sum to $1$, and the $u_n$ are unitaries in
$e_n C(X) e_n$ which sum to $u$.
For each $n$, choose some $\lambda_n\in X_n.$
Then
\[
\|u-\sum_{n=1}^N \lambda_n e_n\|<\ep.
\]
Let $\varphi_0',\, \varphi_1': C(X)\otimes \OA{m} \to B$
be the homomorphisms
defined by
\[
\varphi_i'(f \otimes a) =
          \sum_{n=1}^N f(\lambda_k)\varphi_i(e_n \otimes a).
\]
We then have
\[
\| \varphi_i'(u \otimes 1) - \varphi_i (u \otimes 1) \|
    = \| \varphi_i (u-\sum_{n=1}^N \lambda_n e_n) \|<\ep
\]
and
\[
\| \varphi_i'(1 \otimes s_j) - \varphi_i (1 \otimes s_j) \| = 0.
\]
It therefore suffices to prove the result for
 $\varphi_0'$ and $\varphi_1'$ in place of  $\varphi_0$ and $\varphi_1$.
Since  $\varphi_0'$ and $\varphi_1'$ define injective homomorphisms from
$C(X')\otimes \OA{m} \to B$,
with $X' = \{ \lambda_1, \dots, \lambda_N\}$,
this follows from the case already done.

Now we consider the general case. The connected components of $X$
are all either points or closed arcs.  Let $I_1,I_2,\dots,I_N$ be the
connected components of $X$ which are arcs of length at least
$\ep/2$. Note that no $I_n$ is equal to $S^1$. For each $n$, let
$J_n$ be a closed arc which contains $I_n$, which is at most $\ep/2$
longer than $I_n$, and whose endpoints are not in $X$. We further
require that the $J_n$ be disjoint. Let $h : X \to X$ be the continuous
function which is the identity on
\[
X \setminus \bigcup_{n = 1}^N (J_n \setminus I_n)
\]
and which sends each point of $X \cap (J_n \setminus I_n)$ to the
nearest
endpoint of $I_n$. Note that $\| h(\lambda) - \lambda \| < \ep$ for
all $\lambda \in X$.
Let $\varphi_0',\, \varphi_1': C(X)\otimes \OA{m} \to B$
be the homomorphisms defined by
\[
\varphi_i'(f \otimes a) = \varphi_i((f \circ h) \otimes a).
\]
We then have, as in the previous case,
\[
\| \varphi_i'(u \otimes 1) - \varphi_i (u \otimes 1) \| <\ep \andeqn
   \varphi_i'(1 \otimes s_j) = \varphi_i (1 \otimes s_j).
\]
Therefore, as before, we need only prove the result for
 $\varphi_0'$ and $\varphi_1'$. We may regard them as injective maps
$C(X')\otimes \OA{m} \to B$, with
\[
X' = \left(X \setminus \bigcup_{n = 1}^N J_n \right) \cup
    \bigcup_{n = 1}^N I_n \subset X.
\]

Now let $t \mapsto h_t$ be a homotopy of continuous maps from $X'$
to $X'$ such that $h_0 = \id_{X'}$ and $h_1$ is constant with
value $\lambda_n \in I_n$ on each
arc $I_n$. Then $\varphi_i'$ is homotopic to the \hm\  $\varphi_i''$
given by
\[
\varphi_i''(f \otimes a) = \varphi_i'((f \circ h_1) \otimes a).
\]
As we saw in the proof of the case of finite $X$, we can replace \hm s
by homotopic ones. Therefore it suffices to obtain the result for
 $\varphi_0''$ and $\varphi_1''$.
These \hm s  may be  regarded as injective maps
$C(X'')\otimes \OA{m} \to B$, with
\[
X'' = \left(X' \setminus \bigcup_{n = 1}^N I_n \right) \cup
   \{ \lambda_1, \dots, \lambda_N\}
 \subset X'.
\]
The set $X''$ contains no arcs of length greater than $\ep/2$, so we are
reduced to the previous case. \QED

\vspace{0.6\baselineskip}

{\bf 3.7 Proposition} Let $m$ be even, and let $B$ be a purely
 infinite simple \CA .
Let $X$ be a closed subset of $S^1$, and let
$\varphi_0$ and $\varphi_1$ be two injective \aab\  \hm s from
$C(X) \otimes M_n(\OA{m})$ to $B$, with the same class in
 $KK^0 (C(X) \otimes M_n (\OA{m}), B)$
and either both unital or
both nonunital. Then $\varphi_0$ and $\varphi_1$ are \ayue .

\vspace{0.6\baselineskip}

{\it Proof:} By Lemma 3.1, we may assume that $n =1.$

We first do the nonunital case.
Let $\ep > 0$. Choose $L$, $v_0, \dots, v_L \in U(B)$, and
$\psi_0, \dots, \psi_L : \OA{m} \to B$ as in Lemma 3.6, using
$\ep/7$ for $\ep$ and with $\delta$ chosen so small that it works
in Lemma 3.3 for the choice $\ep/7$ for $\ep$. Define $v \in M_{4L} (B)$
by
\[
 v = v_0^* \oplus v_0 \oplus v_0^* \oplus v_0 \oplus  \cdots
 \oplus v_{L-1}^* \oplus v_{L-1} \oplus v_{L-1}^* \oplus v_{L-1}
\]
(two copies of everything, except no $v_L$ or $v_L^*$).
Further define $\psi: \OA{m} \to M_{4L} (B)$ by
\beqr
\lefteqn{ \psi(a) = \psi_0 (a) \oplus  \psi_0 (a) \oplus
        \psi_0 (a) \oplus \psi_0 (a) \oplus  \cdots } \\
 & & \hspace{5em} \cdots
 \oplus  \psi_{L-1} (a) \oplus  \psi_{L-1} (a) \oplus
      \psi_{L-1} (a) \oplus  \psi_{L-1} (a)
\eeqr
(four copies of everything, except no $\psi_L (a)$).
Apply Lemma 3.3 to the unitary
\[
v_0 \oplus v_0^* \oplus v_1 \oplus v_1^* \oplus  \cdots
 \oplus v_{L-1} \oplus v_{L-1}^* \in U_0 (M_{2L} (B))
\]
and the \hm\   %
\[
a \mapsto \psi_0 (a) \oplus  \psi_0 (a) \oplus  \psi_1 (a) \oplus
     \psi_1 (a) \oplus  \cdots
 \oplus  \psi_{L-1} (a) \oplus  \psi_{L-1} (a)
\]
from $\OA{m} $ to $M_{2L} (B)$, and conjugate everything by a suitable
permutation matrix, to get a \hm\  $\eta_0 : \SO{m} \to
M_{4L} (B)$ such that:

(1) $\eta_0 (u \otimes 1) $ has finite spectrum.

(2) $\| v - \eta_0 (u \otimes 1) \| < \ep/7$.

(3) $\| \psi(s_j) - \eta_0 (1 \otimes s_j) \| < \ep/7$.

Write
\[
\eta_0 (f \otimes a) = \sum_{i = 1}^l f(x_i^{(0)} ) \eta_i (a)
\]
where
${\rm sp}(\eta_0 (u \otimes 1)) = \{ x_0^{(0)} , \dots, x_l^{(0)}  \}$
and $ \eta_1, \dots, \eta_l$ are homomorphisms from $\OA{m}$ to
$M_{4L} (B)$. Since ${\rm sp}(v) \subset X$, condition (2) above
implies that for each $i$ there is $x_i \in X$ with
$| x_i - x_i^{(0)}| < \ep/7$. Define
\[
\eta (f \otimes a) = \sum_{i = 1}^l f(x_i) \eta_i (a).
\]
Then (1) -- (3) above still hold, but with $\ep/7$ replaced by
$2 \ep/7$.

Choose \hm s $\mu_i : \OA{m} \to M_{4L} (B)$ with mutually orthogonal
ranges such that $[\mu_i] = - [\eta_i]$ in $KK^0 (\OA{m}, B)$. (For
example, let $\sigma: \OA{m} \to \OA{m}$ be a homomorphism whose class
in $KK^0 (\OA{m}, \OA{m})$ is $-[\id]$; such a homomorphism can be
constructed by using Theorem 3.1 of \cite{Rr2} with
$D = \OA{m} \otimes K$.
Then set $\mu_i = \eta_i \circ \sigma$.)
Then define $\mu: C(X) \otimes \OA{m}$ by
\[
\mu (f \otimes a) = \sum_{i = 1}^l f(x_i) \mu_i (a).
\]
Since $\varphi_0$ and $\varphi_1$ are \aab , we have
\[
\varphi_0 \aueeps{\ep/7} \varphi_0 \tdsum (\eta \oplus \mu)  \andeqn
\varphi_1 \aueeps{\ep/7} \varphi_1 \tdsum (\eta \oplus \mu).
\]
To complete the proof, we therefore show that
\[
\varphi_0 \oplus \eta \aueeps{5\ep/7} \varphi_1 \oplus \eta.
\]

Let $U$ be the permutation matrix
\[
U = {\rm diag}\left(
\left( \begin{array}{cc} 0 & 1 \\ 1 & 0 \end{array} \right),
\dots,
\left( \begin{array}{cc} 0 & 1 \\ 1 & 0 \end{array} \right)
\right).
\]
Note that $U v U^* = v^*$ and $U \psi(a) U^* = \psi(a)$.  Define
 $\tilde{\eta}(b) = U \eta(b) U^*$ for $b \in C(X) \otimes \OA{m}$.
Then $\varphi_1 \oplus \eta$ is unitarily equivalent
 to $\tilde{\eta} \oplus \varphi_1$,
so we have to show
$ \varphi_0 \oplus \eta \aueeps{5\ep/7} \tilde{\eta} \oplus \varphi_1$.
Now
\begin{eqnarray*}
\lefteqn{ \| (\varphi_0 \oplus \eta) ( u \otimes 1)   -
        (\tilde{\eta} \oplus \varphi_1) ( u \otimes 1)\| } \\
 & & \leq \| (\varphi_0 \oplus \eta) ( u \otimes 1) - v_0 \oplus v \| +
          \| v_0 \oplus v - v^* \oplus v_L \| \\
 & & \,\,\,\, \,\,\,\, \mbox{} + \| v^* \oplus v_L -
                     (\tilde{\eta} \oplus \varphi_1) (u \otimes 1) \| \\
 & & < 2\ep/7 + \ep/7 + 2\ep/7 = 5\ep/7.
\end{eqnarray*}
Similarly,
\begin{eqnarray*}
\lefteqn{ \| (\varphi_0 \oplus \eta) ( 1 \otimes s_j)   -
        (\tilde{\eta} \oplus \varphi_1) ( 1 \otimes s_j)\| } \\
 & & \leq \| (\varphi_0 \oplus \eta) ( 1 \otimes s_j) -
                        \psi_0 (s_j) \oplus \psi(s_j)  \| +
          \| \psi_0 (s_j) \oplus \psi(s_j) -
                      \psi(s_j) \oplus \psi_L (s_j)  \| \\
 & & \,\,\,\, \,\,\,\, \mbox{} + \| \psi(s_j) \oplus \psi_L(s_j) -
                   (\tilde{\eta} \oplus \varphi_1) (1 \otimes s_j) \| \\
 & & < 2\ep/7 + \ep/7 + 2\ep/7 = 5\ep/7.
\end{eqnarray*}
This shows that we do indeed have
$\varphi_0 \oplus \eta \aueeps{5\ep/7} \varphi_1 \oplus \eta$,
and completes the proof that $\varphi_0$ is \ayue\ to $\varphi_1$.

Now we do the unital case. Let $\iota : B \to M_2 (B)$
be the standard embedding in the upper left corner.
Then for any $\delta > 0$, we have
 $\iota \circ \varphi_0 \aueeps{\delta} \iota \circ \varphi_1$.
Let $U$ be an implementing unitary. Then
\[
\left\|U \left(\begin{array}{cc} 1 & 0 \\ 0 & 0 \end{array} \right) U^*-
\left( \begin{array}{cc} 1 & 0 \\ 0 & 0 \end{array} \right) \right\|
\]
is small.
Therefore there exists a unitary $V$ such that $\| U - V \|$ is small
and
\[
V \left( \begin{array}{cc} 1 & 0 \\ 0 & 0 \end{array} \right) V^* =
\left( \begin{array}{cc} 1 & 0 \\ 0 & 0 \end{array} \right).
\]
If $\delta$ is chosen small enough, then $V$ will implement an
\aue\ $\iota \circ \varphi_0 \aueeps{\ep} \iota \circ \varphi_1$.
Now $V$ must have the form ${\rm diag}(V_1, V_2)$ for
$V_1, V_2 \in U(B)$, and it follows that $V_1$ implements an
\aue\ $\varphi_0 \aueeps{\ep} \varphi_1$.
\QED

\vspace{0.6\baselineskip}

In order to handle \hm s which are not necessarily injective, we
introduce the following definition.

\vspace{0.6\baselineskip}

{\bf 3.8 Definition} Let
$X$ be a compact metric space, let  $D$ be a simple \CA ,
and let $\varphi:C(X)\otimes D \to C$ be a homomorphism to another
\CA\ $C$.
Then ${\rm ker}(\varphi) = C_0 (X \setminus Y) \otimes D$
for some closed subset $Y\subset X$.
Define
$$
\delta(\varphi)=\sup\{ {\rm dist}(Y,x): x\in X \}.
$$
Let $F$ be a finite subset of $A.$
 We regard elements of $C(X)\otimes D$ as continuous
 functions from $X$ to $D.$
Then we say that $\varphi$ is
$\ep$-{\it approximately injective with respect to} $F$
if, for every $f\in F,$
\[
\|f(x_1)-f(x_2)\|<\ep
\]
whenever ${\rm dist}(x_1,x_2)\le \delta(\varphi).$

Let $A=\bigoplus_{i=1}^k C(X_i)\otimes D_i,$ where each
$D_i$ is simple. Let $\varphi: A \to C$ be a homomorphism and
let $F$ be a finite subset of $A$.
Let $\pi_i : A \to  C(X_i)\otimes D_i$ be the projection on the $i$-th
summand.
We say $\varphi$ is $\ep$-approximately injective with respect
to $F$ if each $\varphi|_{C(X_i)\otimes D_i }$
 is $\ep$-approximately injective with respect to $\pi_i (F).$

\vspace{0.6\baselineskip}

{\bf 3.9 Lemma} Let $C$ be an even Cuntz-circle algebra, let $G$ be
a finite subset of $C,$ and let $B$ be a purely infinite simple
\CA.
Let $0 < \ep_0 < \ep$, and let
$\varphi: C\to B$ be \aab\  and $\ep_0$-approximately
injective with respect to $G$.
Then there is an injective homomorphism
$\varphi_0: C\to B$
with $[\varphi_0]=0$ in $KK^0 (C,B)$ such that
$\varphi\aueeps{\ep} \varphi\tdsum\varphi_0$  with respect to $G$.

\vspace{0.6\baselineskip}

{\it Proof:} By considering each summand separately, we
may assume, without loss of generality, that
$C=C(X)\otimes M_n (\OA{m}).$
Let $e\in B$ be a nonzero projection with $[e]=0$
in $K_0(B).$
Let $\varphi_0: C \to eBe$ be a unital injective homomorphism having
the same properties as $\vph$ in Lemma 1.15.
Let $\psi : C\to eBe$ be
a unital homomorphism as in the conclusion of Lemma 1.15, using
$(\ep - \ep_0)/2$ in place of $\ep$.
In particular, $\psi$ has the form
\[
\psi (f \otimes a) = \sum_{i = 1}^l f(x_i) \psi_i (a)
\]
for suitable $x_i \in X$ and $\psi_i : M_n (\OA{m}) \to B$, and
$$
\|\vph_0(g)-\psi(g)\|<(\ep- \ep_0 )/2
$$
for all $g\in G$.
Let $\ker(\vph)=C_0(X\setminus Y)\otimes M_n (\OA{m}).$ Since $\vph$
is $ \ep_0 $-approximately injective with respect to $G$,
there are $y_i\in Y\i X$
such that $ \|g(x_i)-g(y_i)\|< \ep_0$ for $i=1,2,...,l$ and
$g \in G$. Define $\psi' : C(X)\otimes M_n (\OA{m}) \to B$  by
\[
\psi' (f \otimes a) = \sum_{i = 1}^l f(y_i) \psi_i (a)
\]
for $f \in C(X)$ and $a \in M_n (\OA{m})$.
We can rewrite the formula for $\psi$ as
$\psi (g) = \sum_{i = 1}^l \psi_i (g(x_i))$ for
$g \in C(X)\otimes M_n (\OA{m})$, and similarly for $\psi'$
(using $y_i$). Therefore
$$
\|\psi(g)-\psi' (g)\| =
 \left\|\sum_{i=1}^l [\psi_i (g(x_i))-\psi_i (g(y_i))] \right\|< \ep_0 .
$$
(The terms in the sum are in orthogonal corners of $eBe$.)
Since $\vph$ is \aab, $\vph\aueeps{(\ep- \ep_0 )/2}\vph\tdsum\psi'$
with respect to $G$. Therefore  $\vph\aueeps{\ep} \vph\tdsum\vph_0$
with respect to $G$.  \QED

\vspace{0.6\baselineskip}

{\bf 3.10 Remark} It is important to note that the
\hm\  $\varphi\tdsum\varphi_0$ in the previous lemma is
\aab\  and injective.

\vspace{0.6\baselineskip}

The following result is the analog of Proposition 3.7 for approximately
injective \hm s.

\vspace{0.6\baselineskip}

{\bf 3.11 Proposition} Let $\varphi_0$
and $\varphi_1$ be two \aab\   homomorphisms from
 an even Cuntz-circle algebra
$C$ to a \pisca\  $B,$ with the same class in $KK^0 (C,B)$
and either both unital or both nonunital.
Let $F \subset C$ be a finite generating set, let $\ep > 0$,
and let $0 < \ep_0 < \ep / 2.$ Suppose that both
$\varphi_0$ and $\varphi_1$ are $\ep_0$-approximately
injective with respect to $F$.
Then $\varphi_0\aueeps{\ep} \varphi_1$
with respect to $F.$

\vspace{0.6\baselineskip}

{\it Proof:} Let $\ep_0 < \eta < \ep/2$.
By Lemma 3.9, there are injective \aab\   homomorphisms
$\psi_1,\,\psi_2: C\to B$ such that $[\psi_i]=[\vph_i]$
in $KK^0 (C,B)$, such that $\vph_i\aueeps{\eta} \psi_i$
with respect to $F$, and such that either both are unital or
both are nonunital.
Then Proposition  3.7 implies that
$\psi_0\aueeps{\ep-2\eta}\psi_1.$ Therefore $\vph_0\aueeps{\ep}\vph_1.$
\QED

\newpage

\section{ The existence theorem}

\vspace{\baselineskip}

The purpose of this section is to prove the following existence theorem,
to be used in the last section in the construction of approximate
intertwinings of direct systems.
We will also prove several related lemmas
that will be needed in the last section.
We will denote the Kasparov product of $\alpha$ and $\bt$ by
$\alpha \times \bt$ whenever it is defined. In particular, we use
this notation for $\alpha \in KK^* (A, B)$ and $\bt \in KK^* (B, C)$
(yielding $\alpha \times \bt \in KK^* (A, C)$), and
for $\alpha \in KK^* (A_1, B_1)$ and $\bt \in KK^* (A_2, B_2)$
(yielding
$\alpha \times \bt \in KK^* (A_1 \otimes A_2, B_1 \otimes B_2)$).

\vspace{0.6\baselineskip}

{\bf 4.1 Theorem} Let $A$ be an even \CSalg , and let $B = \dirlim B_k$
be a \DECSalg , with maps $\psi_k : B_k \to B$.
Let $\alpha \in KK^0 (A, B)$.
Then for every sufficiently large $k$, there exists a
permanently \aab\  homomorphism
$\varphi: A \to B_k$
such that the class $[\varphi] \in KK^0 (A, B_k)$
satisfies  $[\varphi] \times [\psi_k] = \alpha$.
Moreover, if $B$ is unital, and the Kasparov product
 $[1_A] \times \alpha$ is $[1_B]$, then $\varphi$ can be chosen to be
unital.

\vspace{0.6\baselineskip}

To keep down the size of some of the formulas (so that they fit on the
page), we will use the following notation throughout this section.

\vspace{0.6\baselineskip}

{\bf 4.2 Notation}  We denote by $S$ the \CA\ $C(S^1)$.

\vspace{0.6\baselineskip}

Before starting any $K$-theory calculations, we recall from \cite{Cu1}
that $K_0 (\OA{m}) \cong {\bf Z}/(m-1){\bf Z}$ and $K_1 (\OA{m}) = 0$.
Since we will make extensive use of the Universal Coefficient
Theorem \cite{RS}, we also note the following (well known) fact.

\vspace{0.6\baselineskip}

{\bf 4.3  Lemma} Let $X \subset S^1$ be compact.
Then $C(X) \otimes \OA{m}$ is in the bootstrap category
${\cal N}$ defined in \cite{RS} just before 1.17.

\vspace{0.6\baselineskip}

{\em Proof:} It is shown in 2.1
of \cite{Cu2} that $\OA{m}$ is stably isomorphic to a crossed
product of an AF algebra by an action of ${\bf Z}$. So $\OA{m}$ is in
${\cal N}$. Therefore so is $C(X) \otimes \OA{m}$.
(See 22.3.5 (d) and (f) and Chapter 23 in \cite{Bl2}.) \QED

\vspace{0.6\baselineskip}

{\bf 4.4 Lemma} Let $m,n \in {\bf N}, m,n \geq 2$, and let
$\alpha \in KK^0 (\OA{m}, \OA{n})$.
Then there exists a homomorphism $\varphi: \OA{m} \to \OA{n}$ such that
$[\varphi] = \alpha$.

\vspace{0.6\baselineskip}

{\it Proof:} Theorem 3.1 of
\cite{Rr2}  provides $\lambda: \OA{m} \to {\cal K}
\otimes \OA{n}$ such that $[\lambda] = \alpha$.
Let $e_{11} \in {\cal K}$ be the
projection on the first standard basis vector, and choose a partial
isometry $v \in {\cal K} \otimes \OA{n}$
such that $vv^* = \lambda (1)$ and
$v^*v \le e_{11} \otimes 1$. Then $\varphi(a) = v^*\lambda(a)v$,
regarded as a homomorphism from
$\OA{m}$ to $(e_{11} \otimes 1)({\cal K} \otimes \OA{n})
(e_{11} \otimes 1) \cong \OA{n}$, satisfies $[\varphi] = \alpha$.
\QED

\vspace{0.6\baselineskip}

{\bf 4.5 Remark}
It is also easy to construct $\varphi$ directly, since the Universal
Coefficient Theorem (\cite{RS}, Theorem 1.17) implies that the natural
map $\gamma: KK^0 (\OA{m}, \OA{n}) \to \Hom(\KO{0}{m}, \KO{0}{n})$ is an
isomorphism.

\vspace{0.6\baselineskip}

The next two steps are to show that elements of
$ KK^0 (\OA{m}, \sO{n})$  and $KK^0 (\sO{m}, \OA{n})$
are representable by homomorphisms (in the second case, assuming
$n$ even). In both steps, the computation of some of the $KK$-classes
 will
be done by reduction to the identification of elements of
$\Ext({\cal O}_A, B)$ in \cite{Cu4}. We therefore introduce appropriate
notation.

\vspace{0.6\baselineskip}

{\bf 4.6 Notation}
For any \CA\  $B$,
let $Q(B)$ denote the stable outer multiplier algebra
$M({\cal K} \otimes B)/({\cal K} \otimes B)$. In particular,
 let $Q = Q({\bf C})$ be the Calkin algebra.
If
\[ 0 \arrow I \arrow B \arrow C \arrow 0 \]
is a short exact sequence of \CA s, let
 $ {\rm Ind}_i : K_i(C) \to K_{1-i}(I)$ be the
connecting homomorphism in the associated six term exact sequence in
$K$-theory. We write $ {\rm Ind}_1 (u)$ for $ {\rm Ind}_1 ([u])$
for unitaries $u$, and similarly for classes of projections.
We extend the definition of $ {\rm Ind}_1 $ to partial isometries
with the same initial and final projections by setting
$ {\rm Ind}_1 (u) = {\rm Ind}_1 ([u + (1-p)])$
when $uu^* = u^*u = p$.

\vspace{0.6\baselineskip}

{\bf 4.7 Definition} (Compare Cuntz \cite{Cu3}, Section 3.)
Let $A$ be an $m \times m$
matrix with entries in $\{0,1\}$, satisfying the
condition (I) of \cite{CK}, and let
 ${\cal O}_A$ be the corresponding Cuntz-Krieger algebra.
Call its canonical generating partial isometries $s_1, \dots , s_m$,
and set $p_j = s_j s_j^*$.
 Let $B$ be a \CA , and let $\sigma$ and $\tau$ be two extensions of
 ${\cal O}_A$ by $B$, regarded as homomorphisms from
 ${\cal O}_A$ to  $Q(B)$. Suppose that $\sigma (p_j) = \tau(p_j)$
for all $j$. Define
\[ d_{\sigma , \tau} \in K_0 (B)^m/(1-A)K_0 (B)^m \]
to be the image there of
\[ ({\rm Ind}_1 (\sigma (s_1) \tau(s_1)^*), \dots ,
     {\rm Ind}_1 (\sigma (s_m) \tau(s_m)^*) . \]

\vspace{0.6\baselineskip}

{\bf 4.8 Theorem} There is an isomorphism
\[ d: \Ext( {\cal O}_A , B) \to K_0 (B)^m/(1-A)K_0 (B)^m, \]
determined as follows: If $\sigma$ and $\tau$ are as in the previous
definition, and $\tau$ lifts to a homomorphism
$\tilde{\tau} : {\cal O}_A \to M({\cal K} \otimes B)$
such that the projections
$\tilde{\tau}(p_j)$ (for $j = 1, \dots , m$) and  $1 - \tilde{\tau} (1)$
are all Murray-von Neumann equivalent to $1$ in  $M(K \otimes B)$, then
$d([\sigma]) = d_{\sigma , \tau}$.

\vspace{0.6\baselineskip}

{\it Proof:} See Section 3 of \cite{Cu3}. \QED

\vspace{0.6\baselineskip}

{\bf 4.9 Lemma}
 The formula for $d([\sigma])$ in the previous theorem remains
valid if it is merely assumed that $\sigma$ and $\tau$ are related as in
Definition 4.7, and that $[\tau] = 0$ in $KK^1 ( {\cal O}_A , B)$.

\vspace{0.6\baselineskip}

{\it Proof:} Let $\tau_0  :{\cal O}_A \to  Q(B)$
be an absorbing trivial
 extension (\cite{Ks}, Section 7, Definition 2; note that it
is also shown in the same section that $\tau_0$ exists).
Then $\sigma \oplus \tau_0$ and $\tau \oplus \tau_0$  satisfy the
hypotheses of the previous theorem. It is trivial to check that
\[
{\rm Ind}_1 ((\sigma \oplus \tau_0)(s_j) (\tau \oplus \tau_0)(s_j)^*) =
     {\rm Ind}_1 (\sigma (s_j) \tau(s_j)^*).
\]
Therefore
\[ d([\sigma]) = d([\sigma \oplus \tau_0]) =
 d_{\sigma \oplus \tau_0, \tau \oplus \tau_0} = d_{\sigma , \tau}.
\,\,\,  \QED \]

{\bf 4.10 Lemma} Let $m,n \in {\bf N} , m,n \geq 2$, and let
$\alpha \in KK^0 (\OA{m}, \sO{n})$.
Then there exists a homomorphism $\varphi: \OA{m} \to \sO{n}$ such that
$[\varphi] = \alpha$.

\vspace{0.6\baselineskip}

{\it Proof:} The Universal Coefficient Theorem (\cite{RS}, Theorem 1.17)
yields a short exact sequence
\beqr
0 \arrow \ext{ \KO{0}{m}}{ \KsO{1}{n}} & \longrightarrow &
KK^0 (\OA{m}, \sO{n}) \\
 &  \stackrel{\gamma}{ \longrightarrow} & \Hom(\KO{0}{m}, \KsO{0}{n})
 \arrow 0 .
\eeqr
(The two missing terms are both zero, since $\KO{1}{m} = 0$.)

Let $H$ be the set of classes in $KK^0 (\OA{m}, \sO{n})$
of homomorphisms
from $\OA{m}$ to $\sO{n}$. We claim that $H$ is a subgroup of
$KK^0 (\OA{m}, \sO{n})$.
It follows from the exact sequence above
 that $KK^0 (\OA{m}, \sO{n})$ is a finite group, and
therefore it suffices to prove that $H$ is nonempty (which is trivial)
and closed under addition. So let $\varphi$ and $\psi$ be homomorphisms
from $\OA{m}$ to $\sO{n}$. Let $v, w \in \OA{n}$ be isometries with
orthogonal ranges. Then

\[ \sigma (a) = (1 \otimes v)\varphi(a)(1 \otimes v^*) +
(1 \otimes w)\psi(a)(1 \otimes w^*) \]
defines a homomorphism from $\OA{m}$ to $\sO{n}$ whose class in
$KK(\OA{m}, \sO{n})$ is $[\varphi] + [\psi]$.

We now show that $\gamma (H) =  \Hom(\KO{0}{m}, \KsO{0}{n})$.
Let $\iota : \OA{n} \to \sO{n}$ be the map $\iota (a) = 1 \otimes a$.
Then $\iota_*$ is an isomorphism from $\KO{0}{n}$ to $ \KsO{0}{n}$.
Given $\alpha \in \Hom(\KO{0}{m}, \KsO{0}{n})$, we get $\iota_*^{-1}
\circ \alpha \in \Hom(\KO{0}{m}, \KO{0}{n})$. Lemma 4.4 and Remark 4.5
provide a homomorphism $\varphi : \OA{m} \to \OA{n}$
such that $\varphi_* =  \iota_*^{-1}\circ \alpha$. Then
$\gamma ([\iota \circ \varphi ]) = \alpha$.

To complete the proof, it now suffices to show that $H$ contains
$ \ext{ \KO{0}{m}}{ \KsO{1}{n}} $.
Let $[\sigma_0]$ be the standard generator of
$KK^1 (S, {\bf C} )$,
namely the class of the extension given by the $C^*$-algebra of the
unilateral shift. (Thus, $\sigma_0$ is to be regarded as the
homomorphism from $S$ to $Q$ which sends the standard unitary in $S$
to the image in $Q$ of the unilateral shift.)
 Let $[\sigma]$ be its image in $KK^1 (\sO{n}, \OA{n})$.
Note that $[\sigma_0]$ has a left inverse $[\tau_0]$ in
$KK^1 ({\bf C}, S ) = K_1 (S)$.
Therefore $[\sigma]$ has a left inverse $[\tau]$. Hence the Kasparov
product with $[\sigma]$ defines  surjective homomorphisms from
$\KsO{j}{n}$ to $\KO{1-j}{n}$. Since $\KsO{1}{n}$ and $\KO{0}{n}$
are finite groups with the same cardinality, Kasparov product
with $[\sigma]$ is actually an isomorphism. Since the Universal
Coefficient Theorem is natural with respect to Kasparov products
(see \cite{RS}), we obtain the following
commutative diagram with exact columns, in which the horizontal maps are
induced by Kasparov product on the right with $[\sigma]$, and
the top horizontal map is an isomorphism:
$$\begin{array}{ccc}
0 & & 0 \\
\downarrow & & \downarrow \\
 \ext{ \KO{0}{m}}{ \KsO{1}{n}} &
                     \longrightarrow &  \ext{ \KO{0}{m}}{ \KO{0}{n}} \\
\downarrow & & \downarrow \\
 KK^0 ( \OA{m},\sO{n}) & \longrightarrow &    KK^1 ( \OA{m},\OA{n})  \\
\,\,\,\,\, \downarrow \gamma &  & \,\,\,\,\,\,\,  \downarrow \gamma' \\
\Hom(\KO{0}{m},\KsO{0}{n}) & \longrightarrow
                                       & \Hom(\KO{0}{m},\KO{1}{n})   \\
\downarrow & & \downarrow \\
0 & & 0
\end{array} $$

Since $\KO{1}{n} = 0$, the map from
$ \ext{ \KO{0}{m}}{ \KsO{1}{n}} $ to
$ KK^1 ( \OA{m},\OA{n}) $ is an
isomorphism. Thus, we need to find homomorphisms
$\varphi: \OA{m} \to \sO{n}$ which induce the zero map on $K_0$ and such
that the classes $[\varphi] \times [\sigma]$
exhaust $ KK^1 ( \OA{m},\OA{n}) $.
To do this, we identify  $ KK^1 ( \OA{m},\OA{n}) $ with
 $\Ext ( \OA{m},\OA{n})$, and use Theorem 4.8 and Lemma 4.9.

Let $A$ be the $m \times m$ integer matrix with all
entries equal to 1. Then ${\cal O}_A$ is
just $\OA{m}$.
Let $s_1, \dots, s_m \in {\cal O}_A$
be the canonical generating isometries, and let $p_j = s_j s_j^*$, as
in Definition 4.7.
Choose $m$ nonzero
orthogonal projections $q_1, \dots, q_m \in \OA{n}$ such that
the $K_0$ class of each $q_j$ is zero, and such that the projection
$q = q_1 + \cdots + q_m$ is strictly less than $1$. Then the $K_0$ class
of $q$ is also zero. Since $\OA{n}$ is purely infinite, there are
partial isometries $t_1, \dots, t_m \in \OA{n}$ such that
$t_j t_j^* =q_j$ and $t_j^* t_j = q$.
Define a homomorphism $\varphi_0 : \OA{m} \to \sO{n}$ by
$\varphi_0 (s_j) = 1 \otimes t_j$.

Now let $\eta_1, \dots, \eta_m \in \KO{0}{n}$. We construct a
homomorphism
 $\varphi: {\cal O}_A \to \sO{n}$ which induces the zero map on
 $K_0$ and such that
\[ d([\varphi] \times [\sigma]) =
    (\eta_1, \dots, \eta_m) + (1-A)K_0 (\OA{n})^m. \]
Taking an idea from \cite{Cu4},
 let $u$ be the canonical unitary generator
of $S$, choose projections $e_j \le q_j$ such that
$[e_j] =- \eta_j$
in $\KO{0}{n}$, and define
$\varphi(s_j) = [u \otimes e_j + 1 \otimes (q_j - e_j)]t_j$.

The Kasparov product
 $[\varphi] \times [\sigma]$ is represented by the homomorphism
$\sigma \circ \varphi = (\sigma_0 \otimes {\rm id}) \circ \varphi$,
put together as follows:
\[  {\cal O}_A \stackrel{\varphi}{\arrow}  \sO{n}
   \stackrel{\sigma_0 \otimes {\rm id}}{\longrightarrow}
     Q \otimes \OA{n}
 \hookrightarrow Q(\OA{n}). \]
To apply Lemma 4.9, we need a comparison extension whose class is
trivial. To obtain it, we merely use $\varphi_0$ in place of $\varphi$.
Note that in fact  $(\sigma_0 \otimes {\rm id}) \circ \varphi_0$
lifts.  Then
\[ [\varphi] \times [\sigma] =
d_{(\sigma_0 \otimes {\rm id}) \circ \varphi,
  (\sigma_0 \otimes {\rm id}) \circ \varphi_0} =
 (d_1, \dots, d_m) + (1-A)\KO{0}{n}^m, \]
where
\[ d_j = {\rm Ind}_1 ((\sigma \circ \varphi)(s_j)
                  ( \sigma \circ \varphi_0) (s_j^*))
   ={\rm Ind}_1  (\sigma_0 (u) \otimes e_j + 1 \otimes (q_j - e_j))
   = -[e_j] = \eta_j. \]
(The minus sign appears because the unilateral shift $\sigma_0 (u)$ has
index $-1$. The index is most easily computed in $Q \otimes \OA{n}$.) So
\[ d([\varphi] \times [\sigma]) =
             (\eta_1, \dots, \eta_m) + (1-A)\KO{0}{n}^m, \]
as desired.  \QED

\vspace{0.6\baselineskip}

{\bf 4.11 Lemma} Let $m,n \in {\bf N} , m,n \geq 2$, and assume $n$
is even. Let
$\alpha \in KK^0 (\sO{m}, \OA{n})$.
Then there exists a homomorphism $\varphi: \sO{m} \to \OA{n}$ such that
$[\varphi] = \alpha$.

\vspace{0.6\baselineskip}

{\it Proof:} The Universal Coefficient Theorem now
yields a short exact sequence
\beqr
0 \arrow \ext{ \KsO{1}{m}}{ \KO{0}{n}} & \longrightarrow  &
KK^0 (\sO{m}, \OA{n}) \\
  & \stackrel{\gamma}{ \longrightarrow} & \Hom(\KsO{0}{m}, \KO{0}{n})
 \arrow 0.
\eeqr
(The two missing terms are again both zero.)

As before, let $H$ be the set of classes in $KK^0 (\sO{m}, \OA{n})$
of homomorphisms
from $\sO{m}$ to $\OA{n}$. It follows, just as in the proof of
the previous lemma, that $H$ is a subgroup of $KK^0 (\sO{m}, \OA{n})$.
Furthermore, one checks that $\gamma (H) =  \Hom(\KsO{0}{m}, \KO{0}{n})$
in essentially the same way as there, composing maps from $\OA{m}$ to
 $\OA{n}$ on the right with a point evaluation map from $\sO{m}$ to
 $\OA{m}$ instead of on the left with  $\iota : \OA{n} \to \sO{n}$.

It remains to show that $H$ contains
$\ext{ \KsO{1}{m}}{ \KO{0}{n}} $.
We will reduce this proof to the result of the previous lemma.

We begin by constructing a homomorphism
$\omega : S \otimes \sO{n} \to \OA{n}$ which sends the Bott element
to a generator of $\KO{0}{n}$.
We take the Bott element to be the class
$b=[u] \times [u] \in K_0 (S \otimes S) = K_0 (C(S^1 \times S^1))$,
 obtained as the product
of two copies of the $K_1$-class of the standard unitary generator
$u$ of $S$. If $B$ is a $C^*$-algebra,  we then refer for
convenience to $b \times [1]$ as the Bott class in
$K_0 (S \otimes S \otimes B)$.

  The construction is essentially done in the proof
of Proposition 4.2 of \cite{Lr3}. Assume $n > 2$. (Otherwise
$KK^0 (\sO{m}, \OA{n}) = 0$, and there is nothing to prove.)
 Let $B$ be the unital AF algebra whose ordered
$K_0$-group is the group $G_{n-2}$ (defined before Corollary 3.5 in
\cite{Lr3}) and such that the class of the identity is the image of
the element $(1,1)$
in the first term of the direct limit. Let $\eta : S \otimes S
\to B$ be the unitization of the map in \cite{Lr3}, gotten from Theorem
7.3 of \cite{EL}. It sends the Bott element $b \in K_0 (S \otimes S)$
to $(1, -1)$. As in \cite{Lr3}, we tensor this map with the identity on
$\OA{n}$, observe that (by  \cite{Rr1}) $B \otimes \OA{n} \cong \OA{n}$,
and calculate $K$-theory to see that $b \times [1]$ goes to a
generator of $\KO{0}{n}$.

The automorphisms of ${\bf Z}/(n-1){\bf Z}$ act transitively on the
generators. We can therefore compose this homomorphism with a suitable
homomorphism from $\OA{n}$ to $\OA{n}$, so as to obtain a homomorphism
$\omega : S \otimes \sO{n} \to \OA{n}$ which sends the Bott element
to the standard generator $[1]$ of $\KO{0}{n}$.

Let $\sigma_0 : S \to Q$ send the standard unitary to the unilateral
shift, as in the proof of Lemma 4.10. Thus $[\sigma_0]$ is a generator
of $KK^1(S,{\bf C})$, and has a left inverse
 $[\tau_0] \in KK^1 ({\bf C}, S)$, for some homomorphism
 $\tau_0 :{\bf C} \to Q(S)$.
(We will make an explicit choice for $\tau_0$ below.)
Let $\sigma = \sigma_0 \otimes {\rm id}_{\OA{n}}$, as before. Let
$\tau = \tau_0  \otimes {\rm id}_{\OA{m}}$. (Note that we use $m$
 instead of $n$.)
Naturality of the Universal Coefficient Theorem now gives the following
commutative diagram with exact columns, in which the horizontal maps are
Kasparov product on the left with $[\tau]$:
\[ \begin{array}{ccc}
0 & & 0 \\
\downarrow & & \downarrow \\
 \ext{ \KsO{1}{m}}{ \KO{0}{n}} & \longrightarrow &
                                   \ext{ \KO{0}{m}}{ \KO{0}{n}}  \\
\downarrow & & \downarrow \\
KK^0 ( \sO{m},\OA{n}) & \longrightarrow & KK^1 ( \OA{m},\OA{n})  \\
\,\,\,\,\, \downarrow \gamma &  & \,\,\,\,\,\,\,  \downarrow \gamma' \\
\Hom(\KsO{0}{m},\KO{0}{n}) & \longrightarrow &
                                    \Hom(\KO{1}{m},\KO{0}{n})     \\
\downarrow & & \downarrow \\
0 & & 0
\end{array} \]
Left multiplication by $\tau$ is an isomorphism from $\ker(\gamma)$ to
$KK^1(\OA{m}, \OA{n})$, by an argument similar to the one in Lemma 4.10.
To complete the proof of the lemma, it therefore suffices to let
$\beta \in KK^1 (\OA{m}, \OA{n})$, and find $\psi : \sO{m} \to \OA{n}$
which is zero on $K$-theory and
such that $[\tau] \times  [\psi] = \beta$.

As in the proof of Lemma 4.10, let
$\varphi: \OA{m} \to \sO{n}$ be a homomorphism  which is zero on
$K$-theory and such that $[\varphi] \times [\sigma] = \beta$.
Further let
$\varphi_0$ be the comparison homomorphism (with trivial class in
$KK$-theory) used in that proof.
Then define $\psi = \omega \circ ({\rm id}_{S} \otimes \varphi)$
and $\psi_0 = \omega \circ ({\rm id}_{S} \otimes \varphi_0)$.
We claim this $\psi$ works. To do this, we must show that
\[ [\tau] \times  ([{\rm id}_{S}] \otimes [\varphi]) \times  [\omega]
     = \beta .\]
It is well known that there is a rank one projection
 $e \in M_2 \otimes S \otimes S$ which represents the
class $b + [1_{S \otimes S}]$.
Then (depending on sign conventions) we may take  $\tau_0$
to be the extension given by the composite
\[ {\bf C} \arrow M_2 \otimes S \otimes S
\stackrel{{\rm id} \otimes \sigma_0 \otimes {\rm id}_S}{\longrightarrow}
    M_2 \otimes Q \otimes S
    \stackrel{\cong}{\longrightarrow}  Q \otimes S
     \hookrightarrow Q(S), \]
in which the first map sends $1$ to $e$, and the third map is induced
by an isomorphism $M_2 \otimes Q \cong Q$.
Then $\tau$ is the tensor product of this with $\id_{\OA{m}}$,
which is a homomorphism from
$\OA{m}$ to $Q(\sO{m})$.

To compute the Kasparov pairing with $\psi$, we compose on the left
with $Q(\psi)$. This composition can be slightly rearranged to give:
\[ \OA{m} \arrow M_2 \otimes Q \otimes \sO{m} \arrow
  M_2 \otimes Q \otimes S \otimes \sO{n}
   \arrow  M_2 \otimes Q \otimes \OA{n} \]
\[ \hspace{1.5in} \hookrightarrow  M_2 \otimes Q(\OA{n})
   \stackrel{\cong}{\longrightarrow} Q(\OA{n}). \]
{}From now on, we will take $\tau_0$ and $\tau$ to be maps to
$M_2 \otimes Q \otimes S$
 and $M_2 \otimes Q \otimes \sO{m}$. Then
the first map is $\tau$, the second is
 ${\rm id}_{Q \otimes S} \otimes \varphi$, and the third is
 ${\rm id}_Q \otimes \omega$.
As in the proof of Lemma 4.10, we use Lemma 4.9 to compute
 $d([ Q(\psi) \circ \tau] \in \KO{0}{n}^m/(1-A)\KO{0}{n}^m$,
 where $A$ is an $m \times m$ matrix of $1$'s.
Let $(\eta_1, \dots, \eta_m)$, with $\eta_j \in \KO{0}{n}$,
 be a representative of the image of
$\beta$ in this group. Since
$[ Q(\psi_0) \circ \tau] = 0$, it suffices to compute
\[{\rm Ind}_1 ((Q(\psi) \circ \tau)(s_j)
   (Q(\psi_0) \circ \tau)(s_j)^*),\]
where $s_1, \dots, s_m$ are the generating isometries of $\OA{m}$.
Since ${\rm Ind}_1$ is natural,  the expression above is equal to
\beqr
\lefteqn{
 \omega_* ({\rm Ind}_1
 ([(({\rm id} \otimes \varphi) \circ \tau)(s_j)
 (({\rm id} \otimes \varphi_0) \circ \tau)(s_j)^*])) } \\
  & & \hspace{8em}  =
 \omega_* ({\rm Ind}_1
 ([(\tau_0 \otimes \varphi)(s_j)(\tau_0 \otimes \varphi_0)(s_j)^*]).
\eeqr
The elements whose indices we want are
$\tau_0 (e) \otimes \varphi (s_j)^*\varphi_0 (s_j)$.
Since ${\rm Ind}_i$   respects tensor products, the required index is
\[ {\rm Ind}_0 (\tau_0 (e)) \times
 [\varphi (s_j)\varphi_0 (s_j)^* + (1-\varphi (s_j s_j^*)]. \]

Now recall that $b =  [u] \times [u]$.  Let $s \in Q$ be the
unilateral shift. Then from $[e] = 1 + b$ we get
\[  {\rm Ind}_0 (\tau_0 (e)) = {\rm Ind}_0 ([1] \times [1]
      + [s] \times [u])
   =  {\rm Ind}_0 ([1]) \times [1] + {\rm Ind}_1 ( [s]) \times [u]
   = -[u]. \]
Also, $\varphi$ was constructed so as to have
\[  [\varphi (s_j)\varphi_0 (s_j)^* + (1-\varphi (s_j s_j^*)]
  = [u] \times (-\eta_j). \]
(See the proof of Lemma 4.10.) Therefore
\beqr
\lefteqn{
{\rm Ind}_1 ([(Q(\psi) \circ \tau)(s_j)(Q(\psi_0) \circ \tau)(s_j)^*])
             } \\
 & & \hspace{5em}  = \omega_* (-[u] \times [u] \times (-\eta_j))
  = \omega_* (b \times \eta_j) = \eta_j.
\eeqr
Since $d(\beta)$ is the image in $\KO{0}{n}^m/(1-A)\KO{0}{n}^m$ of
$(\eta_1, \dots , \eta_m)$, we have shown that
$d([ Q(\psi) \circ \tau])$ is equal to the image of $\beta$, as desired.
\QED

\vspace{0.6\baselineskip}

{\bf 4.12 Theorem} Let $m,n \in {\bf N} , m,n \geq 2$, and assume that
$n$ is even. Let
$\alpha \in KK^0 (\sO{m}, \sO{n})$.
Then there exists a homomorphism $\varphi: \sO{m} \to \sO{n}$ such that
the class $[\varphi] \in  KK^0 (\sO{m}, \sO{n})$ is equal to $\alpha$.

\vspace{0.6\baselineskip}

{\it Proof:} This time, the Universal Coefficient Theorem
yields the following short exact sequence, which we write vertically
so that it will fit on the page:
\[
\begin{array}{c}
0 \\
\downarrow \\
\ext{ \KsO{0}{m}}{ \KsO{1}{n}} \oplus \ext{ \KsO{1}{m}}{ \KsO{0}{n}}  \\
\downarrow \\
KK^0 (\sO{m}, \sO{n})   \\
\,\,\,\,\, \downarrow \gamma \\
\Hom(\KsO{0}{m}, \KsO{0}{n}) \oplus \Hom(\KsO{1}{m}, \KsO{1}{n})  \\
\downarrow \\
0
\end{array}
\]
Let $H$ be the set of classes in $KK^0(\sO{m}, \sO{n})$
of homomorphisms
from $\sO{m}$ to $\sO{n}$. Then $H$ is a subgroup of
$KK^0 (\sO{m}, \sO{n})$ for the same reason as in the proof of
Lemma 4.10.

\par
We observe that $\gamma (H)$ contains
$\Hom(\KsO{0}{m}, \KsO{0}{n}) \oplus  0$, by considering homomorphisms
that factor as

\[ \sO{m} \longrightarrow \OA{m} \longrightarrow
  \OA{n} \longrightarrow \sO{n} , \]
where the first map is evaluation at some point of $S^1$, the second one
is taken from Lemma 4.4, and the third one is given by
$a \mapsto 1 \otimes a$. Furthermore, if
$\alpha \in \Hom( {\bf Z}/(m-1){\bf Z},  {\bf Z}/(n-1){\bf Z})$,
and $\varphi \in \Hom (\OA{m}, \OA{n})$ satisfies $\varphi_* = \alpha$,
then $\gamma ([{\rm id} \otimes \varphi]) = (\alpha,   \alpha)$.
Since the elements we have exhibited  as being in $\gamma (H)$ generate

\[\Hom(\KsO{0}{m}, \KsO{0}{n}) \oplus \Hom(\KsO{1}{m}, \KsO{1}{n}),\]
we have shown that $\gamma(H)$ is equal to this entire group.
It remains to show that $H$ contains the first term of the exact
sequence above.

Let
\[
\alpha \in \ext{ \KsO{0}{m}}{ \KsO{1}{n}} \subset KK^0 (\sO{m}, \sO{n}).
\]
Let $\ep : \sO{m} \to \OA{m}$
be evaluation at some point of $S^1$. Then $\ep_*$ is an
isomorphism on $K_0$. Therefore we can form
\[
(\ep^*)^{-1} (\alpha) \in  \ext{ \KO{0}{m}}{ \KsO{1}{n}}
 \subset KK^0 (\OA{m}, \sO{n}).
\]
Choose by Lemma 4.10 a homomorphism
$\varphi : \OA{m} \to \sO{n}$ such that
$[\varphi] = (\ep^*)^{-1} (\alpha)$. Then
$\varphi \circ \ep : \sO{m} \to \sO{n}$ satisfies
$[\varphi \circ \ep] = \alpha$.
So
$\ext{ \KsO{0}{m}}{ \KsO{1}{n}} \subset H$.

A similar argument using Lemma 4.11 and $\iota : \OA{n} \to \sO{n}$,
defined by $ \iota (a) = 1 \otimes a$, shows that
$\ext{ \KsO{1}{m}}{ \KsO{0}{n}} \subset H$.
This shows $H = KK^0 (\sO{m}, \sO{n})$.  \QED

\vspace{0.6\baselineskip}

We now prove results in which we impose conditions on the \hm s.

\vspace{0.6\baselineskip}

{\bf 4.13 Theorem} Let $X_1$ and $X_2$ be compact connected subsets
of $S^1$. (Thus, each is either a point, a closed arc, or all of $S^1$.)
Let $m_1, m_2 \geq 2$ be even integers, let $n_1, n_2$ be
integers, and let
\[
\alpha \in KK^0 (M_{n_1} \otimes C(X_1) \otimes \OA{m_1},
       M_{n_2} \otimes C(X_2) \otimes \OA{m_2}).
\]
Then for every nonzero projection
$p \in M_{n_2} \otimes C(X_2) \otimes \OA{m_2}$ satisfying
\[
[1_{M_{n_1} \otimes C(X_1) \otimes \OA{m_1}} ] \times \alpha = [p],
\]
there exists a permanently \aab\  \hm
\[
\varphi: M_{n_1} \otimes C(X_1) \otimes \OA{m_1} \to
       M_{n_2} \otimes C(X_2) \otimes \OA{m_2}
\]
such that $[\varphi] = \alpha$ and $\varphi(1) = p.$

\vspace{0.6\baselineskip}

{\em Proof:} Using Lemma 1.11, we can without loss of generality take
$p = 1$.

Let $e_{11}$ be a rank one projection in $M_{n_1}$, and set
$e = e_{11} \otimes 1 \otimes 1$. Then
$n_1 [e] \times \alpha = [1]$ in
$K_0 (M_{n_2} \otimes C(X_2) \otimes \OA{m_2}) \cong
       K_0 (M_{n_2} \otimes \OA{m_2})$.
Since $M_{n_2} \otimes \OA{m_2}$  is purely infinite and simple,
there exist $n_1$ Murray-von Neumann equivalent \mops\  in
$M_{n_2} \otimes \OA{m_2}$ with
$K_0$-classes equal to $[e] \times \alpha$ and which sum to $1$.
Let $p$ be one of these, and regard $p$ as a projection in
$M_{n_2} \otimes C(X_2) \otimes \OA{m_2}$. Then
\[
M_{n_2} \otimes C(X_2) \otimes \OA{m_2} \cong M_{n_1}
      (p [ M_{n_2} \otimes C(X_2) \otimes \OA{m_2} ] p).
\]
It suffices to construct a suitable \hm\  from $C(X_1) \otimes \OA{m_1}$
to $p [ M_{n_2} \otimes C(X_2) \otimes \OA{m_2} ] p$. Thus, without
loss of generality, we may assume $n_1 = 1$. Applying Lemma 1.11
again, we reduce again to the case that $p = 1$.

We now construct a \hm\  $\varphi_0$ such that $[\varphi_0] = \alpha$,
but without requiring that $\varphi_0 (1) = 1$ or that
$\varphi_0$ be permanently \aab .
For this step, it
suffices to construct a \hm\  to $C(X_2) \otimes \OA{m_2}$.
 If $X_1$ and $X_2$ are each either a point or $S^1$, then the
existence of the required \hm\  follows from Lemma 4.4, Lemma 4.10,
Lemma 4.11, or Theorem 4.12. If $X_1$ is a closed arc, we compose
the map of evaluation at some point of $X_1$ (which is a homotopy
equivalence) with a suitable map obtained from the case in which
$X_1$ is a point. If now $X_2$ is a closed arc, we compose on the
other side with the map from $\OA{m_2}$ to
$C(X_2) \otimes \OA{m_2}$ which sends each element of $\OA{m_2}$
to the corresponding constant function. (This map is also a
homotopy equivalence.)

We now observe that, just as in the proof of Lemma 1.11, Theorem B
of \cite{Zh3} shows that $\varphi_0 (1)$ is unitarily equivalent to
a constant projection. So we can assume it is a constant projection.
If $\varphi_0 (1) = 1$, we choose a constant proper isometry
$v \in  M_{n_2} \otimes C(X_2) \otimes \OA{m_2}$ (constant as a
function from $X_2$ to $M_{n_2} \otimes \OA{m_2}$), and replace
$\varphi_0$ by $a \mapsto v \varphi_0 (a) v^*$. So we can assume
$\varphi_0 (1)$ is a constant projection different from $1$.
Let its constant value be $f$.

Note that in $K_0 (M_{n_2} \otimes \OA{m_2})$ we have
\beqr
[1 - f] & = & ({\rm ev}_{x_0})_* ([1] - [\varphi_0 (1)])  \\
 & = &
   ({\rm ev}_{x_0})_* ( [1_{C(X_1) \otimes \OA{m_1}} ] \times \alpha  -
            [1_{C(X_1) \otimes \OA{m_1}} ] \times \alpha ) = 0.
\eeqr
Therefore Lemma 1.15 yields a unital permanently \aab\  \hm
\[
\varphi_1: C(X_1) \otimes \OA{m_1} \to
              (1-f) (M_{n_2} \otimes \OA{m_2}) (1-f)
\]
such that $[\varphi_1] = 0$ in
$KK^0 ( C(X_1) \otimes \OA{m_1},  M_{n_2} \otimes \OA{m_2})$.
Then we define
$\varphi (a) = \varphi_0 (a) + 1_{C(X_2)} \otimes \varphi_1 (a)$.
This is the required \hm .  \QED

\vspace{0.6\baselineskip}

{\bf 4.14 Lemma} Let $A$ be a \CSalg .  Then $KK^* (A, -)$
commutes naturally with countable direct limits.

\vspace{0.6\baselineskip}

{\em Proof:} The algebra $A$ is in the bootstrap category $\cal N$ of
\cite{RS} by Lemma 4.3. Clearly $K_* (A)$ is finitely generated (even
 finite). Therefore Proposition 7.13
of \cite{RS} implies that $KK^* (A, -)$
is an additive homology theory. The desired result now follows from
Section 5 of \cite{Sch}.  \QED

\vspace{0.6\baselineskip}

{\em Proof of Theorem 4.1:}  As in the statement of the theorem, let
\[
A = \bigoplus_{i = 1}^r M_{n(i)} \otimes C(X_i) \otimes \OA{m(i)}
\]
 be an even \CSalg , let $B = \dirlim B_k$ be a \DECSalg ,
with maps $\psi_k : B_k \to B$, and let
 $\alpha \in KK^0 (A, B)$. By the definition of a \CSalg ,
the spaces $X_i$ have only finitely many connected components.
Replacing each one by its connected components (and correspondingly
increasing the number of summands), we may assume each $X_i$ is
connected.
It follows from the
previous lemma that for every sufficiently large $k$,
there is $\alpha_0 \in KK^0(A, B_k)$
such that $\alpha = \alpha_0 \times [\psi_k]$.
If $B$ is unital, then we may assume the maps
$\psi_k$ are unital. Choose $\alpha_0$ as above for some
fixed $k$. Then $[1_A] \times \alpha_0$ and $[1_{B_{k}}]$
have the same image in $K_0 (B)$, and so also have the same image in
$K_0 (B_l)$ for all sufficiently large $l$. We replace $\alpha_0$
by its product with the map from $B_k$ to $B_l$.

Write
$B_k =
    \bigoplus_{j = 1}^{r'} M_{n'(j)} \otimes C(Y_j) \otimes \OA{m'(j)}$,
with the $Y_j$ compact connected subsets of $S^1$. Let $B_k^{(j)}$
be the $j$-th summand in this expression, and let $A^{(i)}$ be
the $i$-th summand of $A$. Write
$\alpha_0 = \sum_{i, j} \alpha_0^{(i, j)}$ with
$\alpha_0^{(i, j)} \in KK^0 (  A^{(i)}, B_k^{(j)} )$.
The K\"{u}nneth formula \cite{Sch0} shows that the map
\[
M_{n'(j)}  \otimes \OA{m'(j)} \to
      M_{n'(j)} \otimes C(Y_j) \otimes \OA{m'(j)},
\]
given by tensoring with $1_{C(Y_j)}$, is an isomorphism on $K_0$.
Since $M_{n'(j)}  \otimes \OA{m'(j)}$ is purely infinite and simple,
it follows that we can find nonzero \mops\  $p^{(i, j)} \in B_k^{(j)}$
such that $[p^{(i, j)}] = [1_{A^{(i)}}] \times \alpha_0^{(i, j)}$.
In the unital case, we have
$\sum_i [1_{A^{(i)}}] \times \alpha_0^{(i, j)} = [1_{B_k^{(j)}}]$,
and we can require that
$\sum_i p^{(i, j)} = 1_{B_k^{(j)}}$.
Now use Theorem 4.13 to choose permanently \aab\  \hm s
$\varphi^{(i, j)} : A^{(i)} \to B_k^{(j)}$ such that
$[\varphi^{(i, j)}] = \alpha_0^{(i, j)}$. Define
$\varphi = \bigoplus_j \sum_i \varphi^{(i, j)}$.
Since each $\varphi^{(i, j)}$
is permanently \aab , so is $\varphi$. Also, in the unital case
$\varphi$ is unital. \QED

\vfill
\newpage

\section{ The main results}

\vspace{\baselineskip}

In this section, we work with algebras in the following class.

\vspace{0.6\baselineskip}

{\bf 5.1 Definition}
Let $\Class$ be the class  of simple \CA s $A$ which are
direct limits $A \cong \dirlim A_k$, in which
each $A_k$ is an even \CSalg\  and each
map $A_k \to A$ is \aab .

\vspace{0.6\baselineskip}

Our main result (see Theorems 5.4 and 5.17)
is that algebras
$A \in \Class$ are classified up to isomorphism by the $K$-theoretic
 invariant
$(K_0 (A), [1_A], K_1 (A))$ in the unital case and
$(K_0 (A),  K_1 (A))$ in the nonunital case. (In the first of these
expressions, $[1_A]$ is the class in $K_0 (A)$ of the identity of $A$.)
The class $\Class$ contains all \CA s of the form $B \otimes \OA{m}$,
with $m$ even and
$B$ a simple \CA\  obtained as a direct limit of finite direct sums of
matrix algebras over $C(S^1)$ or $C([0,1])$.
In particular, it contains the tensor products of irrational rotation
algebras with even Cuntz algebras. These facts give us Corollaries
5.9 through 5.13. The class $\Class$ is also closed under the
formation of hereditary subalgebras, countable direct limits (provided
that the direct limit is simple), and
tensor products with simple AF algebras.

In Theorem 5.24, we give a classification theorem for direct limits in
which the building blocks are certain simple \CA s, but in which no
restriction is made on the maps associated with the direct systems.

If $A \in \Class$, then
$K_0 (A)$ and $K_1 (A)$  are countable abelian groups
in which every element has finite
odd order. We also show in this section that the class  $\Class$
is large enough that all possible values of $(K_0 (A), [1_A], K_1 (A))$,
with $K_0 (A)$ and $K_1 (A)$  countable odd torsion groups, are
realized by unital algebras $A \in \Class$, and that similarly all
such values of $(K_0 (A),  K_1 (A))$ are realized by nonunital
algebras $A \in \Class$. (See Theorems 5.26 and 5.27.)

We begin by establishing our notation for direct limits.

\vspace{0.6\baselineskip}

{\bf 5.2 Notation}  The notation
$A = \dirlim (A_k, \varphi_{k, k+1})$
will be taken to mean that $(A_k, \varphi_{k, k+1})_{k = 1}^{\infty}$
is a direct system of \CA s, with \hm s
$\varphi_{k, k+1}: A_k \to A_{k + 1}$ and direct limit
$A = \dirlim A_k$. We will further implicitly define
$\varphi_{k, l}: A_k \to A_{l}$, for $l \geq k$, to be the composite
$\varphi_{l-1, l} \circ \varphi_{l-2, l-1} \circ
     \cdots \circ \varphi_{k, k+1}$,
and $\varphi_{k, \infty} : A_k \to A$ to be the map to the direct limit
induced by the system.

A system of finite generating sets for the system
$(A_k, \varphi_{k, k+1})_{k = 1}^{\infty}$
consists of finite subsets $G_k \subset A_k$ such that $G_k$ generates
$A_k$ as a \CA\  and $\varphi_{k, k+1} (G_k) \subset G_{k+1}$ for
all $k$.
(In most cases of interest, each $A_k$ will be finitely generated,
and so such systems will exist.)

\vspace{0.6\baselineskip}

We now show that the \CA s in $\Class$ are purely infinite.

\vspace{0.6\baselineskip}

{\bf 5.3 Lemma } Let $A= \dirlim  (A_k, \varphi_{k,k+1})$ be a
direct limit of \CA s $A_k$, and assume that $A$ is simple.
If either

(1) all $A_k$ are even \CSalg s, or

(2) all $A_k$ are finite direct sums of \pisca s,

\noindent
then $A$ is a \pisca .

\vspace{0.6\baselineskip}

{\it Proof:}
(1) Proposition 7.7 of \cite{Rr1} implies that
 $\OA{m}$ is approximately divisible.
Therefore each  $A_k$ is approximately
divisible (see the remark after 1.4 of \cite{BKR}).
It is obvious from the definition of approximate divisibility \cite{BKR}
that a unital direct limit of approximately divisible \CA s is
approximately divisible. If $A$ is unital, it therefore follows that
$A$ is approximately divisible. Clearly $A$ is infinite, so it is
purely infinite by Theorem 1.4 (a) of \cite{BKR}.

If $A$ is  not unital, choose $k$ and an
infinite projection $p \in A_k$ such
that $\varphi_{k, \infty} (p) \neq 0$. Corollary 2.9 of \cite{BKR}
implies that
$\varphi_{k, l} (p) A_l \varphi_{k, l} (p)$
is approximately divisible for
$l \geq k$. Therefore
$\varphi_{k, \infty} (p) A \varphi_{k, \infty} (p) =
   \dirlim \varphi_{k, l} (p) A_l \varphi_{k, l} (p)$
is approximately divisible and infinite simple, hence purely infinite.
Now $A$ is stably isomorphic to
$\varphi_{k, \infty} (p) A \varphi_{k, \infty} (p)$ (see \cite{Bn1}),
and hence also purely infinite.

(2) We may assume that each $\varphi_{k, \infty}$ is injective. (If not,
we replace $A_k$ by $A_k / \ker(\varphi_{k, \infty})$, which is again a
finite direct sum of \pisca s.)

By Theorem 1.2 (i) of \cite{Zh1}, every purely infinite simple \CA\  has
real rank zero. It clearly follows from 3.1 of \cite{BP} that $A$ has
real rank zero. Therefore
it suffices to show that every nonzero projection $p\in A$ is
infinite. Choose $k$ and a projection $q\in A_k$  such that
$\varphi_{k,\infty}(q)$ is unitarily equivalent to $p.$
It suffices to prove that
that $\varphi_{k, \infty} (q)$ is infinite. But it is immediate that
$q$ is infinite in $A_k$, and infiniteness of $\varphi_{k, \infty} (q)$
now follows from injectivity of $\varphi_{k, \infty}$.  \QED

\vspace{0.6\baselineskip}

The following result is our main theorem, from which most of the other
results in this section will follow.

\vspace{0.6\baselineskip}

{\bf 5.4 Theorem } Let $A= \dirlim\,(A_n,\varphi_{n,n+1})$ and
$ B= \dirlim\,(B_n,\psi_{n,n+1})$ be two simple \CA s which are
direct limits of
even Cuntz-circle algebras, and assume that the maps
$\varphi_{n,\infty}$ and $\psi_{n,\infty}$ are all unital and
\aab.  If
$$
(K_0(A),[1_A],K_1(A))\cong (K_0(B),[1_B],K_1(B)),
$$
then $A\cong B.$

\vspace{0.6\baselineskip}

The proof consists of constructing an approximate intertwining of the
direct systems, as was first done in \cite{Ell2}. We will use the
particular statement given in \cite{Thn}.

If all the maps of the system were injective, this construction would
be fairly direct from Theorems 3.7 and 4.1. Lack of injectivity
causes some technical problems; in particular, we must replace the
system $(A_k, \varphi_{k, k+1})_{k = 1}^{\infty}$ by a system
 $(A_k', \varphi_{k, k+1}')_{k = 1}^{\infty}$ in which the maps
$\varphi_{k, \infty}'$ are approximately injective (see Definition 3.8),
and similarly with the system $(B_k, \psi_{k, k+1})_{k = 1}^{\infty}$.
We have to construct the $A_k'$ and $B_k'$ at the same time as the
approximate intertwining, which makes the proof somewhat complicated.
We will isolate the actual computational steps as the following two
lemmas; the proof of the theorem will then consist mainly of keeping
track of many maps and indices.

In the first of these lemmas, we refer to compact subsets of $S^1$ with
finitely many components. Note that such a subset is either $S^1$ itself
or a finite disjoint union of closed arcs and points.

\vspace{0.6\baselineskip}

{\bf 5.5 Lemma}
 Let $ A = \dirlim (A_k, \varphi_{k, k+1})$ be a simple
\DECSalg . Thus, in particular we can
write $A_1 = \bigoplus_{i=1}^{r} C(X_{i}) \otimes D_{i}$,
where the $D_{i}$ are matrix algebras over even Cuntz algebras
and the $X_{i}$ are compact subsets
of $S^1$ with finitely many components.
Assume
that each $\varphi_{k, k+1}$ is unital and that $\varphi_{1, \infty}$
is \aab .
Let $G \subset A_1$ be finite, and let $\ep > 0$.

Then there exist compact subsets $X_i'$ of $X_i$ with finitely many
components such that the
\CA\  $A_1' = \bigoplus_{i=1}^{r} C(X_{i}') \otimes D_{i}$ and the
restriction map $\alpha : A_1 \to A_1'$ satisfy the following.
There is $l > 1$ and  a \hm\  $\alpha': A_1' \to A_l$ such that
$\alpha' \circ \alpha = \varphi_{1,l}$, the map
$\varphi_{l, \infty} \circ \alpha'$ is $\ep$-approximately injective
with respect to the finite set $\alpha (G)$ (see Definition 3.8), and
$\varphi_{l, \infty} \circ \alpha'$ is \aab .

\vspace{0.6\baselineskip}

{\em Proof:}
For simplicity of notation, we will assume that
$A_1 = C(X, D)$, with $X \subset S^1$ and $D$ simple.
(Thus, we are assuming $A_1$ has only one summand. As will be
clear from the proof, if it has more, we will be able to
use the largest of the values of $l$ associated to
the summands.) Choose $\dt > 0$ such that whenever $x_1, x_2 \in X$
satisfy $ |x_1 - x_2| \leq \dt$ and $f \in G$, then
$\| f(x_1) - f(x_2) \| < \ep$. Let $u \otimes 1 \in C(X, D)$
be the usual standard generator of $C(S^1) \otimes {\bf C}$,
and note that
\[
{\rm sp}(\varphi_{1, \infty} (u \otimes 1)) =
   \bigcap_{l = 1}^{\infty} {\rm sp}(\varphi_{1, l} (u \otimes 1)).
\]
A standard compactness argument yields $l \geq 1$ such that
${\rm sp}(\varphi_{1, l} (u \otimes 1))$ is contained in a
$\dt /2$-neighborhood of ${\rm sp}(\varphi_{1, \infty} (u \otimes 1))$.
Furthermore, the complement
$X \setminus {\rm sp}(\varphi_{1, l} (u \otimes 1))$ is
a countable disjoint union of arcs open in $X$, so repeating the
compactness argument gives $W \subset X$, a finite union of arcs
open in $X$, such that $X' = X \setminus W$ is contained in a
$\dt /2$-neighborhood of ${\rm sp}(\varphi_{1, l} (u \otimes 1))$.

Define $A_1' = C(X', D)$, and let $\alpha$ be the
restriction map. Note that
\[
\ker(\alpha) = C_0 (X \setminus X', D) \subset
                C_0 (X \setminus {\rm sp}(\varphi_{1, l} (u \otimes 1)), D)
   = \ker(\varphi_{1, l}).
\]
Therefore $\varphi_{1, l}$ factors through $\alpha$; let
$\alpha' : A_1' \to A_l$ be the resulting map.
Then
\[
\ker(\varphi_{l, \infty} \circ \alpha' ) =
     C_0 (X' \setminus {\rm sp}(\varphi_{1, \infty} (u \otimes 1)), D),
\]
and $X'$ is contained in a $\dt$-neighborhood of
${\rm sp}(\varphi_{1, \infty} (u \otimes 1))$, so the choice of $\dt$
ensures that $\varphi_{l, \infty} \circ \alpha' $ is
$\ep$-approximately injective.

It remains to prove that $\varphi_{l, \infty} \circ \alpha' $ is \aab .
Since $\varphi_{1, \infty}$ is \aab , it is easy to check this directly
from the definition, using the relation
$(\varphi_{l, \infty} \circ \alpha') \circ \alpha = \varphi_{1, \infty}$
and the fact that $\alpha$ is a restriction map.  \QED

\vspace{0.6\baselineskip}

{\bf 5.6 Notation} Let $A$ and $B$ be \CA s, and let
$\varphi , \psi : A \to B$ be \hm s. If $G \subset A$ and $\ep > 0$,
then we write $\varphi \aeps{\ep} \psi$ (with respect to $G$) to mean
$\| \varphi(a) - \psi(a) \| < \ep$ for all $a \in G$. (Compare with
the relation  $\varphi \aueeps{\ep} \psi$ in Definition 1.1.)

\vspace{0.6\baselineskip}

{\bf 5.7 Lemma}
 Let $ A = \dirlim (A_k, \varphi_{k, k+1})$ be a simple
\DECSalg , with all $\varphi_{k, k+1}$ unital,
 let $D$ be an even \CSalg , and let
$\rho: A_1 \to D$ and $\theta_0 : D \to A_2$ be unital \hm s.
Let $0 < \ep \leq 1/2$, and let $F \subset A_1$ be a finite generating
set.  Assume that
\[
[\rho] \times [\theta_0] \times [\varphi_{2, \infty}]
   = [\varphi_{1, \infty}]
    \,\,\,\,\,\, {\rm in} \,\,\, KK^0 (A_1, A),
\]
that $\varphi_{2, \infty} \circ \theta_0 \circ \rho$
is injective and \aab ,
that $\theta_0$ is permanently \aab ,
and that $\varphi_{1, \infty}$ is \aab\  and $\ep$-approximately
injective with respect $F$. Finally, assume that, with respect to
some realization of $A_1$ as a finite direct sum as in the definition
of a \CSalg , the set $F$ contains the identities of all the summands.

Then there is $l \geq 2$ and a unital permanently
\aab\  \hm\  $\theta : D \to A_l$ such that
$[\theta] = [\theta_0] \times [\varphi_{2,l}]$ in $KK^0 (D, A_l)$,
such that
$[\rho] \times [\theta] = [\varphi_{1,l}]$ in $KK^0 (A_1, A_l)$,
and such that
$ \theta \circ \rho \aeps{4 \ep} \varphi_{1,l}$ with respect to $F$.

\vspace{0.6\baselineskip}

{\em Proof:}
Without loss of generality, assume the elements of $F$ all have norm at
most $1$.

Since $KK^0 (A_1, - )$ commutes with countable direct limits
(by Lemma 4.14), there is $m \geq 2$ such that
\[
[\rho] \times [\theta_0] \times [\varphi_{2, m}]
   = [\varphi_{1, m}]
    \,\,\,\,\,\, {\rm in} \,\,\, KK^0 (A_1, A_m).
\]
Also, the hypotheses and Theorem 3.11 imply that
$\varphi_{2, \infty} \circ \theta_0 \circ \rho \aueeps{2 \ep}
                                              \varphi_{1, \infty}$
with respect to $F$. (We have $A$ purely infinite by Lemma 5.3 (1).
Also, we apply Theorem 3.11 to each summand separately.
That is, if $A_1 = A_{11} \oplus \cdots \oplus A_{1r}$, and $e_i$ is the
identity of $A_{1i}$, then we first observe that
$(\varphi_{2, \infty} \circ \theta_0 \circ \rho) (e_i)$ is close to
$\varphi_{1, \infty} (e_i)$. Thus these projections are unitarily
equivalent for each $i$; without loss of generality, we assume they
are equal for each $i$. Now use Theorem 3.11 on the restrictions of
$\varphi_{2, \infty} \circ \theta_0 \circ \rho$ and
$\varphi_{1, \infty}$ to each $A_{1i}$, regarded as maps to
$\varphi_{1, \infty} (e_i) A \varphi_{1, \infty} (e_i)$.)

There is thus a unitary $v \in A$ such that
\[
\| v^* \varphi_{1, \infty} (a) v -
     (\varphi_{2, \infty} \circ \theta_0 \circ \rho ) (a) \| < 2 \ep
\]
for all $a \in F$.
Choose $m' \geq m$ such that there is a unitary $w \in A_{m'}$ with
$\| \varphi_{m', \infty} (w) - v \| < \ep /2$. Then for each $a \in F$,
we have
\[
\| \varphi_{m', \infty} \left( w^* \varphi_{1, m'} (a) w -
    (\varphi_{2, m'} \circ \theta_0 \circ \rho ) (a) \right) \| < 3 \ep.
\]
Since $F$ is finite, there is therefore $l \geq m'$ such that
\[
\| \varphi_{m', l} \left( w^* \varphi_{1, m'} (a) w -
     (\varphi_{2, m'} \circ \theta_0 \circ \rho ) (a) \right) \| < 4 \ep
\]
for all $a \in F$.

Let $z = \varphi_{m', l} (w)$, and   define
$\theta (b) = z (\varphi_{2, l} \circ \theta_0) (b) z^*$ for
$b \in D$.  Then
$ \theta \circ \rho \aeps{4 \ep} \varphi_{1,l}$ with respect to $F$.
 Note that $\theta$ is the composite of a unital \hm\  with
the permanently \aab\  \hm\  $\theta_0$, and hence still
permanently \aab\  by Lemma 1.14. Also, conjugation by a unitary does
not change the class in $KK$-theory. Therefore
$[\theta] = [\theta_0] \times [\varphi_{2, l}]$. Furthermore,
\[
[\rho] \times [\theta] =
[\rho] \times [\theta_0] \times [\varphi_{2, m}] \times [\varphi_{m, l}]
 = [\varphi_{1, l}]
\]
in $KK^0 (A_1, A)$, by the choice of $m$ at the beginning of the proof.
\QED

\vspace{0.6\baselineskip}

{\em Proof of Theorem 5.4:} Since $K_i(A)\cong K_i(B),$ the Universal
Coefficient Theorem (Theorem 1.17 of \cite{RS}) and Proposition
7.3 of \cite{RS} yield
an invertible $\sigma\in KK^0(A,B).$
(Lemma 4.3  implies that
both $A$ and $B$ are in the bootstrap category ${\cal N}$
of \cite{RS}.)
Let $\sigma^{-1}\in KK^0(B,A)$ be the inverse.

Fix realizations of each $A_k$ and $B_k$ as
direct sums as in the definition of a \CSalg . (Quotients of
$A_k$ and $B_k$ will then be realized as direct sums in the same way.)
 Let $(F_k)_{k = 1}^{\infty}$ and
$(G_k)_{k = 1}^{\infty}$ be systems of finite generating sets for
$(A_k, \varphi_{k, k+1})_{k = 1}^{\infty}$ and
$(B_k, \psi_{k, k+1})_{k = 1}^{\infty}$ respectively. We require that
$F_k$ contain the identities of the summands of $A_k$, and similarly
for $G_k$ and $B_k$.
We construct a diagram as follows:

\begin{picture}(330, 195)(-15, -25)

\put( 0,150){\makebox(0,0){$A_{m(1)}$}}
\put( 100,150){\makebox(0,0){$A_{m(2)}$}}
\put( 200,150){\makebox(0,0){$A_{m(3)}$}}
\put( 300,150){\makebox(0,0){$\cdots$}}
\put( 0,100){\makebox(0,0){$A_1'$}}
\put( 100,100){\makebox(0,0){$A_2'$}}
\put( 200,100){\makebox(0,0){$A_3'$}}
\put( 300,100){\makebox(0,0){$\cdots$}}
\put( 0,50){\makebox(0,0){$B_1'$}}
\put( 100,50){\makebox(0,0){$B_2'$}}
\put( 200,50){\makebox(0,0){$B_3'$}}
\put( 300,50){\makebox(0,0){$\cdots$}}
\put( 0,  0){\makebox(0,0){$B_{n(1)}$}}
\put( 100,  0){\makebox(0,0){$B_{n(2)}$}}
\put( 200,  0){\makebox(0,0){$B_{n(3)}$}}
\put( 300,  0){\makebox(0,0){$\cdots$}}

\put(15,150){\vector(1,0){70}}
\put(115,150){\vector(1,0){70}}
\put(215,150){\vector(1,0){75}}
\put(10,100){\vector(1,0){80}}
\put(110,100){\vector(1,0){80}}
\put(210,100){\vector(1,0){80}}
\put(10,50){\vector(1,0){80}}
\put(110,50){\vector(1,0){80}}
\put(210,50){\vector(1,0){80}}
\put(15,  0){\vector(1,0){70}}
\put(115,  0){\vector(1,0){70}}
\put(215,  0){\vector(1,0){75}}

\put(50, 156){\makebox(0,0)[b]{$\varphi_{{m(1)},{m(2)}}$}}
\put(150, 156){\makebox(0,0)[b]{$\varphi_{{m(2)},{m(3)}}$}}
\put(250, 156){\makebox(0,0)[b]{$\varphi_{{m(3)},{m(4)}}$}}
\put(50, 106){\makebox(0,0)[b]{$\varphi_{1,2}'$}}
\put(150,106 ){\makebox(0,0)[b]{$\varphi_{2,3}'$}}
\put(250,106 ){\makebox(0,0)[b]{$\varphi_{3,4}'$}}
\put(50,  56){\makebox(0,0)[b]{$\psi_{1,2}'$}}
\put(150,  56){\makebox(0,0)[b]{$\psi_{2,3}'$}}
\put(250,  56){\makebox(0,0)[b]{$\psi_{3,4}'$}}
\put(50, -2){\makebox(0,0)[t]{$\psi_{{n(1)},{n(2)}}$}}
\put(150, -2){\makebox(0,0)[t]{$\psi_{{n(2)},{n(3)}}$}}
\put(250, -2){\makebox(0,0)[t]{$\psi_{{n(3)},{n(4)}}$}}

\put(  0,142){\vector(0,-1){34}}
\put(100,142){\vector(0,-1){34}}
\put(200,142){\vector(0,-1){34}}
\put(  0, 92){\vector(0,-1){34}}
\put(100, 92){\vector(0,-1){34}}
\put(200, 92){\vector(0,-1){34}}
\put(  0,  8){\vector(0, 1){34}}
\put(100,  8){\vector(0, 1){34}}
\put(200,  8){\vector(0, 1){34}}

\put( -2,125){\makebox(0,0)[r]{$\alpha_1$}}
\put( -2, 75){\makebox(0,0)[r]{$\rho_1$}}
\put(  2, 25){\makebox(0,0)[l]{$\beta_1$}}
\put( 98,125){\makebox(0,0)[r]{$\alpha_2$}}
\put( 98, 75){\makebox(0,0)[r]{$\rho_2$}}
\put(102, 25){\makebox(0,0)[l]{$\beta_2$}}
\put(198,125){\makebox(0,0)[r]{$\alpha_3$}}
\put(198, 75){\makebox(0,0)[r]{$\rho_3$}}
\put(202, 25){\makebox(0,0)[l]{$\beta_3$}}

\put( 10, 105){\vector(2,1){80}}
\put(110, 105){\vector(2,1){80}}
\put(210, 105){\vector(2,1){80}}
\put( 10,  55){\vector(2,1){80}}
\put(110,  55){\vector(2,1){80}}
\put(210,  55){\vector(2,1){80}}
\put( 10, 45 ){\vector(2, -1){75}}
\put(110, 45 ){\vector(2, -1){75}}
\put(210, 45 ){\vector(2, -1){80}}

\put( 48, 127){\makebox(0,0)[br]{$\alpha_1'$}}
\put(148, 127){\makebox(0,0)[br]{$\alpha_2'$}}
\put(248, 127){\makebox(0,0)[br]{$\alpha_3'$}}
\put( 50,  77){\makebox(0,0)[br]{$\theta_1$}}
\put(150,  77){\makebox(0,0)[br]{$\theta_2$}}
\put(250,  77){\makebox(0,0)[br]{$\theta_3$}}
\put( 50, 27){\makebox(0,0)[bl]{$\beta_1'$}}
\put(150, 27){\makebox(0,0)[bl]{$\beta_2'$}}
\put(250, 27){\makebox(0,0)[bl]{$\beta_3'$}}

\end{picture}

\noindent
Here, all triangles in the top and bottom sections are supposed to
commute. In the middle section, we require that
$\theta_k \circ \rho_k \aeps{1/2^k} \varphi_{k, k+1}'$ with respect to
$\alpha_k (F_{m(k)}) \cup (\theta_{k-1} \circ \beta_{k-1})(G_{n(k-1)})$,
and
$ \rho_{k+1} \circ \theta_k \aeps{1/2^k} \psi_{k, k+1}'$
with respect to
$\beta_k (G_{n(k)}) \cup (\rho_{k} \circ \alpha_{k}) (F_{m(k)})$.
In order to make this happen, we will also construct systems
$(F_k')_{k = 1}^{\infty}$ and $(G_k')_{k = 1}^{\infty}$ of finite
generating sets for $(A_k', \varphi_{k, k+1}')_{k = 1}^{\infty}$ and
$(B_k', \psi_{k, k+1}')_{k = 1}^{\infty}$ respectively, such that
\[
\alpha_k (F_{m(k)}) \cup (\theta_{k-1} \circ \beta_{k-1})(G_{n(k-1)})
       \subset F_k' \,\,\,{\rm and}\,\,\,
\beta_k (G_{n(k)}) \cup (\rho_{k} \circ \alpha_{k}) (F_{m(k)})
       \subset G_k'.
\]
We will then require that
$\theta_k$ and $\rho_k$ be permanently \aab , that
$\varphi_{k, \infty}'$ be \aab\  and
$(1/2^{k+2})$-approximately injective
with respect to $F_k'$, and that
$\psi_{k, \infty}'$ be \aab\  and $(1/2^{k+2})$-approximately injective
with respect to $G_k'$.

We construct this diagram by induction on the column number. Some of the
steps will start by going further out in one of the original direct
systems, and so we will have to construct temporary versions of some of
the maps in the diagram. They will be distinguished from the final ones
with tildes.

We show the details only for columns 1 and 2; the remaining steps are
essentially the same as for column 2.

{\em Step 1, part A:} Set $m(1) = 1$. Use Lemma 5.5 to construct
$A_1'$, a surjective
\hm\  $\alpha_1 : A_{m(1)} \to A_1'$, an integer $i(1) > m(1)$,
and a \hm\  $\tilde{\alpha}_1' : A_1' \to A_{i(1)}$ such that
$\tilde{\alpha}_1' \circ \alpha_1 = \varphi_{1, i(1)}$,
and the map
$\varphi_{i(1), \infty} \circ \tilde{\alpha}_1'$ is
$1/8$-approximately injective
with respect to the finite set $\alpha_1 (F_1)$ and is \aab .
Define $F_1' = \alpha_1 (F_1)$. Since $\alpha_1$ is a direct sum of
restriction maps, and
 $F_1$ contains the identities of the summands of $A_1$, the set
 $F_1'$ contains the identities of the summands of $A_1'$.

{\em Step 1, part B:}
Use Theorem 4.1 to find $n(1)$ and a permanently \aab\  unital
\hm\  $\tilde{\rho}_1: A_1' \to B_{n(1)}$ such that
\[
[\tilde{\rho}_1] \times [\psi_{n(1), \infty}] =
  [\tilde{\alpha}_1'] \times [\varphi_{i(1), \infty}] \times \sigma
    \,\,\,\,\,\, {\rm in} \,\,\, KK^0 (A_1', B).
\]

Use Lemma 5.5 to find
$B_1'$,
a \hm\  $\beta_1 : B_{n(1)} \to B_1'$, an integer $j(1) > n(1)$,
and a \hm\  $\tilde{\beta}_1' : B_1' \to B_{j(1)}$ such that
$\tilde{\beta}_1' \circ \beta_1 = \psi_{n(1), j(1)}$,
and the map
$\varphi_{j(1), \infty} \circ \tilde{\beta}_1'$ is
$1/8$-approximately injective
with respect to the finite set
$\beta_1 (G_{n(1)} \cup \tilde{\rho}_1 (F_1'))$ and is \aab .

Define $G_1' = \beta_1 (G_{n(1)} \cup \tilde{\rho}_1 (F_1'))$.
Then $G_1'$ contains the identities of the summands of $B_1'$.
Further define $\rho_1 = \bt_1 \circ \tilde{\rho}_1$. Notice that
the surjectivity of $\bt_1$ implies that $\rho_1$ is still
permanently \aab . Furthermore, we have
$\tilde{\beta}_1' \circ \rho_1 = \psi_{n(1), j(1)} \circ
                                         \tilde{\rho}_1$,
whence
\[
[\rho_1] \times [\tilde{\bt}_1] \times [\psi_{j(1), \infty}] =
   [\tilde{\alpha}_1'] \times [\varphi_{i(1), \infty}] \times \sigma
    \,\,\,\,\,\, {\rm in} \,\,\, KK^0 (A_1', B).
\]

{\em Step 2, part A:} Use Theorem 4.1 to find $l \geq i(1)$,
where $i(1)$ is as chosen in Step 1 part A,
and a permanently \aab\  unital
\hm\  $\tilde{\theta}_1: B_1' \to A_{l}$ such that
\[
[\tilde{\theta}_1] \times [\varphi_{l, \infty}] =
  [\tilde{\beta}_1'] \times [\psi_{j(1), \infty}] \times \sigma^{-1}
    \,\,\,\,\,\, {\rm in} \,\,\, KK^0 (B_1', A).
\]

Now compute:
\beqr
\lefteqn{  [\tilde{\alpha}_1' ] \times [ \varphi_{i(1), l}]
                                     \times [ \varphi_{l, \infty}]
    = [\tilde{\alpha}_1' ] \times [\varphi_{i(1), \infty}] \times
                                  \sigma  \times \sigma^{-1}   } \\
 & & = [\rho_1 ] \times [\tilde{\bt}_1'] \times [\psi_{j(1), \infty}]
                                                \times \sigma^{-1} =
       [\rho_1 ] \times [\tilde{\theta}_1] \times [ \varphi_{l, \infty}]
\eeqr
in $KK^0 ( A_1', A)$. We next apply Lemma 5.7 with $D = B_1'$,
$\rho = \rho_1$, $\theta_0 = \tilde{\theta}_1$, and the direct
system being
\[
A_1'
   \stackrel{\varphi_{i(1), l} \circ \tilde{\alpha}_1'}{\longrightarrow}
A_{l} \rightarrow  A_{  l  + 1} \rightarrow A_{  l  + 2}
            \rightarrow \cdots.
\]
Note that
$\varphi_{l, \infty} \circ \varphi_{i(1), l} \circ \tilde{\alpha}_1'$ is
$1/8$-approximately injective with respect to $F_1'$, by Step 1
part A, and that
$\varphi_{l, \infty} \circ \tilde{\theta}_1 \circ \rho_1$
is injective and \aab\  by Lemma 1.14 (because $\rho_1$ is permanently
\aab\  and $\varphi_{l, \infty} \circ \tilde{\theta}_1$ is unital).
Therefore
Lemma 5.7 yields $m(2) \geq l$ and a unital permanently
\aab\  \hm\  $\tilde{\tilde{\theta}}_1 : B_1' \to A_{m(2)}$
such that
$[ \tilde{\tilde{\theta}}_1 ] =
          [\tilde{\theta}_1 ] \times [\varphi_{l, m(2)}]$ and
$[\rho_1 ] \times [ \tilde{\tilde{\theta}}_1] =
         [\tilde{\alpha}_1' ] \times [ \varphi_{i(1), m(2)}]$,
and also
$\tilde{\tilde{\theta}}_1 \circ \rho_1 \aeps{1/2}
         \varphi_{i(1), m(2)} \circ \tilde{\alpha}_1'$
with respect to $F_1'$.

Define $\alpha_1' = \varphi_{i(1), m(2)} \circ \tilde{\alpha}_1'$.

Use Lemma 5.5 as before to construct $A_2'$, a surjective
\hm\  $\alpha_2 : A_{m(2)} \to A_2'$, an integer $i(2) > m(2)$,
and a \hm\  $\tilde{\alpha}_2' : A_2' \to A_{i(2)}$ such that
$\tilde{\alpha}_2' \circ \alpha_2 = \varphi_{m(2), i(2)}$,
the map
$\varphi_{i(2), \infty} \circ \tilde{\alpha}_2'$ is
$1/16$-approximately injective with respect to
$\alpha_2 (F_{m(2)} \cup \tilde{\tilde{\theta}}_1 (G_1'))$,
and this map is \aab .
Define $\theta_1 = \alpha_2 \circ \tilde{\tilde{\theta}}_1$, which is
still permanently \aab\  since $\alpha_2$ is surjective.
Define $\varphi_{1,2}' = \alpha_2 \circ \alpha_1'$, and
define
$F_2' = \alpha_2 (F_{m(2)} \cup \tilde{\tilde{\theta}}_1 (F_1'))$.

Note that
$\theta_1 \circ \rho_1 = \alpha_2 \circ \tilde{\tilde{\theta}}_1
                                 \circ \rho_1$
and
$\varphi_{1,2}' =
      \alpha_2 \circ \varphi_{i(1), m(1)} \circ \tilde{\alpha}_1'$.
It follows that
$\theta_1 \circ \rho_1 \aeps{1/2} \varphi_{1,2}'$
and $[\rho_1] \times [\theta_1] = [\varphi_{1,2}']$.
Furthermore,
\beqr
[\theta_1] \times [\tilde{\alpha}_2'] \times [\varphi_{i(2), \infty}]
  & = & [\tilde{\tilde{\theta}}_1] \times [\varphi_{m(2), \infty}]  \\
  & = & [\tilde{\theta}_1] \times [\varphi_{l, \infty}]
    = [\tilde{\beta}_1'] \times [\psi_{j(1), \infty}] \times \sigma^{-1}
\eeqr
in $KK^0 (B_1', A)$. Finally,  $F_2'$ contains the identities of the
summands of $A_2'$.

{\em Step 2 part B:} This step is similar to part A, so we will be
briefer. Use Theorem 4.1 to find $l \geq j(1)$ and a permanently
\aab\  unital \hm\  $\tilde{\rho}_2 : A_2' \to B_l$ such that
$[\tilde{\rho}_2] \times [\psi_{l, \infty}] =
     [\tilde{\alpha}_2'] \times [\varphi_{i(2), \infty}] \times \sigma$.
Use Lemma 5.7 to find $n(2) \geq l$ and a unital permanently
\aab\  \hm\  $\tilde{\tilde{\rho}}_2 : A_2' \to B_{n(2)}$
such that
$[\tilde{\tilde{\rho}}_2] = [\tilde{\rho}_2] \times [\psi_{l, n(2)}]$
and
$[\theta_1] \times [\tilde{\tilde{\rho}}_2] =
    [\tilde{\bt}_1'] \times [\psi_{j(1), n(2)}]$,
and also
$\tilde{\tilde{\rho}}_2 \circ \theta_1  \aeps{1/2}
    \psi_{j(1), n(2)} \circ \tilde{\bt}_1' $
with respect to $G_1'$.
Define $ \bt_1' = \psi_{j(1), n(2)} \circ \tilde{\bt}_1' $.

Now use Lemma 5.5 to produce a surjective
\hm\  $\bt_2 : B_{n(2)} \to B_2'$
and a \hm\  $\tilde{\bt}_2' : B_2' \to B_{j(2)}$ (with $j(2) > n(2)$)
such that $\tilde{\bt}_2' \circ \bt_2  = \psi_{n(2), j(2)}$,
the map $\psi_{j(2), \infty} \circ  \tilde{\bt}_2' $
is $1/16$-approximately injective with respect to
$\bt_2 (G_{n(2)} \cup \tilde{\tilde{\rho}}_2 (F_2'))$,
and $\psi_{j(2), \infty} \circ \tilde{\bt}_2' $ is \aab .
Define $\rho_2 = \bt_2  \circ  \tilde{\tilde{\rho}}_2$ and
$\psi_{1,2}' = \bt_2 \circ \bt_1'$, and set
$G_2' = \bt_2 (G_{n(2)} \cup \tilde{\tilde{\rho}}_2 (F_2'))$.
Then check that $\rho_2 \circ \theta_1  \aeps{1/2}  \psi_{1,2}'$, that
$[\theta_1] \times [\rho_2] = [\psi_{1,2}']$, that
$[\rho_2] \times [\tilde{\bt}_2'] \times [\psi_{j(2), \infty}] =
     [\tilde{\alpha}_2'] \times [\varphi_{i(2), \infty}] \times \sigma$,
and that $G_2'$ contains the identities of the summands of $B_2'$.
This completes part B of Step 2.

Each successive half step consists of a repetition of the argument
already used in the two halves of Step 2, using, in order, Theorem
4.1, Lemma 5.7, and Lemma 5.5. In Step $k$ part A, we choose
$\varphi_{i(k), \infty} \circ \tilde{\alpha}_{k}'$ to be
$(1/2^{k + 2})$-approximately injective, and we get
$\theta_{k - 1} \circ \rho_{k - 1} \aeps{1/2^{k - 1}}
            \varphi_{k - 1, k}'$
with respect to $F_{k-1}'$. In part B, we choose
$\psi_{j(k), \infty} \circ \tilde{\beta}_{k}'$ to be
$(1/2^{k + 2})$-approximately injective, and we get
$\rho_k \circ \theta_{k - 1} \aeps{1/2^{k - 1}}
            \psi_{k - 1, k}'$
with respect to $G_{k-1}'$.

With this construction, the top and bottom rows of triangles in our
diagram commute. Therefore they induce isomorphisms
$\dirlim A_k \cong \dirlim A_k'$ and
$\dirlim B_k \cong \dirlim B_k'$
The middle row of triangles is an approximate
intertwining in the sense of Definition 2 of \cite{Thn}.
By Theorem 3 of \cite{Thn}, it therefore induces an isomorphism
$\dirlim A_k' \cong \dirlim B_k'$. So $A \cong B$. \QED

\vspace{0.6\baselineskip}

One easily checks that the isomorphism constructed in this proof induces
the same map on $K$-theory as $\sigma$. However, it is not clear that
its class in $KK^0 (A, B)$ is equal to $\sigma$.

A one sided version of the previous proof, using Lemma 1 of \cite{Thn},
establishes the following:

\vspace{0.6\baselineskip}

{\bf 5.8 Proposition } Let $A$ and $B$ as in Theorem 5.4. If there are
homomorphisms
$$
\alpha_0: K_0(A)\to K_0(B)\andeqn \alpha_1: K_1(A)\to K_1(B)
$$
such that $\alpha_0([1_A])=[1_B],$ then there is a unital
homomorphism $\varphi: A\to B$ such that
$$
\varphi_*^{(0)}=\alpha_0\andeqn \varphi_*^{(1)}=\alpha_1.\,\,\,  \QED
$$

\vspace{0.6\baselineskip}

We now give a number of corollaries of Theorem 5.4.

\vspace{0.6\baselineskip}

{\bf 5.9 Corollary } Let $A = \dirlim\,(A_k,\varphi_{k,k+1})$ and
$B = \dirlim\,(B_k,\psi_{k,k+1})$ be unital simple \CA s,
such that each $A_k$ and each $B_k$ is a finite direct
sum of matrix algebras over algebras $C(X)$, with each $X$ being a
point, a compact interval, or a circle. ($A$ and $B$ can independently
have real rank either 0 or 1.)
If $m$ is even and
\[
\left( K_0 (A \otimes \OA{m}), [1_{A \otimes \OA{m}}],
                 K_1 (A \otimes \OA{m}) \right)    \cong
\left( K_0 (B \otimes \OA{m}), [1_{B \otimes \OA{m}}],
                   K_1 (B \otimes \OA{m}) \right),
\]
then
$$
A\otimes \OA{m}\cong B\otimes \OA{m}.
$$

{\it Proof:}
We may assume that the $\varphi_{k,k+1}$ and $\psi_{k,k+1}$
are unital. Let $\Phi_{k,k+1}=\varphi_{k,k+1}\otimes {\rm id}_{\OA{m}}$
 and
$\Psi_{k,k+1}=\psi_{k,k+1}\otimes {\rm id}_{\OA{m}}$
Then
$$
A \otimes \OA{m} \cong \dirlim\,(A_k\otimes \OA{m}, \Phi_{k,k+1})
             \andeqn
B \otimes \OA{m} \cong \dirlim\,(B_k\otimes \OA{m}, \Psi_{k,k+1}).
$$
Both  $A\otimes \OA{m}$ and $B\otimes \OA{m}$ are simple, hence
(by Lemma 5.3 (1))
both are purely infinite simple.
It follows from Corollary 1.9 that
 $\Phi_{k,\infty}$ and $\Psi_{k,\infty}$ are \aab.
Therefore Theorem 5.4 applies.  \QED

\vspace{0.6\baselineskip}

{\bf 5.10 Corollary} Let $A$ and $B$ be simple direct limits of
circle algebras etc.\  as in Corollary 5.9, and let $m$ be even.
If
\[
(K_0(A),[1_A], K_1(A))\cong(K_0(B),[1_B],K_1(B))
\]
(ignoring order),
then
$$
A\otimes \OA{m}\cong B\otimes \OA{m}.
$$

\vspace{0.6\baselineskip}

{\em Proof:}
The K\"{u}nneth formula \cite{Sch0} and its
splitting (see Remark 7.11 of \cite{RS}) imply
that
\[
\left(K_0(A \otimes \OA{m}),[1_{A \otimes \OA{m}}],
                    K_1(A \otimes \OA{m}) \right) \cong
\left(K_0(B \otimes \OA{m}),[1_{B \otimes \OA{m}}],
                      K_1(B \otimes \OA{m}) \right).
\]
\QED

\vspace{0.6\baselineskip}

{\bf 5.11 Corollary } Let $A = \dirlim\,(A_k,\varphi_{k,k+1})$ and
$B = \dirlim\,(B_k,\psi_{k,k+1})$ be unital simple \CA s,
such that each $A_k$ and each $B_k$ is a finite direct
sum of matrix algebras over algebras $C(X)$, with each $X$ being a
compact subset of $S^1$. ($A$ and $B$ can independently
have real rank either 0 or 1.) If $m$ is even and
\[
\left( K_0 (A \otimes \OA{m}), [1_{A \otimes \OA{m}}],
                 K_1 (A \otimes \OA{m}) \right)    \cong
\left( K_0 (B \otimes \OA{m}), [1_{B \otimes \OA{m}}],
                   K_1 (B \otimes \OA{m}) \right),
\]
then
$$
A\otimes \OA{m}\cong B\otimes \OA{m}.
$$

\vspace{0.6\baselineskip}

{\it Proof:} We only need to show that $A$ and $B$ can be rewritten
as
direct limits as above, but with the restriction that the subsets
$X \subset S^1$ have only finitely many components.
Now $A$ and $B$ clearly satisfy condition (ii) of Theorem
4.3 of \cite{Ell2}. That theorem therefore implies they
are direct limits of finite direct sums of matrix algebras over
$C(S^1)$. \QED

\vspace{0.6\baselineskip}

{\bf 5.12 Corollary } Let $A_{\theta_1}$ and $A_{\theta_2}$ be two
irrational rotation algebras and let $m$ be even. Then
$$
A_{\theta_1}\otimes \OA{m}\cong A_{\theta_2}\otimes \OA{m}.
$$

\vspace{0.6\baselineskip}

{\it Proof:} It follows from \cite{EE} that every irrational
rotation algebra is a simple (unital) direct limit of finite direct
sums of matrix algebras over $C(S^1)$.
The Pimsner-Voiculescu exact sequence \cite{PV2}
shows that, ignoring the order,
 $K_0 (A_{\theta})\cong K_1(A_{\theta})\cong {\bf Z}\oplus{\bf Z}$
for all $\theta$.
\QED

\vspace{0.6\baselineskip}

By contrast, recall that by \cite{Rf} and \cite{PV}, we have
$A_{\theta_1}  \cong  A_{\theta_2}$ only when
$\theta_1 = \pm \theta_2 \pmod{{\bf Z}}$.

\vspace{0.6\baselineskip}

{\bf 5.13 Corollary} (Corollary 3.6 of \cite{Ln5})
 Let $A$ be a simple direct limit of circle algebras etc., as in
Corollary 5.9. Then $ A\otimes \OA{2}\cong \OA{2}.$ \QED

\vspace{0.6\baselineskip}

We now prove that the class $\Class$ is closed under several natural
operations. These results will enable us to extend the classification
results above to the nonunital case.

\vspace{0.6\baselineskip}

{\bf 5.14 Lemma} If $A\in {\Class},$ then $A\otimes
{\cal K}\in \Class.$

\vspace{0.6\baselineskip}

{\it Proof:} Suppose that $A = \dirlim\,(A_k,\varphi_{k,k+1}),$ where
the $A_k$ are even Cuntz-circle algebras and the $\varphi_{k,k+1}$ are
\aab. Define $\psi_{k,k+1}: M_k(A_k)\to M_{k+1}(A_{k+1})$ by
$\psi_{k,k+1} = {\rm id}_{M_k}\otimes \varphi_{k,k+1}\oplus 0.$
Then
\[
\dirlim\,(M_k(A_k),\psi_{k,k+1})\cong A\otimes {\cal K}.
\]
Note that $\psi_{k,\infty}$ maps $M_k(A_k)$ into $M_k(A).$
Corollary 1.10 implies that
 $\psi_{k,\infty}$ is \aab. Therefore $A\otimes {\cal K}\in
\Class$. \QED

\vspace{0.6\baselineskip}

{\bf 5.15 Lemma } Let $A$ be a simple \CA\  in $\Class$,  and
let $p\in A$ be a nonzero projection. Then $pAp\in \Class$.

\vspace{0.6\baselineskip}

{\bf Proof:} Write $A = \dirlim\,(A_k,\varphi_{k,k+1}),$ where each
$A_k$ is an even Cuntz-circle algebra and each
 $\varphi_{k,\infty}$ is
\aab. There is $l$ and a projection $q\in A_l$ such that
$ \|\varphi_{l,\infty}(q)-p\|<1.$ Therefore
$ \varphi_{l,\infty}(q)$ is unitarily equivalent to $p.$
It follows that
\[
pAp \cong
   \varphi_{l,\infty}(q) A \varphi_{l,\infty}(q)  \cong
   \mbox{$\dirlim$}_{k \geq l} \varphi_{l,k}(q) A_k \varphi_{l, k}(q).
\]
Lemma 1.11 implies that the algebras in this direct system are
even \CSalg s, and Lemma 1.12 implies that the maps from them to
the direct limit are \aab . Certainly $pAp$ is simple, so
$pAp \in \Class$. \QED

\vspace{0.6\baselineskip}

{\bf 5.16 Corollary } Let $A\in \Class$ and let $B$ be a
hereditary $C^*$-subalgebra of $A.$ Then $B\in \Class.$

\vspace{0.6\baselineskip}

{\it Proof:} If $B$ is unital, then $B=pAp$ for some projection
$p\in A.$ So $B \in \Class$ by Lemma 5.15. Otherwise, note
that $A$ is purely infinite by Lemma 5.3 (1). Therefore $B$ is stable
by Theorem 1.2 (i) of \cite{Zh1}.
Using \cite{Bn1},
it follows that $B \cong A \otimes {\cal K}$, whence
$B \in \Class$ by Lemma 5.14.  \QED

\vspace{0.6\baselineskip}

{\bf 5.17 Theorem } Let $A$ and $B$ be two simple \CA s in
$\Class$. Suppose that both $A$ and $B$ are nonunital and that
$$
(K_0 (A), K_1 (A)) \cong (K_0 (B), K_1 (B)).
$$
Then $A\cong B.$

\vspace{0.6\baselineskip}

{\bf Proof:} By Lemma 5.3 and Theorem 1.2 (i) of \cite{Zh1},
both $A$ and $B$ are stable.
Let $p\in A$ be a nonzero projection. Then \cite{Bn1} implies that
$A\cong pAp\otimes {\cal K}.$ Let $\alpha: K_0(A)\to K_0(B)$ be
an isomorphism. Then there is a nonzero projection $q\in B$ such that
$\alpha([p])=[q].$ By \cite{Bn1} again, $B\cong qBq\otimes {\cal K}.$
It follows from Lemma 5.15 that $pAp,\, qBq\in \Class$.
Theorem 5.4 therefore gives $ pAp\cong qBq.$
Consequently
$$
A\cong pAp\otimes {\cal K}\cong qBq\otimes {\cal K}\cong B. \,\,\, \QED
$$

\vspace{0.6\baselineskip}

{\bf 5.18 Theorem } Let $A = \dirlim\,(A_k, \alpha_{k, k+1}),$
where each $A_k$ is a finite direct sum of (simple) \CA s in
$\Class$.
Assume that $A$ is simple. Then $A \in \Class$.

\vspace{0.6\baselineskip}

{\it Proof:} We first note that $A$ is purely infinite by Lemma 5.3 (2).

We now reduce to the unital case. If $A$ is not unital, choose a nonzero
projection $p \in A$ such that $p =  \alpha_{l, \infty} (q)$ for
some $l$ and some projection $q \in A_l$. Then
$pAp \otimes {\cal K} \cong A$, as in the proof of Corollary
5.16. So it suffices to show that $pAp \in \Class$.
Now
$pAp \cong
   \mbox{$\dirlim$}_{k \geq l}  \alpha_{l, k} (q) A  \alpha_{l, k} (q)$,
and Lemma 5.15 implies that each algebra in this direct system is a
finite direct sum of algebras in $\Class$.   Thus, we
may assume that $A$ is unital. Therefore we may assume that all
the maps $ \alpha_{k, k+1}$ are unital too.

For each $k$, write
$A_k = \bigoplus_{j = 1}^{r(k)} A_{k}^{(j)}$,
with each $A_{k}^{(j)} \in \Class$.
(To keep the notation in this proof straight, we will write indices
associated with direct sums as superscripts.)
Note that each
$ \alpha_{k, \infty}|_{A_{k}^{(j)}}$
is either injective or zero, and that
we can drop all $A_{k}^{(j)}$ for which this map is zero. Thus, without
loss of generality, we can assume that each $ \alpha_{k, \infty} $ is
injective; then so is each $ \alpha_{k, k + 1}$.
Let $\pi_{k}^{(j)} : A_k \to A_{k}^{(j)}$ be the projection map.
Now choose finite sets $F_k \subset A_k$ such that
$F_k = \bigcup_{j = 1}^{r(k)} F_{k}^{(j)}$
with $F_{k}^{(j)} \subset A_{k}^{(j)}$,
such that
$ \alpha_{k, k+1} (F_k) \subset F_{k + 1}$ for each $k$, and
such that
$\bigcup_{k = 1}^{\infty}  \alpha_{k, \infty} (F_k)$ is dense in the
unit ball of $A$.

We will construct even \CSalg s
$D_k = \bigoplus_{j = 1}^{r(k)} D_{k}^{(j)}$, finite generating subsets
$H_{k}^{(j)}$ contained in the unit ball of $D_{k}^{(j)}$ and
with $1_{D_{k}^{(j)}} \in H_{k}^{(j)}$,
and \hm s
$\rho_k = \bigoplus_{j = 1}^{r(k)} \rho_{k}^{(j)} : D_k \to A_k$ and
$\dt_{k, k+1} : D_k \to D_{k + 1}$, as in the following
approximately commutative diagram:

\begin{picture}(330, 85)(-15, 40)

\put( 0,100){\makebox(0,0){$D_1$}}
\put( 100,100){\makebox(0,0){$D_2$}}
\put( 200,100){\makebox(0,0){$D_3$}}
\put( 300,100){\makebox(0,0){$\cdots$}}
\put( 0,50){\makebox(0,0){$A_1$}}
\put( 100,50){\makebox(0,0){$A_2$}}
\put( 200,50){\makebox(0,0){$A_3$}}
\put( 300,50){\makebox(0,0){$\cdots$}}

\put(10,100){\vector(1,0){80}}
\put(110,100){\vector(1,0){80}}
\put(210,100){\vector(1,0){80}}
\put(10,50){\vector(1,0){80}}
\put(110,50){\vector(1,0){80}}
\put(210,50){\vector(1,0){80}}

\put(50, 106){\makebox(0,0)[b]{$\dt_{1,2}$}}
\put(150,106 ){\makebox(0,0)[b]{$\dt_{2,3}$}}
\put(250,106 ){\makebox(0,0)[b]{$\dt_{3,4}$}}
\put(50,  56){\makebox(0,0)[b]{$ \alpha_{1,2}$}}
\put(150,  56){\makebox(0,0)[b]{$ \alpha_{2,3}$}}
\put(250,  56){\makebox(0,0)[b]{$ \alpha_{3,4}$}}

\put(  0, 92){\vector(0,-1){34}}
\put(100, 92){\vector(0,-1){34}}
\put(200, 92){\vector(0,-1){34}}

\put( -2, 75){\makebox(0,0)[r]{$\rho_1$}}
\put( 98, 75){\makebox(0,0)[r]{$\rho_2$}}
\put(198, 75){\makebox(0,0)[r]{$\rho_3$}}

\end{picture}

We will require:

(1) The $k$-th square approximately commutes up to $1/2^k$, that is,
$\rho_{k + 1} \circ \dt_{k, k+1} \aeps{1/2^k}
              \alpha_{k, k+1} \circ \rho_k$
with respect to $H_k = \bigcup_{j = 1}^{r(k)} H_{k}^{(j)}$.

(2) The squares commute in $KK$-theory, that is,
\[
[\rho_{k + 1} \circ \dt_{k, k+1}] = [ \alpha_{k, k+1} \circ \rho_k]
          \,\,\,\,\,\, {\rm in}  \,\,\, KK^0 (D_k, A_{k + 1}).
\]

(3) Each $\rho_{k}^{(j)}$ is injective and \aab .

(4) $\dt_{k, k+1}$ is permanently \aab .

(5) $\dt_{k, k+1} (H_k) \subset H_{k + 1}$ for all $k$.

(6) Every point of $F_k$ is within $1/2^{k - 1}$ of some point of
$\rho_{k} (H_k)$.

Suppose the diagram has been constructed. Let
$D = \dirlim (D_k, \dt_{k, k+1})$. Lemma 1 of \cite{Thn} gives
a \hm\  $\rho : D \to A$ such that for every $k$ and every $a \in D_k$,
we have
\[
(\rho \circ \dt_{k, \infty}) (a) =
         \lim_{l \to \infty}
                ( \alpha_{l, \infty} \circ \rho_l \circ \dt_{k, l}) (a).
\]
Since each $\alpha_{l, \infty} \circ \rho_l \circ \dt_{k, l}$
is injective, this implies $\rho$ is isometric,
hence injective. One furthermore checks that for $a \in H_k$, we
have
\[
\| (\rho \circ \dt_{k, \infty}) (a) -
               ( \alpha_{k, \infty} \circ \rho_k) (a) \|
\leq \sum_{l = k}^{\infty} 1/2^l = 1/2^{k - 1}.
\]
It follows from (6) that every point of $ \alpha_{k, \infty} (F_k)$
is within $1/2^{k - 2}$ of a point of $\rho (D)$. Since
$ \alpha_{k, k+1} (F_k) \subset F_{k + 1}$, and the union of the images
of these sets in $A$ in dense in the unit ball of $A$, we conclude that
$\rho$ is surjective. Thus $D \cong A$. In particular, $D$ is purely
infinite and simple.  Condition (4) now implies that each
$\dt_{k, \infty}$ is \aab . It follows that $D$, and hence $A$, is in
$\Class$.

We construct the squares in the diagram one at a time, using induction.
We start with the construction of $\rho_1$.
Since each $A_{1}^{(j)}$ is a unital algebra in $\Class$, we can write
$A_{1}^{(j)} = \dirlim (C_{l}^{(j)}, \gamma_{l, l + 1}^{(j)} )$,
where each
$C_{l}^{(j)}$ is an even Cuntz-circle algebra and
$\gamma_{l, \infty}^{(j)} $ is \aab .
Choose $l$ so large that for each $j$ there is a subset $G^{(j)}$ of the
unit ball of $C_l^{(j)}$ whose image in $A_{1}^{(j)}$ approximates each
point of
$F_{1}^{(j)}$ to within $1/4$.
Increasing the size of $G^{(j)}$, we may assume
it contains the identities of the summands and generates $C_l^{(j)}$.
Use Lemma 5.5 to produce \CSalg s
$C_l^{(j) \prime}$, obtained as quotients
of $C_l^{(j)}$ with quotient maps $\kappa^{(j)}$, and unital \hm s
$\gamma_{l, \infty}^{(j) \prime} : C_l^{(j) \prime} \to A_{1}^{(j)}$
which are
$1/8$-approximately injective with respect to the image
$G^{(j) \prime}$
of $G^{(j)}$ in $C_l^{(j) \prime}$, which are \aab , and such that
$\gamma_{l, \infty}^{(j) \prime} =
               \gamma_{l, \infty}^{(j) \prime} \circ \kappa^{(j)}$.
Using Lemma 3.9, find injective \aab\  unital \hm s
$\gamma_{l, \infty}^{(j) \prime \prime} :
                        C_l^{(j) \prime} \to A_{1}^{(j)}$
such that
$\gamma_{l, \infty}^{(j) \prime \prime} \aeps{1/2}
                               \gamma_{l, \infty}^{(j) \prime}$
with respect to $G^{(j) \prime}$. Set
$D_{1}^{(j)} = C_l^{(j) \prime}$, $H_{1}^{(j)} =  G^{(j) \prime}$ and
$\rho_{1}^{(j)} = \gamma_{l, \infty}^{(j) \prime \prime}$.
Then set $D_1 = \bigoplus_{j = 1}^{r(1)} D_{1}^{(j)}$ and
$H_1 = \bigcup_{j = 1}^{r(1)} H_{1}^{(j)}$.
Note that, with these
definitions, $\rho (H_1)$ generates $D_1$ and approximates each point
of $F_1$ to within $1$.

We now assume that $D_k$ and $\rho_k$ have been constructed, and we
construct $D_{k + 1}$, $\rho_{k + 1}$, and $\dt_{k, k+1}$.
Since $A_{k + 1}^{(j)}$ is a unital algebra in $\Class$, we can write
$A_{k + 1}^{(j)} = \dirlim (C_{l}^{(j)}, \gamma_{l, l+1}^{(j)})$,
where each
$C_{l}^{(j)}$ is an even Cuntz-circle algebra.
There is $l$ such that there are, for each $j$,
\mops\  $q^{(i,j)} \in C_{l}^{(j)}$ with
$\gamma_{l, \infty}^{(j)} ( q^{(i,j)})$ close enough to
$(\pi_{k + 1}^{(j)} \circ  \alpha_{k, k+1} ) (1_{A_{k}^{(i)}})$
that there is a unitary
$v \in A_{k + 1}$ with $\| v - 1 \| < \ep_1$ and
$v [ \gamma_{l, \infty}^{(j)} ( q^{(i,j)}) ] v^* =
          (\pi_{k + 1}^{(j)} \circ  \alpha_{k, k+1} ) (1_{A_{k}^{(i)}})$
for all $i$ and $j$. (The number $\ep_1 > 0$ will be specified later.)
Using Theorem 4.1, find $l' \geq l$ and \hm s
$\mu^{(i, j)} : D_{k}^{(i)} \to C_{l'}^{(j)}$ with
$\mu^{(i, j)} (1) = q^{(i,j)}$ and
\[
[\mu^{(i, j)}] \times [\gamma_{l', \infty}^{(j)}] =
     [\rho_{k}^{(i)}] \times [ \alpha_{k, k+1}] \times [\pi_{k+1}^{(j)}]
        \,\,\,\,\,\, {\rm in} \,\,\, KK^0 (D_{k}^{(i)}, A_{k+1}^{(j)}).
\]
We can require that $\mu^{(i, j)}$ be permanently \aab\  whenever
$q^{(i,j)} \neq 0$.

For $q^{(i, j)} \neq 0$, we now apply an argument similar to, but
simpler than, the proof of Lemma 5.7, to the two \hm s
$\mu^{(i, j)} : D_k^{(i)} \to q^{(i, j)} C_{l'}^{(j)} q^{(i, j)}$
and
\[
v^* [(\pi_{k+1}^{(j)} \circ \alpha_{k,k+1} \circ \rho_{k}^{(i)}) (-)] v:
D_k^{(i)} \to
[ \gamma_{l', \infty}^{(j)} (q^{(i, j)}) ] A_{k + 1}^{(j)}
                         [\gamma_{l', \infty}^{(j)} (q^{(i, j)})],
\]
noting that
\[
[\gamma_{l', \infty}^{(j)} (q^{(i, j)})] A_{k + 1}^{(j)}
                 [\gamma_{l', \infty}^{(j)} (q^{(i, j)})] =
\mbox{$\dirlim$}_m
[\gamma_{l', m}^{(j)} (q^{(i, j)})] C_{l'}^{(j)}
                              [\gamma_{l', m}^{(j)} (q^{(i, j)})]
\]
is a \DECSalg . We obtain a number $l'' \geq l'$ (obtained as the
maximum over $i$ and $j$ of suitable numbers depending on $i$ and
$j$) and  permanently \aab\  \hm s
\[
\mu^{(i, j) \prime} :
D_k^{(i)} \to
[\gamma_{l', l''}^{(j)} (q^{(i, j)})] C_{l''}^{(j)}
                         [\gamma_{l', l''}^{(j)} (q^{(i, j)})]
\]
such that
\[
\gamma_{l'', \infty}^{(j)} \circ \mu^{(i, j) \prime}  \aeps{\ep_2}
 v^*[(\pi_{k+1}^{(j)} \circ  \alpha_{k, k+1} \circ \rho_{k}^{(i)}) (-)]v
\]
with respect to $H_k^{(i)}$ for all $i$ and $j$. (The number $\ep_2 > 0$
will be specified below.)

The rest of the induction step is essentially the same as the
construction of $\rho_1$ in the initial step. Choose $l''' \geq l''$
and finite generating subsets $G^{(j)}$
of the unit ball of $C_{l'''}^{(j)}$
containing the identities of the summands and
$\bigcup_i (\gamma_{l'', l'''}^{(j)} \circ \mu^{(i, j) \prime} )
                                                       (H_k^{(i)})$,
and whose image in $A_{k + 1}^{(j)}$ approximates every element
of $F_{k + 1}^{(j)}$ to within $\ep_3$. Use Lemma 5.5 to produce
suitable quotients $D_{k + 1}^{(j)}$ of $C_{l'''}^{(j)}$, with quotient
maps $\kappa^{(j)}$ and $\ep_4$-approximately injective
\aab\  \hm s $\gamma_{l''', \infty}^{(j) \prime}$ to $A_{k + 1}^{(j)}$.
The approximate injectivity is with respect to the image
$H_{k + 1}^{(j)}$ of $G^{(j)}$. Then use Lemma 3.9 to replace these
\hm s by injective \aab\  \hm s
$\gamma_{l''', \infty}^{(j) \prime \prime}$
which agree to within $2 \ep_4$ on $H_{k + 1}^{(j)}$.
Define
$\rho_{k + 1}^{(j)} (a) =
              v \gamma_{l''', \infty}^{(j) \prime \prime} (a) v^*$
and
\[
\dt_{k, k + 1} = \bigoplus_j \sum_i
  \kappa^{(j)} \circ \gamma_{l'', l'''}^{(j)} \circ \mu^{(i, j) \prime}.
\]
Then the square in the diagram commutes to within
$2 \ep_1 + \ep_2 + 2 \ep_4$ on $H_k$, and the image of $H_{k + 1}$
approximates every element of $F_k$ to within
$2 \ep_1 + \ep_3 + 2 \ep_4$. So we choose the $\ep_m$ such that these
numbers are both less than $1/2^{k}$. \QED

\vspace{0.6\baselineskip}

{\bf 5.19 Corollary } Let $A\in \Class$ and let $B$ be a
(separable) AF algebra.  Then $A\otimes B \in \Class$.

\vspace{0.6\baselineskip}

{\it Proof:} It is trivial that $A\otimes M_n\in  \Class$
for any $n.$  Apply Theorem 5.18.  \QED

\vspace{0.6\baselineskip}

We now use Theorem 5.18 to give a classification theorem for a class
of direct limits in which no restrictions are imposed on the
maps of the direct systems. The building blocks will be matrix
algebras over even Cuntz algebras, and the even algebras from the
following definition, which do for $K_1$ what the the Cuntz algebras
do for $K_0$.

\vspace{0.6\baselineskip}

{\bf 5.20 Definition} Let $D$ be the Bunce-Deddens algebra
\cite{BD}  whose
ordered $K_0$-group is ${\bf Q}$.
(See, for example, 10.11.4 of \cite{Bl2}.)
For $2 \leq m < \infty$, we define the {\em co-Cuntz algebra}
${\cal Q}_m$ to be $D \otimes \OA{m}$.

\vspace{0.6\baselineskip}

{\bf 5.21 Lemma} We have
\[
K_0 (\QA{m}) = 0 \andeqn K_1 (\QA{m}) \cong {\bf Z} / (m - 1){\bf Z}.
\]
If $m$ is even, then $\QA{m} \in \Class$. Moreover, if $A$ is any unital
\CA\  in $\Class$ with $K_0 (A) = 0$ and
$K_1 (A) \cong {\bf Z} / (m - 1){\bf Z}$, then $A \cong \QA{m}$.

\vspace{0.6\baselineskip}

{\em Proof:} Let $D$ be as in Definition 5.20. Note that
$K_1 (D) \cong {\bf Z}$ (see 10.11.4 of \cite{Bl2}),
and that if $G$ is any torsion group,
then ${\bf Q} \otimes G = \Tor_1^{\bf Z} ({\bf Q}, G) = 0$.
The computation of
$K_* (\QA{m})$ now follows from the K\"{u}nneth formula \cite{Sch0}.

It is well known (Theorem 2 of \cite{BD}) that $D$ is a simple
direct limit of \CA s of the form $C(S^1) \otimes M_n$. If $m$ is
even, then  $\QA{m} \in \Class$ as in the proof of Corollary 5.9.
The last sentence follows from Theorem 5.4. \QED

\vspace{0.6\baselineskip}

{\bf 5.22 Notation} Let $\Class_0$ denote the class of all simple
direct limits $A = \dirlim (A_k, \varphi_{k, k + 1})$ in which each
$A_k$ is a finite direct sum of finite matrix algebras over even
Cuntz algebras and even co-Cuntz algebras.

\vspace{0.6\baselineskip}

{\it Note that no conditions are imposed on the maps in the system.}

\vspace{0.6\baselineskip}

{\bf 5.23 Proposition} $\Class_0 \subset \Class$.

\vspace{0.6\baselineskip}

{\em Proof:} This is immediate from Theorem 5.18 and Lemma 5.21. \QED

\vspace{0.6\baselineskip}

We therefore get the following classification theorem for $\Class_0$.

\vspace{0.6\baselineskip}

{\bf 5.24 Theorem}
(1) Let $A, B \in \Class_0$ be unital, and suppose that
\[
\left( K_0 (A), [1_{A}], K_1 (A) \right)    \cong
\left( K_0 (B), [1_{B}], K_1 (B) \right).
\]
Then $ A \cong B$.

(2) Let $A, B \in \Class_0$ be nonunital, and suppose that
\[
\left( K_0 (A),  K_1 (A) \right) \cong \left( K_0 (B),  K_1 (B) \right).
\]
Then $ A \cong B$. \QED

\vspace{0.6\baselineskip}

It is clear that for every algebra $A \in \Class$, the groups
$K_0 (A)$ and $K_1 (A)$ are countable torsion groups in which
every element has odd order. We now prove a converse. Our result will
also show that $\Class_0 = \Class$. We need a lemma first.

\vspace{0.6\baselineskip}

{\bf 5.25 Lemma} Let
\[
B = \bigoplus_{i = 1}^m B_i \andeqn C = \bigoplus_{j = 1}^n C_j
\]
be finite direct sums of matrix algebras over even Cuntz algebras
and  even co-Cuntz algebras. Let $\lambda: B \to C$ be a unital
\hm . Then there is a unital \hm\  $\varphi : B \to C$ such that
$\varphi_* = \lambda_*$ as maps from $K_* (B)$ to $K_* (C)$, and
such that every partial map $\varphi_{i, j}: B_i \to C_j$ is nonzero.
(The maps $\varphi_{i,j}$ are defined to be
$\kappa_j \circ \varphi \circ \mu_i$, where $\kappa_j : C \to C_j$ is
the quotient map and $\mu_i : B_i \to B$ is the inclusion.)

\vspace{0.6\baselineskip}

{\em Proof:}
It suffices to prove this when $n = 1$, that is, when $C$ is simple.
Let $v \in C$ satisfy $v^* v = 1$ and $q = v v^* < 1$. Then
$[1 - q] = 0$ in $K_0 (C)$. Since $C$ is purely infinite simple,
there are nonzero
\mops\  $p_1, \dots, p_m \in C$ with $\sum_i p_i = 1-q$
and each $[p_i] = 0$ in $K_0 (C)$.
It suffices to construct unital \hm s $\varphi_i : B_i \to p_i C p_i$
with $[\varphi_i] = 0$ in $KK^0 (B_i, C)$. We will then define,
for $a = (a_1, \dots, a_m) \in B$,
\[
\varphi (a) = v \lambda (a) v^* + \sum_{i = 1}^m \psi_i (a_i).
\]

If $C$ is a
matrix algebra over an even Cuntz algebra,
then so is $p_i C p_i$ (by Lemma 1.11). If $C$ is a
matrix algebra over an even co-Cuntz algebra,
then $p_i C p_i$ is a co-Cuntz algebra (since
it is in $\Class$ and has the right $K$-theory).
Also, $B_i$ is given as a matrix algebra over something,
but
an argument as in the
beginning of the proof of Theorem 4.13 enables us to assume that
the matrix size is $1$.
So we need to prove that if $B$ is an
even Cuntz algebra or co-Cuntz algebra,
and $C = M_k (C_0)$  where $C_0$ is an
even Cuntz algebra or co-Cuntz algebra,
with $[1_C] = 0$ in $K_0 (C)$, then there is an (injective)
unital \hm\  $\varphi$ from $B$ to $C$. Moreover, if $C_0$ is a
co-Cuntz algebra, we can vary $k$ at will.

If both $B$ and $C_0$ are even Cuntz algebras, use Lemma 1.15. If
both are co-Cuntz algebras, let $D$ be as in Definition 5.20,
take $C = D \otimes M_{n - 1} (\OA{n})$, and take $\varphi$
to be the tensor product of ${\rm id}_D$ with a suitable map
$\OA{m} \to M_{n - 1} (\OA{n})$ from Lemma 1.15.
Now suppose $B = \QA{m}$ and $C_0 = \OA{n}$.
Choose
$\alpha : \OA{m} \to \OA{2}$ and
$\gamma : \OA{2} \to M_k (\OA{n})$ as in Lemma 1.15, and use
Corollary 5.13 to choose an isomorphism
$\bt : D \otimes \OA{2} \to \OA{2}$. Then let $\varphi$ be
the composite
\[
\QA{m} \stackrel{\id \otimes \alpha}{\longrightarrow}
   D \otimes \OA{2} \stackrel{\bt}{\longrightarrow}
   \OA{2} \stackrel{\gamma}{\longrightarrow}  M_k (\OA{n}).
\]
Finally, suppose $B = \OA{m}$ and $C = M_k (\QA{n})$. We may assume
$k = n - 1$. Let $\alpha$, $\beta$, and $\gamma$ be as above,
and set $\varphi = (\id_D \otimes \gamma) \circ \bt^{-1} \circ \alpha$.
\QED

\vspace{0.6\baselineskip}

{\bf 5.26 Theorem} Let $G_0$ and $G_1$ be countable abelian
torsion groups in which element has odd order,
and let $g \in G_0$.
Then:

(1) There is a unital algebra $A \in \Class_0$ such that
\[
\left( K_0 (A), [1_{A}], K_1 (A) \right) \cong (G_0, g, G_1).
\]

(2) There is a nonunital algebra $A \in \Class_0$ such that
\[
\left( K_0 (A), K_1 (A) \right) \cong (G_0, G_1).
\]

\vspace{0.6\baselineskip}

{\em Proof:}
Part (2) follows from part (1) by tensoring with
the compact operators. (Note that the proof of Lemma 5.14 shows
equally well that $\Class_0$ is closed under tensoring with ${\cal K}$.)
So it suffices to prove (1).

Theorem 2.6 of \cite{Rr1} gives simple direct limits
$B^{(i)} \cong \dirlim (B^{(i)}_k, \psi_{k, k+1}^{(i)})$, in
which each $\psi_{k, k+1}^{(i)}$ is unital and each $B^{(i)}_k$
is a finite direct sum of even Cuntz algebras, such that
$K_0 (B^{(i)}) \cong G_i$ and $[1_{B^{(0)}}] = g$.
Let $D$ be as in Definition 5.20. Then define
$A_k = B_k^{(0)} \oplus (D \otimes B_k^{(1)})$, which is a direct sum
of matrix algebras over even Cuntz algebras and even co-Cuntz algebras.
Use Lemma 5.25 to find a unital
\hm\  $\varphi_{k, k + 1} : A_k \to A_{k + 1}$
which does the same thing on $K$-theory as
$\psi_{k, k + 1}^{(0)} \oplus
                    ({\rm id}_D \otimes \psi_{k, k + 1}^{(1)})$,
and such that all the partial maps between simple summands are nonzero.
Set $A = \dirlim (A_k, \varphi_{k, k + 1})$. Then $A$ is easily seen to
have the right $K$-theory and $[1] = g$ in $K_0 (A) \cong G_0$.
It is easy to check, using the condition on $\varphi_{k, k + 1}$,
that the algebraic direct limit of the
$A_k$ is simple. It follows from a standard argument that $A$ is
simple. So $A \in \Class_0$.  \QED

\vspace{0.6\baselineskip}

{\bf 5.27 Corollary} $\Class_0 = \Class$.

\vspace{0.6\baselineskip}

{\em Proof:} It follows from Theorem 5.26 that every possible value
of our invariant for $\Class$ is already attained for some algebra
in $\Class_0$.  \QED

\noindent
Huaxin Lin\\
Department of Mathematics\\
University of Oregon\\
Eugene, Oregon 97403-1222\\
and\\
Department of Mathematics\\
East China Normal University\\
Shanghai 200062, China\\

\noindent
N. Christopher Phillips\\
Department of Mathematics\\
University of Oregon\\
Eugene, Oregon 97403-1222

\end{document}